\documentclass{amsart}
\usepackage[all]{xy}
\usepackage{latexsym}
\usepackage{amssymb}
\usepackage{amsmath}
\usepackage{amsthm}
\usepackage{amscd}
\usepackage{fancyhdr}
\usepackage[dvips]{graphicx}
\usepackage{a4wide}
\usepackage{mathrsfs}
\input cyracc.def

 \newfam\cyrfam

  \font\tencyr=wncyr10

  \font\sevencyr=wncyr7

  \font\fivecyr=wncyr5

  \def\cyr{\fam\cyrfam\sevencyr\cyracc}

  \textfont\cyrfam=\tencyr \scriptfont\cyrfam=\sevencyr

    \scriptscriptfont\cyrfam=\fivecyr

  \newfam\cyifam

  \font\tencyi=wncyi10

  \font\sevencyi=wncyi7

  \font\fivecyi=wncyi5

  \def\cyi{\fam\cyifam\sevencyi\cyracc}

  \textfont\cyifam=\tencyi \scriptfont\cyifam=\sevencyi

    \scriptscriptfont\cyifam=\fivecyi

\def\id{{\mbox{1 \hskip -8pt 1}}}
 \newcommand{\lon}{\longrightarrow}
 
 \newcommand{\rar}{\rightarrow}
 \newcommand{\hook}{\hookrightarrow}
 \newcommand{\Proof}{{\bf Proof}.\, }

\newcommand{\p}{{\partial}}
\newcommand{\Id}{{\mathrm I\mathrm d}}
\newcommand{\no}{{\noindent}}

\newcommand{\Q}{{\mathbb Q}}
 \newcommand{\Z}{{\mathbb Z}}
 \newcommand{\bS}{{\mathbb S}}

 \renewcommand{\P}{{\mathbb P}}
 
 \newcommand{\R}{{\mathbb R}}
 
 \newcommand{\K}{{\mathbb K}}
 \newcommand{\ot}{\otimes}

\newcommand{\ULB}{{\cU \caL \mathit i \mathit e^1\hspace{-0.3mm} \mathit \cB}}
\newcommand{\Po}{\cU\cP\hspace{-0.5mm}\mathit o \mathit i \mathit s \mathit s
\mathit o \mathit n  }
  \newcommand{\Poly}{{\mathcal P}{\mathit o}{\mathit l}{\mathit y}}
  \newcommand{\CoLie}{{\mathit C \mathit o \mathit L \mathit i\mathit e}}

\newcommand{\eff}{{\mathit e\mathit f\mathit f}}
\newcommand{\LB}{{\mathcal L} ie^1\hspace{-0.3mm} {\mathcal B}}

\newcommand{\Ber}{{\mathit B}{\mathit e} {\mathit r}}
\newcommand{\Koz}{{\mbox {\scriptsize !`}}}
 \newcommand{\Beq}{\begin{equation}}
 \newcommand{\Eeq}{\end{equation}}
 \newcommand{\Beqr}{\begin{eqnarray}}
 \newcommand{\Eeqr}{\end{eqnarray}}
 \newcommand{\Beqrn}{\begin{eqnarray*}}
 \newcommand{\Eeqrn}{\end{eqnarray*}}
 \newcommand{\Ba}{\begin{array}}
 \newcommand{\Ea}{\end{array}}
 \newcommand{\Bi}{\begin{itemize}}
 \newcommand{\Ei}{\end{itemize}}
 \newcommand{\Bc}{\begin{center}}
 \newcommand{\Ec}{\end{center}}

 \newcommand{\fG}{{\mathfrak G}}

 \newcommand{\f}{{\mathcal O}}
 \newcommand{\cA}{{\mathcal A}}
 \newcommand{\cB}{{\mathcal B}}
 \newcommand{\cC}{{\mathcal C}}
 \newcommand{\caD}{{\mathcal D}}
 \newcommand{\cE}{{\mathcal E}}
 \newcommand{\cF}{{\mathcal F}}

 \newcommand{\cI}{{\mathcal I}}
 
 \newcommand{\caL}{{\mathcal L}}
 \newcommand{\cM}{{\mathcal M}}
 \newcommand{\cN}{{\mathcal N}}
 \newcommand{\cP}{{\mathcal P}}
 \newcommand{\cR}{{\mathcal R}}
 \newcommand{\cS}{{\mathcal S}}
 \newcommand{\cT}{{\mathcal T}}

 \newcommand{\cU}{{\mathcal U}}
 \newcommand{\cV}{{\mathcal V}}

 \newcommand{\al}{\alpha}
 \newcommand{\be}{\beta}

 \newcommand{\Ga}{\Gamma}
 \newcommand{\bGa}{{\mathbf \Gamma}}
 
 \newcommand{\var}{\varepsilon}
 \newcommand{\la}{\lambda}
 \newcommand{\om}{\omega}

 \newcommand{\tc}{{\mathbf c}}
 \newcommand{\ta}{{\mathbf a}}
 \newcommand{\tb}{{\mathbf b}}

 \newcommand{\zal}{{\bar{\al}}}

 \newcommand{\Hom}{{\mathrm H\mathrm o\mathrm m}}

 \def\hgw{{\mbox {\large $\circlearrowright$}}}
 \newcommand{\sip}{\smallskip}
 \newcommand{\bip}{\bigskip}
 \newcommand{\mip}{\vspace{2.5mm}}

\theoremstyle{plain}
\swapnumbers

\newtheorem{prop-def}[theorem]{Proposition-definition}

\newtheorem{f-theorem}{Formality Theorem}[section]
\newtheorem{main-theorem}{Main~Theorem}[section]
\newtheorem{section-theorem}{Theorem}[section]
\newtheorem{section-corollary}{Corollary}[section]

\theoremstyle{definition}

\newtheorem{fact-me}{Fact \cite{Me1}}[subsection]

\begin{document}

 \sloppy

\long\def\symbolfootnote[#1]#2{\begingroup%
\def\thefootnote{\fnsymbol{footnote}}\footnote[#1]{#2}\endgroup}

 \title{Wheeled pro(p)file of Batalin-Vilkovisky formalism
}
 \author{ S.A.\ Merkulov}
\address{Sergei~A.~Merkulov: Department of Mathematics, Stockholm University, 10691 Stockholm, Sweden}
\email{sm@math.su.se}
 \date{}
 \begin{abstract}
Using technique of wheeled props we establish a correspondence between the homotopy theory
of unimodular Lie 1-bialgebras and the famous Batalin-Vilkovisky formalism.
Solutions of the so called quantum master equation satisfying certain boundary conditions
are proven to be in 1-1 correspondence with representations of a wheeled dg prop which,
on the one hand, is isomorphic to the cobar construction of the prop of unimodular
Lie 1-bialgebras and, on the other hand, is  quasi-isomorphic to the dg wheeled prop of unimodular
Poisson structures. These results allow us to apply properadic methods for computing formulae
for a homotopy transfer of a unimodular Lie 1-bialgebra structure on an arbitrary complex to
the associated quantum master function
on its cohomology. It is proven that in the category of quantum BV manifolds associated with the homotopy
theory of unimodular  Lie 1-bialgebras quasi-isomorphisms are equivalence relations.

It is shown that   Losev-Mnev's  BF theory for unimodular Lie algebras can be naturally
extended to the case of unimodular Lie 1-bialgebras (and, eventually, to the case of unimodular Poisson structures).
Using a finite-dimensional version of the Batalin-Vilkovisky quantization formalism
it is rigorously proven  that the Feynman integrals computing
the effective action of this new BF theory describe precisely homotopy transfer formulae obtained
within the wheeled properadic approach to the quantum master equation. Quantum corrections
(which are present
 in our BF model to all orders
of the Planck constant) correspond precisely to what are often
called ``higher Massey products" in the homological algebra.

\end{abstract}
 \maketitle

\vspace{-5mm}

\begin{center}
\sc Contents
\end{center}
{\Small
1. {\bf Introduction}\hfill 2
\vspace{1mm}

2. {\bf Quadratic wheeled properads and homotopy transfer formulae}\hfill 4
\vspace{-0.1mm}

2.1. Wheeled operads, properads and props\hfill 4
\vspace{-0.1mm}

2.2. Morphisms of wheeled props\hfill 7
\vspace{-0.1mm}

2.3. Coprop(erad)s\hfill 7
\vspace{-0.1mm}

2.4. Bar construction\hfill 8
\vspace{-0.1mm}

2.5. Cobar construction\hfill 10
\vspace{-0.1mm}

2.6. Quadratic wheeled (co)properads and  Koszul duality\hfill 10
\vspace{-0.1mm}

2.7. Homotopy transfer formulae\hfill 12\ \

\vspace{1mm}

3. {\bf Geometry of quantum Batalin-Vilkovisky manifolds}\hfill 14
\vspace{-0.1mm}

 3.1. $\Z$-graded formal manifolds\hfill 14
\vspace{-0.1mm}

3.2. Odd Poisson structures \hfill 16
\vspace{-0.1mm}

3.3. Polyvector fields\hfill 16
\vspace{-0.1mm}

3.4. Odd symplectic structures\hfill 17
\vspace{-0.1mm}

3.5. Densities and semidensities\hfill 19
\vspace{-0.1mm}

3.6. Batalin-Vilkovisky manifolds\hfill 20
\vspace{-0.1mm}

3.7. Sheaves of Gerstenhaber-Batalin-Vilkovisky algebras\hfill 21
\vspace{-0.1mm}

3.8. Quantum master equation\hfill 22
\vspace{-0.1mm}

3.9. Quantum BV manifiolds\hfill 23
\vspace{-0.1mm}

3.10. Homotopy classification of quantum BV manifolds\hfill 24
\vspace{-0.1mm}

3.11. Quantum morphisms of BV manifolds\hfill 29\ \

\vspace{1mm}

4. {\bf From unimodular Lie 1-bialgebras to  quantum BV manifolds}\hfill 31
\vspace{-0.1mm}

4.1. Lie $n$-bialgebras\hfill 31
\vspace{-0.1mm}

4.2. Wheeled prop, $\ULB$, of unimodular Lie 1-bialgebras\hfill 33
\vspace{-0.1mm}

4.3. Representations of $\ULB_\infty$ $\Leftrightarrow$ quantum BV manifolds\hfill 34\
\vspace{-0.1mm}


\vspace{1mm}

5.{\bf  Wheeled prop of Poisson structures with vanishing modular class}\hfill 35
\vspace{-0.1mm}

5.1 Modular volume form\hfill 35
\vspace{-0.1mm}

5.2. Wheeled dg prop of unimodular Poisson structures\hfill 36
\vspace{-0.1mm}

5.3 Quasi-isomorphism theorem\hfill 36\ \

\vspace{1mm}

6. {\bf $BF$  theory of quantum BV manifolds}\hfill 37
\vspace{-0.1mm}

6.1. Introduction\hfill 37
\vspace{-0.1mm}

6.2. $BF$-theory of unimodular Lie 1-bialgebras\hfill 37
\vspace{-0.1mm}

\vspace{-0.1mm}
}

\section{Introduction}
The theory of operads and props has grown nowadays from a useful technical tool
into a kind of universal mathematical
language with the help of which topologist, algebraists, homotopy theorists  and  geometers can
fruitfully communicate
with each other. For example, one and the same operad of little 2-disks (i) solves
the recognition problem for based 2-loop spaces in algebraic topology, (ii)
describes homotopy Gerstenhaber
structure on the Hochschild deformation complex in homological algebra, and (iii)
controls diffeomorphism invariant Hertling-Manin's integrability equations \cite{HM}
in differential
geometry.
It is yet to see
whether or not basic concepts and constructions of theoretical physics can be
understood and developed in
the framework of operads and props, but the fact that {space-time},
``the  background of everything",  can be turned
into an ordinary observable --- a certain function (representation) on a prop --- is rather
intriguing.

\mip

This paper attempts to tell a story of the famous theoretical physics
{\em quantum master equation},
\Beq\label{introduction-master-eqn}
\hbar \Delta \Gamma + \frac{1}{2}\{\Gamma, \Gamma\}=0,
\Eeq
in the language of  wheeled prop(erad)s. It is shown that an important class of its
solutions (specified by certain boundary conditions in the quasi-classical limit)
is controlled by a surprisingly simple
wheeled prop of unimodular Lie 1-bialgebras and hence can be understood  as a class of strongly homotopy
algebras. It is proven  that the homotopy classification of this class of quantum master functions
is as simple as, for example,
the homotopy classification of strongly homotopy Lie algebras given in \cite{Ko}. These results allow us
to compare
the standard Feynman technique of producing new quantum master functions (called often in physics literature
``effective actions") by integrating the original ones
along certain Lagrangian submanifolds with the purely properadic homotopy transfer method
which uses Koszul duality theory,  and conclude (in a mathematically rigorous way) that they are identical to each other.

 \mip

Here is  a detailed description of paper's content.
Section 2 gives a self-contained introduction into the theory of wheeled props, their bar and cobar
constructions
\cite{Me-graphs,MMS}. We
introduce and study {\em Koszul duality}\, theory  for quadratic wheeled {properads}\footnote{Koszul
duality
for wheeled {\em operads}\, has been studied earlier in  \cite{MMS}.} having in mind
 applications (in \S 4 and 5) of the Koszul duality technique to two
 important for us examples, the first  of which
 controls the local finite-dimensional Poisson geometry,
and the other one the local geometry of master equation (\ref{introduction-master-eqn}).
The content of this theory is standard (cf.\ \cite{GJ}):
\Bi
\item
For any quadratic wheeled properad $\cP$ there is a naturally
 associated Koszul dual wheeled coproperad
$\cP^\bot$ which comes together with a canonical
 monomorphism
of dg coproperads,
$\imath: \cP^{\bot}\rar B(\cP)$, into the bar construction on $\cP$.
\item The cobar construction, $B^c(\cP^\bot)$, is a dg free wheeled
properad denoted in this paper by $\cP_\infty$.
\item
There exists an epimorphism, $\cP_\infty \rar \cP$,
which is a quasi-isomorphism if $\cP$ is Koszul.
\Ei
 The main result in \S 2 is Theorem 2.7.1
which, if reformulated shortly, says that
{\em given an arbitrary (not necessarily  Koszul)
 quadratic wheeled properad and an arbitrary dg $\cP$-algebra $V$, then every cohomological
 splitting of $V$
 makes canonically its cohomology, $H(V)$, into a $\cP_\infty$-algebra; moreover,
 this induced $\cP_\infty$ structure is given precisely by that sum of decorated graphs which describe
 the image of the canonical monomorphism $\imath: \cP^{\bot}\rar B(\cP)$}.
 This result gives
    a conceptual explanation of the well-known ``experimental" fact
that the homotopy transfer formulae of infinity structures can be given in terms of graphs.
A closely related result (for ordinary operads)  has been  obtained recently in \cite{CL}.
The first explicit graphic formulae have been obtained By Kontsevich and Soibelman
 \cite{KS} who have rewritten in terms of graphs
the homotopy transfer formulae of \cite{Me0} for the case when $\cP$ is an operad of
associative algebras. Another example can be found in the work of Mnev \cite{Mn} who
treated the case when $\cP$ is a wheeled operad of unimodular Lie algebras.
One more example of explicit transfer formulae (related to the master equation (\ref{introduction-master-eqn}))
is given below in \S 6.

\sip


\mip

In \S 3 we introduce and study a category, $\cC at(BV)$, of {\em  (quasi-classically split) quantum BV manifolds}\,
  whose
\Bi
\item {\em objects}, $\cM$, are, roughly speaking, formal solutions of all possible
quantum master equations (\ref{introduction-master-eqn})
with non-degenerate odd Poisson  brackets $\{\ , \}$ which satisfy in the quasiclassical
($\lim_{\hbar\rar 0} + \lim_{\hbar\rar 0}\frac{d}{d\hbar}$) limit certain  boundary
conditions (see \S 3.9 for a precise definition);
these boundary conditions imply that the tangent space, $\cT_*\cM$, to the formal  manifold $\cM$ at
 the distinguished point comes equipped with an induced differential $d$;
 if this induced differential vanishes, then
 $\cM$ is called {\em minimal}; if, on the other hand, $d$ encodes the full information about the corresponding
 solution to (\ref{introduction-master-eqn}) and the complex $(\cT_*\cM, d)$ is acyclic, then such a quantum BV manifold
 $\cM$ is called {\em contractible};

\item {\em morphisms} are generated by symplectomorphisms, natural projections $\cM_1\times \cM_2\rar \cM_1$,
and quantum embeddings, $\cM_1\rar \cM_1\times \cM_2$, depending on a choice of a Lagrangian submanifold in $\cM_2$.
\Ei
One has the following two results in the category $\cC at(BV)$:

(i) {\em
Every quantum BV manifold is isomorphic to the product of a minimal quantum BV manifold
and a contractible one}.

(ii) {\em
Quasi-isomorphisms are equivalence relations}.


\sip

In \S 4 the material of  \S 2 and \S 3 is tied  together. We introduce and study a wheeled
prop, $\ULB$, of unimodular Lie 1-bialgebras
and prove that {\em there is a one-to-one correspondence between
quantum BV manifolds  and representations of the associated
dg free wheeled prop $\ULB_\infty$}.
We do not know at present whether or not
the wheeled prop(erad) $\ULB$ is Koszul, i.e.\ whether or not the natural epimorphism,
$$
(\ULB_\infty, \delta) \lon (\ULB, 0),
$$
is a quasi-isomorphism.
  If it is, then
the wheeled prop quantization machine of \cite{Me-lec}
would apply to deformation quantization  of {\em unimodular}\,
Poisson structures.

\mip

Formal unimodular Poisson structures can be identified with a subclass of
solutions, $\Gamma$, of the master equation (\ref{introduction-master-eqn})
which are independent of $\hbar$. Hence there is a canonical epimorphism of dg wheeled
props,
$$
F: \ULB_\infty \lon \Po,
$$
where $\Po$ is a dg prop whose representations in a vector space $V$ are formal
unimodular Poisson structures on $V$ vanishing at $0$.
It is proven in \S 5 that $F$ is a quasi-isomorphism.

\mip

Section 6 is inspired by the work of Mnev  \cite{Mn} on a remarkable approach
to the homotopy transfer formulae of unimodular
$L_\infty$-algebras which is based on  the BV quantization of an extended $BF$ theory
and the associated  Feynman integrals. We apply in \S 6
 Losev-Mnev's ideas
to unimodular Lie 1-bialgebras
and show that the Feynman integrals technique provides us with exactly the same
formulae for the homotopy transfer of $\ULB_\infty$-structures as the ones which follow
from the Koszul duality theory for quadratic wheeled properads developed in  \S 2.
 These result implies essentially that the Ward identities
in a certain class of quantum field theories can be interpreted as equations for a {\em morphism}\,
of certain dg wheeled
(co)props.

\mip

A few words about notations.
The symbol $\bS_n$ stands for the permutation group, that is the group of all bijections,
$[n]\rar [n]$, where $[n]$ denotes (here and everywhere) the set
$\{1,2,\ldots,n\}$. If $V=\oplus_{i\in \Z} V^i$ is a graded vector space, then
$V[k]$ is a graded vector space with $V[k]^i:=V^{i+k}$.
We work throughout over a field $\K$ of characteristic 0 so that, for an action
of  finite group $G$  on a vector space $V$, the subspace of invariants, $\{v\in V| \sigma(v)=v\ \forall \sigma\in G\}$, is canonically isomorphic to the
quotient space of coinvariants, $V/span\{v-\sigma(v)\}_{v\in V, \sigma \in G}$, so that
we denote them by one and the same symbol $V_G$.

\bip

\section{Quadratic wheeled properads and homotopy transfer formulae}

{\bf 2.1. Wheeled operads, properads and props \cite{Me-graphs, MMS, Me-Perm}.}
Let $\fG^\circlearrowright$ be the family of all possible (not necessarily connected) graphs
constructed
 from the so called directed
$(m,n)$-{\em corollas},
\Beq\label{corolla}
\begin{xy}
 <0mm,0mm>*{\circ};
 <-0.5mm,0.2mm>*{};<-8mm,3mm>*{}**@{-},
 <-0.4mm,0.3mm>*{};<-4.5mm,3mm>*{}**@{-},
 <0mm,0mm>*{};<0mm,2.6mm>*{\ldots}**@{},
 <0.4mm,0.3mm>*{};<4.5mm,3mm>*{}**@{-},
 <0.5mm,0.2mm>*{};<8mm,3mm>*{}**@{-},
<-0.4mm,-0.2mm>*{};<-8mm,-3mm>*{}**@{-},
 <-0.5mm,-0.3mm>*{};<-4.5mm,-3mm>*{}**@{-},
 <0mm,0mm>*{};<0mm,-2.6mm>*{\ldots}**@{},
 <0.5mm,-0.3mm>*{};<4.5mm,-3mm>*{}**@{-},
 <0.4mm,-0.2mm>*{};<8mm,-3mm>*{}**@{-};
<0mm,5mm>*{\overbrace{\ \ \ \ \ \ \ \ \ \ \ \ \ \  }};
<0mm,-5mm>*{\underbrace{\ \ \ \ \ \ \ \ \ \ \ \ \ \ }};
<0mm,7mm>*{^{m\ \ output\ legs}};
<0mm,-7mm>*{_{n\ \ input\ legs}};
 \end{xy}, \ \ \ m,n\geq 0,
\Eeq
by taking their disjoint unions, and gluing some output legs in such a union with the same number of input legs.
The glued legs are called the {\em internal edges of the graph}\, and all the rest retain their name {\em legs of the graph}.
Note that every
internal edge as well as every leg of a graph is naturally directed; unless otherwise is explicitly shown,
we tacitly assume in all our pictures that the direction flow runs from the bottom to the top.
 We have
$\fG^\circlearrowright=\coprod_{m,n\geq 0} \fG^\circlearrowright(m,n)$, where $\fG^\circlearrowright(m,n)
\subset \fG^\circlearrowright$ is the subset of graphs having $m$
output legs and $n$ input legs. We assume from now on that the input legs of each graph $G\in  \fG^\circlearrowright(m,n)$
are labeled by the natural numbers $\{1,\ldots, n\}$ and the output legs are labeled by  $\{1,\ldots, m\}$  so that each
set $\fG^\circlearrowright(m,n)$ comes equipped with a natural action of the group $\bS_m\times \bS_n$.
For example,
$$
\Ba{c}
\begin{xy}
 <0mm,0mm>*{\circ};
 <0.39mm,-0.39mm>*{};<2.4mm,-2.4mm>*{}**@{-},
 <-0.35mm,-0.35mm>*{};<-1.9mm,-1.9mm>*{}**@{-},
 <-2.4mm,-2.4mm>*{\circ};
 <-2.0mm,-2.8mm>*{};<-0.4mm,-4.5mm>*{}**@{-},
 <2.4mm,-2.4mm>*{};<0.4mm,-4.5mm>*{}**@{-},
  <0mm,-5.1mm>*{\circ};
  <0.4mm,-5.5mm>*{};<2mm,-7.7mm>*{}**@{-},
<-0.4mm,-5.5mm>*{};<-2mm,-7.7mm>*{}**@{-},
  <0.4mm,-5.5mm>*{};<2.9mm,-9.7mm>*{^2}**@{},
 <0.4mm,-5.5mm>*{};<-2.9mm,-9.7mm>*{^1}**@{},
%
%
(0,0)*{}
   \ar@{->}@(ul,dl) (-2.4,-2.4)*{}
 \end{xy}
\Ea \in  \fG^\circlearrowright(0,2), \ \
\begin{xy}
 <0mm,-1.3mm>*{};<0mm,-3.5mm>*{}**@{-},
 <0.38mm,-0.2mm>*{};<2.2mm,2.2mm>*{}**@{-},
 <-0.38mm,-0.2mm>*{};<-2.2mm,2.2mm>*{}**@{-},
<0mm,-0.8mm>*{\circ};
 <2.4mm,2.4mm>*{\circ};
 <2.5mm,2.3mm>*{};<4.4mm,-0.8mm>*{}**@{-},
 <2.4mm,2.8mm>*{};<2.4mm,5.2mm>*{}**@{-},
     <0mm,-1.3mm>*{};<0mm,-5.3mm>*{^1}**@{},
     <2.5mm,2.3mm>*{};<5.1mm,-2.6mm>*{^2}**@{},
    <2.4mm,2.5mm>*{};<2.4mm,5.7mm>*{^1}**@{},
    <-0.38mm,-0.2mm>*{};<-2.8mm,2.5mm>*{^2}**@{},
    \end{xy}\in  \fG^\circlearrowright(2,2),
\ \
\Ba{c}
\begin{xy}
 <0mm,0mm>*{\circ};
 <0.39mm,-0.39mm>*{};<2.4mm,-2.4mm>*{}**@{-},
 <-0.35mm,-0.35mm>*{};<-1.9mm,-1.9mm>*{}**@{-},
 <-2.4mm,-2.4mm>*{\circ};
 <-2.0mm,-2.8mm>*{};<-0.4mm,-4.5mm>*{}**@{-},
 <2.4mm,-2.4mm>*{};<0.4mm,-4.5mm>*{}**@{-},
  <0mm,-5.1mm>*{\circ};
  <0.4mm,-5.5mm>*{};<2mm,-7.7mm>*{}**@{-},
<-0.4mm,-5.5mm>*{};<-2mm,-7.7mm>*{}**@{-},
  <0.4mm,-5.5mm>*{};<2.9mm,-9.7mm>*{^2}**@{},
 <0.4mm,-5.5mm>*{};<-2.9mm,-9.7mm>*{^3}**@{},
(0,0)*{}
   \ar@{->}@(ul,dl) (-2.4,-2.4)*{}
 \end{xy}
\Ea
\begin{xy}
 <0mm,-1.3mm>*{};<0mm,-3.5mm>*{}**@{-},
 <0.38mm,-0.2mm>*{};<2.2mm,2.2mm>*{}**@{-},
 <-0.38mm,-0.2mm>*{};<-2.2mm,2.2mm>*{}**@{-},
<0mm,-0.8mm>*{\circ};
 <2.4mm,2.4mm>*{\circ};
 <2.5mm,2.3mm>*{};<4.4mm,-0.8mm>*{}**@{-},
 <2.4mm,2.8mm>*{};<2.4mm,5.2mm>*{}**@{-},
     <0mm,-1.3mm>*{};<0mm,-5.3mm>*{^1}**@{},
     <2.5mm,2.3mm>*{};<5.1mm,-2.7mm>*{^4}**@{},
    <2.4mm,2.5mm>*{};<2.4mm,5.7mm>*{^2}**@{},
    <-0.38mm,-0.2mm>*{};<-2.8mm,2.5mm>*{^1}**@{},
    \end{xy}\in  \fG^\circlearrowright(2,4).
$$
Let $E$ be an $\bS$-{\em bimodule}, that is a family, $\{E(p,q)\}_{p,q\geq 0}$, of vector spaces on which the group
$\bS_p$ act on the left and the group $\bS_q$ act on the right and both actions commute with each other.
Fix an arbitrary graph $G\in \fG^\circlearrowright(m,n)$ and denote by $V(G)$ the set of its {\em vertices}, that is
the set of its
generating corollas (\ref{corolla}).
For each   $v\in V(G)$, denote by  ${\mathit I\mathit n}_v$ (resp.\ ${\mathit O\mathit u\mathit t}_v$)
the set of the input (resp.\ output) legs of the vertex $v$. Assume the cardinality of $In_v$ (resp.\ $Out_v$)
equals $q$ (resp.\ $p$) and note that vector spaces,
$$
\langle In_v\rangle :=\left\{
\Ba{ll}
 \mbox{the $q!$-dimensional vector space spanned by all bijections
$In_v\rar [q]$}
      & \mbox{if}\ q\geq 1\\
\K &  \mbox{if}\ q=0.
\Ea
\right.
$$
and
$$
\langle Out_v\rangle :=\left\{
\Ba{ll}
\hspace{-1mm} \mbox{the $p!$-dimensional vector space spanned by all bijections
$[p]\rar Out_v$}
      & \mbox{if}\ p\geq 1\\
\K &  \mbox{if}\ p=0.
\Ea
\right.
$$
have, respectively,  a natural left $\bS_q$-module structure and a right
  $\bS_p$-module structure. Hence one can form a quotient,
$$
E(Out_v, In_v):=  \langle Out_v\rangle \ot_{\bS_p} E(p,q) \ot_{\bS_q} \langle In_v\rangle,
$$
which is (non-canonically) isomorphic to $E(p,q)$ as a vector space and which carries  natural actions
of the automorphism groups of the sets $Out_v$ and $In_v$. These actions  make a so called
{\em unordered tensor product}\, over the set $V(G)$ (of cardinality, say, $N$),
$$
\bigotimes_{v\in V(G)} E(Out_v, In_v):= \left(\bigoplus_{i:\{1,...\, ,N\}\rar V(G) }
 E(Out_{i(1)}, In_{i(1)})
\ot\ldots \ot
 E(Out_{i(N)}, In_{i(N)})\right)_{\bS_N},
$$
into a representation space of the automorphism group, $Aut(G)$, of the graph $G$ which is,
by definition, the subgroup of the
 symmetry group of the 1-dimensional
$CW$-complex underlying the graph $G$ which fixes its legs. Hence with an arbitrary
graph $G\in \fG^\circlearrowright$
and an arbitrary $\bS$-bimodule $E$ one can associate a vector space,
$$
G\langle E\rangle:= \left(
\otimes_{v\in V(G)} E(Out_v, In_v)\right)_{Aut G},
$$
whose elements are called {\em decorated (by $E$) graphs}. For example, the automorphism
group of the graph
$G_0=\Ba{c}
\begin{xy}
 <0mm,0mm>*{\circ};
<0mm,0.41mm>*{};<0mm,2.9mm>*{}**@{-},
 <0.39mm,-0.39mm>*{};<2.4mm,-2.4mm>*{}**@{-},
 <-0.35mm,-0.35mm>*{};<-2.4mm,-2.4mm>*{}**@{-},
 <-2.4mm,-2.4mm>*{};<-0.4mm,-4.5mm>*{}**@{-},
 <2.4mm,-2.4mm>*{};<0.4mm,-4.5mm>*{}**@{-},
  <0mm,-5.1mm>*{\circ};
  <0.4mm,-5.5mm>*{};<2mm,-7.7mm>*{}**@{-},
<-0.4mm,-5.5mm>*{};<-2mm,-7.7mm>*{}**@{-},
  <0.4mm,-5.5mm>*{};<2.9mm,-9.7mm>*{^2}**@{},
 <0.4mm,-5.5mm>*{};<-2.9mm,-9.7mm>*{^1}**@{},
 \end{xy}
\Ea$
is $\Z_2$ so that $G_0\langle E \rangle= E(1,2)\ot_{\Z_2} E(2,2)$. It is useful to think
of an element in $G_0\langle E\rangle$  as  the graph $G_0$
whose vertices are literarily decorated by some
elements $a\in E(1,2)$ and $b\in E(2,1)$;  this pictorial representation of
$G_0\langle E\rangle$ is correct provided  the
relations,
$$
\begin{xy}
 <0mm,0mm>*{\circ};
<3mm,0mm>*{a};
<0mm,0.41mm>*{};<0mm,2.9mm>*{}**@{-},
 <0.39mm,-0.39mm>*{};<2.4mm,-2.4mm>*{}**@{-},
 <-0.35mm,-0.35mm>*{};<-2.4mm,-2.4mm>*{}**@{-},
 <-2.4mm,-2.4mm>*{};<-0.4mm,-4.5mm>*{}**@{-},
 <2.4mm,-2.4mm>*{};<0.4mm,-4.5mm>*{}**@{-},
  <0mm,-5.1mm>*{\circ};
 <3mm,-5.1mm>*{b};
  <0.4mm,-5.5mm>*{};<2mm,-7.7mm>*{}**@{-},
<-0.4mm,-5.5mm>*{};<-2mm,-7.7mm>*{}**@{-},
  <0.4mm,-5.5mm>*{};<2.9mm,-9.7mm>*{^2}**@{},
 <0.4mm,-5.5mm>*{};<-2.9mm,-9.7mm>*{^1}**@{},
 \end{xy}
=
\begin{xy}
 <0mm,0mm>*{\circ};
<7mm,0mm>*{a\sigma^{-1}};
<0mm,0.41mm>*{};<0mm,2.9mm>*{}**@{-},
 <0.39mm,-0.39mm>*{};<2.4mm,-2.4mm>*{}**@{-},
 <-0.35mm,-0.35mm>*{};<-2.4mm,-2.4mm>*{}**@{-},
 <-2.4mm,-2.4mm>*{};<-0.4mm,-4.5mm>*{}**@{-},
 <2.4mm,-2.4mm>*{};<0.4mm,-4.5mm>*{}**@{-},
  <0mm,-5.1mm>*{\circ};
 <5mm,-5.1mm>*{\sigma b};
  <0.4mm,-5.5mm>*{};<2mm,-7.7mm>*{}**@{-},
<-0.4mm,-5.5mm>*{};<-2mm,-7.7mm>*{}**@{-},
  <0.4mm,-5.5mm>*{};<2.9mm,-9.7mm>*{^2}**@{},
 <0.4mm,-5.5mm>*{};<-2.9mm,-9.7mm>*{^1}**@{},
 \end{xy}, \ \ \sigma\in \Z_2,
$$
$$
\lambda\left(
\Ba{c}
\begin{xy}
 <0mm,0mm>*{\circ};
<3mm,0mm>*{a};
<0mm,0.41mm>*{};<0mm,2.9mm>*{}**@{-},
 <0.39mm,-0.39mm>*{};<2.4mm,-2.4mm>*{}**@{-},
 <-0.35mm,-0.35mm>*{};<-2.4mm,-2.4mm>*{}**@{-},
 <-2.4mm,-2.4mm>*{};<-0.4mm,-4.5mm>*{}**@{-},
 <2.4mm,-2.4mm>*{};<0.4mm,-4.5mm>*{}**@{-},
  <0mm,-5.1mm>*{\circ};
 <3mm,-5.1mm>*{b};
  <0.4mm,-5.5mm>*{};<2mm,-7.7mm>*{}**@{-},
<-0.4mm,-5.5mm>*{};<-2mm,-7.7mm>*{}**@{-},
  <0.4mm,-5.5mm>*{};<2.9mm,-9.7mm>*{^2}**@{},
 <0.4mm,-5.5mm>*{};<-2.9mm,-9.7mm>*{^1}**@{},
 \end{xy}
\Ea
\right)
=
\Ba{c}
\begin{xy}
 <0mm,0mm>*{\circ};
<5mm,0mm>*{\lambda a};
<0mm,0.41mm>*{};<0mm,2.9mm>*{}**@{-},
 <0.39mm,-0.39mm>*{};<2.4mm,-2.4mm>*{}**@{-},
 <-0.35mm,-0.35mm>*{};<-2.4mm,-2.4mm>*{}**@{-},
 <-2.4mm,-2.4mm>*{};<-0.4mm,-4.5mm>*{}**@{-},
 <2.4mm,-2.4mm>*{};<0.4mm,-4.5mm>*{}**@{-},
  <0mm,-5.1mm>*{\circ};
 <3mm,-5.1mm>*{b};
  <0.4mm,-5.5mm>*{};<2mm,-7.7mm>*{}**@{-},
<-0.4mm,-5.5mm>*{};<-2mm,-7.7mm>*{}**@{-},
  <0.4mm,-5.5mm>*{};<2.9mm,-9.7mm>*{^2}**@{},
 <0.4mm,-5.5mm>*{};<-2.9mm,-9.7mm>*{^1}**@{},
 \end{xy}
\Ea
=
\Ba{c}
\begin{xy}
 <0mm,0mm>*{\circ};
<3mm,0mm>*{a};
<0mm,0.41mm>*{};<0mm,2.9mm>*{}**@{-},
 <0.39mm,-0.39mm>*{};<2.4mm,-2.4mm>*{}**@{-},
 <-0.35mm,-0.35mm>*{};<-2.4mm,-2.4mm>*{}**@{-},
 <-2.4mm,-2.4mm>*{};<-0.4mm,-4.5mm>*{}**@{-},
 <2.4mm,-2.4mm>*{};<0.4mm,-4.5mm>*{}**@{-},
  <0mm,-5.1mm>*{\circ};
 <5mm,-5.1mm>*{\lambda b};
  <0.4mm,-5.5mm>*{};<2mm,-7.7mm>*{}**@{-},
<-0.4mm,-5.5mm>*{};<-2mm,-7.7mm>*{}**@{-},
  <0.4mm,-5.5mm>*{};<2.9mm,-9.7mm>*{^2}**@{},
 <0.4mm,-5.5mm>*{};<-2.9mm,-9.7mm>*{^1}**@{},
 \end{xy}
\Ea
\ \ \ \forall \lambda \in \K,
$$
$$
\Ba{c}
\begin{xy}
 <0mm,0mm>*{\circ};
<7mm,0mm>*{a_1\hspace{-0.5mm}+\hspace{-0.5mm} a_2};
<0mm,0.41mm>*{};<0mm,2.9mm>*{}**@{-},
 <0.39mm,-0.39mm>*{};<2.4mm,-2.4mm>*{}**@{-},
 <-0.35mm,-0.35mm>*{};<-2.4mm,-2.4mm>*{}**@{-},
 <-2.4mm,-2.4mm>*{};<-0.4mm,-4.5mm>*{}**@{-},
 <2.4mm,-2.4mm>*{};<0.4mm,-4.5mm>*{}**@{-},
  <0mm,-5.1mm>*{\circ};
 <3mm,-5.1mm>*{b};
  <0.4mm,-5.5mm>*{};<2mm,-7.7mm>*{}**@{-},
<-0.4mm,-5.5mm>*{};<-2mm,-7.7mm>*{}**@{-},
  <0.4mm,-5.5mm>*{};<2.9mm,-9.7mm>*{^2}**@{},
 <0.4mm,-5.5mm>*{};<-2.9mm,-9.7mm>*{^1}**@{},
 \end{xy}
\Ea=
\Ba{c}
\begin{xy}
 <0mm,0mm>*{\circ};
<3mm,0mm>*{a_1};
<0mm,0.41mm>*{};<0mm,2.9mm>*{}**@{-},
 <0.39mm,-0.39mm>*{};<2.4mm,-2.4mm>*{}**@{-},
 <-0.35mm,-0.35mm>*{};<-2.4mm,-2.4mm>*{}**@{-},
 <-2.4mm,-2.4mm>*{};<-0.4mm,-4.5mm>*{}**@{-},
 <2.4mm,-2.4mm>*{};<0.4mm,-4.5mm>*{}**@{-},
  <0mm,-5.1mm>*{\circ};
 <3mm,-5.1mm>*{b};
  <0.4mm,-5.5mm>*{};<2mm,-7.7mm>*{}**@{-},
<-0.4mm,-5.5mm>*{};<-2mm,-7.7mm>*{}**@{-},
  <0.4mm,-5.5mm>*{};<2.9mm,-9.7mm>*{^2}**@{},
 <0.4mm,-5.5mm>*{};<-2.9mm,-9.7mm>*{^1}**@{},
 \end{xy}
\Ea
+
\Ba{c}
\begin{xy}
 <0mm,0mm>*{\circ};
<3mm,0mm>*{a_2};
<0mm,0.41mm>*{};<0mm,2.9mm>*{}**@{-},
 <0.39mm,-0.39mm>*{};<2.4mm,-2.4mm>*{}**@{-},
 <-0.35mm,-0.35mm>*{};<-2.4mm,-2.4mm>*{}**@{-},
 <-2.4mm,-2.4mm>*{};<-0.4mm,-4.5mm>*{}**@{-},
 <2.4mm,-2.4mm>*{};<0.4mm,-4.5mm>*{}**@{-},
  <0mm,-5.1mm>*{\circ};
 <3mm,-5.1mm>*{b};
  <0.4mm,-5.5mm>*{};<2mm,-7.7mm>*{}**@{-},
<-0.4mm,-5.5mm>*{};<-2mm,-7.7mm>*{}**@{-},
  <0.4mm,-5.5mm>*{};<2.9mm,-9.7mm>*{^2}**@{},
 <0.4mm,-5.5mm>*{};<-2.9mm,-9.7mm>*{^1}**@{},
 \end{xy}
\Ea
\ \ \ \mbox{and similarly for $b$}.
$$
are imposed.
It also follows from the definition that
$$
\Ba{c}
\begin{xy}
 <0mm,0mm>*{\circ};
<3mm,0mm>*{a};
<0mm,0.41mm>*{};<0mm,2.9mm>*{}**@{-},
 <0.39mm,-0.39mm>*{};<2.4mm,-2.4mm>*{}**@{-},
 <-0.35mm,-0.35mm>*{};<-2.4mm,-2.4mm>*{}**@{-},
 <-2.4mm,-2.4mm>*{};<-0.4mm,-4.5mm>*{}**@{-},
 <2.4mm,-2.4mm>*{};<0.4mm,-4.5mm>*{}**@{-},
  <0mm,-5.1mm>*{\circ};
 <3mm,-5.1mm>*{b};
  <0.4mm,-5.5mm>*{};<2mm,-7.7mm>*{}**@{-},
<-0.4mm,-5.5mm>*{};<-2mm,-7.7mm>*{}**@{-},
  <0.4mm,-5.5mm>*{};<2.9mm,-9.7mm>*{^2}**@{},
 <0.4mm,-5.5mm>*{};<-2.9mm,-9.7mm>*{^1}**@{},
 \end{xy}
\Ea =
\Ba{c}
\begin{xy}
 <0mm,0mm>*{\circ};
<3mm,0mm>*{a};
<0mm,0.41mm>*{};<0mm,2.9mm>*{}**@{-},
 <0.39mm,-0.39mm>*{};<2.4mm,-2.4mm>*{}**@{-},
 <-0.35mm,-0.35mm>*{};<-2.4mm,-2.4mm>*{}**@{-},
 <-2.4mm,-2.4mm>*{};<-0.4mm,-4.5mm>*{}**@{-},
 <2.4mm,-2.4mm>*{};<0.4mm,-4.5mm>*{}**@{-},
  <0mm,-5.1mm>*{\circ};
 <8mm,-5.1mm>*{b(12)};
  <0.4mm,-5.5mm>*{};<2mm,-7.7mm>*{}**@{-},
<-0.4mm,-5.5mm>*{};<-2mm,-7.7mm>*{}**@{-},
  <0.4mm,-5.5mm>*{};<2.9mm,-9.7mm>*{^1}**@{},
 <0.4mm,-5.5mm>*{};<-2.9mm,-9.7mm>*{^2}**@{},
 \end{xy}
\Ea, \ \ \ (12)\in \Z_2
$$
Thus one can define alternatively $G_0\langle E \rangle$ as a quotient space, $\prod_{v\in V(G)} E(Out_v, In_v)/\sim$,
with respect to the equivalence relation generated by the above pictures.

\sip

Note that if $E$ is a {\em differential
graded}\, (dg, for short) $\bS$-bimodule, then, for any graph $G\in \fG^\circlearrowright(m,n)$,
the associated  graded
vector space
$G\langle E \rangle$ comes equipped with an induced  $\bS_m\times \bS_n$-equivariant differential so
that the collection, $\{\bigoplus_{G\in \fG^\circlearrowright(m,n)} G\langle E \rangle\}_{m,n\geq 0}$,
 is again  a {\em dg}\, $\bS$-bimodule. The differential in $G\langle E\rangle$ induced from a differential $\delta$
on $E$ is denoted  by $\delta_G$ or, when no confusion may arise, simply by $\delta$.

\sip

\no{\bf 2.1.1. Definition.}
A {\em wheeled prop}\,  is an $\bS$-bimodule $\cP=\{\cP(m,n)\}$ together with
a family of linear $\bS_m\times \bS_n$-equivariant maps,
$$
\left\{\mu_G: G\langle \cP\rangle\rar \cP(m,n)\right\}_{G\in \fG^\circlearrowright(m,n), m,n\geq 0},
$$
parameterized by elements $G\in \fG^\circlearrowright$, which
satisfy the  condition
\Beq\label{graph-associativity}
\mu_G=\mu_{G/H}\circ \mu_H'
\Eeq
 for any subgraph $H\subset G$. Here $G/H$ is the graph obtained from $G$ by shrinking
 the whole subgraph $H$ into a single internal vertex, and
  $\mu_H': G\langle E \rangle \rar (G/H)\langle E\rangle$ stands for the map
which equals $\mu_H$ on the decorated vertices lying in $H$ and which is identity on all other vertices of $G$.

\sip

If the $\bS$-bimodule $\cP$ underlying a wheeled prop has a differential $\delta$ satisfying,
for any $G\in \fG^\circlearrowright$, the condition $\delta\circ \mu_G=\mu_G\circ \delta_G$, then the
wheeled prop $\cP$
is called {\em differential}.

\mip

\no{\bf 2.1.2.  Remarks.} {\bf (i)} If ${\mathfrak C}_{m,n}$ denotes $(m,n)$-corolla (\ref{corolla}),
then the $\bS_m\times \bS_n$-module
${\mathfrak C}_{m,n}\langle \cP  \rangle$ is canonically isomorphic to $\cP(m,n)$. Thus the defining linear map
$\mu_G: G\langle \cP\rangle\rar \cP(m,n)$ associated to an arbitrary graph $G\in \fG^\circlearrowright(m,n)$
can  be interpreted as a {\em contraction}\, map, $\mu_G: G\langle \cP\rangle\rar {\mathfrak C}_{m,n}\langle \cP \rangle$,
contracting all the internal edges and all the internal vertices of $G$ into a single vertex.

\mip

{\bf (ii)} Equation (\ref{graph-associativity}) implies $\mu_G= \mu_{G/G}\circ \mu_G$
for any graph $G\in \fG^\circlearrowright $, which in turn implies that $\mu_{{\mathfrak C}_{m,n}}$:
$\cP(m,n)\rar \cP(m,n)$ is the identity map.

\mip

{\bf (iii)}
 Condition (\ref{graph-associativity}) can be equivalently
rewritten as the equality, $\mu_{G/H_1}\circ \mu_{H_1}'=\mu_{G/H_2}\circ \mu_{H_2}'$, for any
subgraphs $H_1, H_2\subset G$, i.e.\ it is a kind of associativity condition for the family
of
 contraction operations $\{\mu_G\}$.

\mip

{\bf
(iv)} Strictly speaking, the notion introduced in \S 2.1.1 should be called a wheeled
prop {\em without unit}. A wheeled prop {\em with unit}\, can be defined as in \S 2.1.1
provided one  enlarges the family of graphs
$\fG^\circlearrowright$ by adding the following graphs without vertices,
$$
{\mathfrak t}_{p,q}:= \underbrace{
\uparrow \ \uparrow \ \uparrow \ \cdots \uparrow\ }_p
\underbrace{
\hgw\hgw\cdots \hgw}_q,\ p, q \geq 0, p+q\geq 1,
$$
 to the family $\fG^\circlearrowright(p,p)$ (see \cite{MMS}). The $\bS$-bimodule spanned by
such graphs without vertices has an obvious structure of wheeled prop with unit
called  the {\em trivial}\, wheeled prop ${\mathfrak t}$.
Similar to the case of an associative algebra,
 any wheeled prop, $\cP$, {\em without unit}\, can
be made into a wheeled prop, $\cP^+:=\cP * {\mathfrak t}$, {\em with unit}\,
by taking the free product of $\cP$ and ${\mathfrak t}$.
All the unital wheeled  props we study in this paper are obtained in this trivial way from
non-unital ones prompting us to work in this paper with non-unital props only. A small bonus
of this choice is that one can avoid bothering
about (co)augmentation (co)ideals
when dealing with  bar-cobar constructions of wheeled (co)props (see \S 2.4 below)

\mip

\no{\bf 2.1.3. Definitions.}
A {\em wheeled properad}, $\cP=\{\cP(m,n)\}$, is defined exactly as in \S 2.1.1 except that
the graphs $G$ and $H$ are required now to belong
to the subfamily, $\fG^\circlearrowright_c$, of $\fG^\circlearrowright$ consisting of
{\em connected}\, graphs.

\sip

A  {\em wheeled operad}\, is a wheeled properad $\cP=\{\cP(m,n)\}$ with $\cP(m,n)=0$ for $m\geq 2$.

\sip

\no{\bf 2.1.4. Generating compositions.} Associativity equations (\ref{graph-associativity})
imply that for an arbitrary
wheeled properad $\cP$  the defining family of contraction
maps, $\left\{\mu_G: G\langle \cP\rangle\rar \cP\right\}_{ G\in
\fG_c^\circlearrowright}$,
is uniquely determined (via iteration) by its subfamily,
$\left\{\mu_G: G\langle \cP\rangle\rar \cP\right\}_{ G\in {\fG}_{gen}^\circlearrowright}$,
where ${\fG}_{gen}^\circlearrowright\subset {\fG}^\circlearrowright_c$
consists of graphs of the form,
\Beq\label{2.1.4-generating compositions-i}
(i)\ \
 \begin{xy}
 <0mm,0mm>*{\circ};<0mm,0mm>*{}**@{},
 <-0.6mm,0.44mm>*{};<-8mm,5mm>*{}**@{-},
 <-0.4mm,0.7mm>*{};<-4.5mm,5mm>*{}**@{-},
 <0mm,0mm>*{};<0mm,5mm>*{\ldots}**@{},
 <0.4mm,0.7mm>*{};<4.5mm,5mm>*{}**@{-},
 <0.6mm,0.44mm>*{};<12.4mm,4.8mm>*{}**@{-},
 <-0.6mm,-0.44mm>*{};<-8mm,-5mm>*{}**@{-},
 <-0.4mm,-0.7mm>*{};<-4.5mm,-5mm>*{}**@{-},
 <0mm,0mm>*{};<-1mm,-5mm>*{\ldots}**@{},
 <0.4mm,-0.7mm>*{};<4.5mm,-5mm>*{}**@{-},
 <0.6mm,-0.44mm>*{};<8mm,-5mm>*{}**@{-},
 <13mm,5mm>*{};<13mm,5mm>*{\circ}**@{},
 <12.6mm,5.44mm>*{};<5mm,10mm>*{}**@{-},
 <12.6mm,5.7mm>*{};<8.5mm,10mm>*{}**@{-},
 <13mm,5mm>*{};<13mm,10mm>*{\ldots}**@{},
 <13.4mm,5.7mm>*{};<16.5mm,10mm>*{}**@{-},
 <13.6mm,5.44mm>*{};<20mm,10mm>*{}**@{-},
 <12.4mm,4.3mm>*{};<8mm,0mm>*{}**@{-},
 <12.6mm,4.3mm>*{};<12mm,0mm>*{\ldots}**@{},
 <13.4mm,4.5mm>*{};<16.5mm,0mm>*{}**@{-},
 <13.6mm,4.8mm>*{};<20mm,0mm>*{}**@{-},
 \end{xy}\ \ \ \ \
\mbox{and}
\ \ \ \ \ (ii)\ \ \
\begin{xy}
 <0mm,0mm>*{\circ};<0mm,0mm>*{}**@{},
 <-0.6mm,0.44mm>*{};<-8mm,5mm>*{}**@{-},
 <-0.4mm,0.7mm>*{};<-4.5mm,5mm>*{}**@{-},
 <0mm,0mm>*{};<0mm,4mm>*{\ldots}**@{},
 <0.4mm,0.7mm>*{};<4.5mm,5mm>*{}**@{-},
 <0.6mm,0.44mm>*{};<6mm,4mm>*{}**@{-},
 <-0.6mm,-0.44mm>*{};<-8mm,-5mm>*{}**@{-},
 <-0.4mm,-0.7mm>*{};<-4.5mm,-5mm>*{}**@{-},
 <0mm,0mm>*{};<0mm,-4mm>*{\ldots}**@{},
 <0.4mm,-0.7mm>*{};<4.5mm,-5mm>*{}**@{-},
 <0.6mm,-0.44mm>*{};<6mm,-4mm>*{}**@{-},
(6,4)*{}
   \ar@{->}@(ur,dr) (6,-4)*{}
 \end{xy}
\Eeq
i.e.\
of one-vertex graphs with precisely one internal edge (forming a loop)
and of  connected two vertex
graphs with precisely one internal edge. The set of graphs  $\fG_{gen}^\circlearrowright$
 lies behind the notion of a {\em quadratic}\, wheeled properad introduced below in \S 2.6.1.
\sip

Generating compositions  of a wheeled prop are given by graphs shown above and the extra ones,
\Beq\label{2.1.4-generating compositions2}
\begin{xy}
 <0mm,0mm>*{\circ};<0mm,0mm>*{}**@{},
 <-0.6mm,0.44mm>*{};<-7mm,4mm>*{}**@{-},
 <-0.4mm,0.7mm>*{};<-3.5mm,4mm>*{}**@{-},
 <0mm,0mm>*{};<0mm,3mm>*{...}**@{},
 <0.4mm,0.7mm>*{};<3.5mm,4mm>*{}**@{-},
 <0.6mm,0.44mm>*{};<7mm,4mm>*{}**@{-},
 <-0.6mm,-0.44mm>*{};<-7mm,-4mm>*{}**@{-},
 <-0.4mm,-0.7mm>*{};<-3.5mm,-4mm>*{}**@{-},
 <0mm,0mm>*{};<0mm,-3mm>*{...}**@{},
 <0.4mm,-0.7mm>*{};<3.5mm,-4mm>*{}**@{-},
 <0.6mm,-0.44mm>*{};<7mm,-4mm>*{}**@{-},
\end{xy}
\ \
\begin{xy}
 <0mm,0mm>*{\circ};<0mm,0mm>*{}**@{},
 <-0.6mm,0.44mm>*{};<-7mm,4mm>*{}**@{-},
 <-0.4mm,0.7mm>*{};<-3.5mm,4mm>*{}**@{-},
 <0mm,0mm>*{};<0mm,3mm>*{...}**@{},
 <0.4mm,0.7mm>*{};<3.5mm,4mm>*{}**@{-},
 <0.6mm,0.44mm>*{};<7mm,4mm>*{}**@{-},
 <-0.6mm,-0.44mm>*{};<-7mm,-4mm>*{}**@{-},
 <-0.4mm,-0.7mm>*{};<-3.5mm,-4mm>*{}**@{-},
 <0mm,0mm>*{};<0mm,-3mm>*{...}**@{},
 <0.4mm,-0.7mm>*{};<3.5mm,-4mm>*{}**@{-},
 <0.6mm,-0.44mm>*{};<7mm,-4mm>*{}**@{-},
\end{xy},
\Eeq
having two vertices and no internal edges.

\sip

\no{\bf 2.1.5. An endomorphism wheeled prop(erad).}
 For any finite-dimensional vector space $V$ the $\bS$-bimodule
$\cE nd_V:=\{ \Hom(V^{\ot n}, V^{\ot m})\}$ is naturally a wheeled prop(erad)
with compositions defined as follows:
\Bi
\item for graphs $G$ of the form (\ref{2.1.4-generating compositions-i})(i) the associated
composition $\mu_G: G\langle \cE nd_V \rangle \rar \cE nd_V$ is the ordinary composition of two linear maps;
\item for graphs $G$ of the form (\ref{2.1.4-generating compositions-i})(ii) the associated
composition $\mu_G$ is the ordinary trace of a linear map;

\item for graphs $G$ of the form (\ref{2.1.4-generating compositions2}) the associated
composition $\mu_G$ is the ordinary tensor product of linear maps.
\Ei
For an arbitrary graph $G\in \fG^\circlearrowright$ the associated composition
$\mu_G: G\langle \cE nd_V \rangle \rar \cE nd_V$ is defined as an iteration of the above ``elementary"
compositions, and it is easy to see that such a $\mu_G$ is independent of a particular choice of an iteration;
this independence means, in fact, that associativity conditions (\ref{graph-associativity}) are fulfilled.
The prop(erad) $\cE nd_V$ is called the {\em endomorphism wheeled prop(erad)}\, of $V$. Note that if $V$ is a complex,
then $\cE nd_V$ is naturally a {\em dg}\, prop(erad).

\sip

\no{\bf 2.1.6. A free wheeled prop(erad).}
Given an arbitrary  $\bS$-bimodule,
 $E=\{E(m,n)\}$, there is an associated $\bS$-bimodule,
$\cF^\circlearrowright\hspace{-0.2mm} \langle E\rangle=
\{\cF^\circlearrowright \hspace{-0.2mm}\langle E\rangle (m,n):=
\bigoplus_{G\in \fG^\circlearrowright(m,n)} G\langle E\rangle\}$,
which has a natural prop
structure  with the contraction maps
$\mu_G: G\langle\cF^\circlearrowright\hspace{-0.2mm} \langle E\rangle\rangle\rar
\cF^\circlearrowright\hspace{-0.3mm} \langle E\rangle
$ being tautological. The wheeled prop  $\cF^\circlearrowright\hspace{-0.2mm} \langle E\rangle$
 is called the {\em free wheeled prop
 generated
by an $\bS$-bimodule $E$}.

\sip

A {\em free wheeled properad, $\cF_c^\circlearrowright\hspace{-0.2mm} \langle E\rangle$,
 generated
by an $\bS$-bimodule $E$}\, is defined as in the previous paragraph but with the symbol $\fG^\circlearrowright$
replaced by  $\fG_c^\circlearrowright$.

\mip

{\bf 2.1.7. Prop(erad)s, dioperads and operads}.
Consider the follows subsets of the set  $\fG^\circlearrowright$:
\Bi
\item[(a)]
$\fG^\uparrow$  is a subset of $\fG^\circlearrowright$ consisting of directed graphs
with no {\em wheeles}, i.e.\
 directed paths of internal edges which begin and end at the same vertex;
\item[(b)]
$\fG_c^\uparrow:= \fG^\uparrow \cap \fG_c^\circlearrowright$;
\item[(c)]
$\fG^\uparrow_{c,0}$ is a subset of $\fG_c^\uparrow$ consisting of graphs of
 {\em genus zero};
\item[(d)]
$\fG^\uparrow_{oper}$ is a subset of $\fG^\uparrow_{c,0}$ built from corollas (\ref{corolla}) of type $(1,n)$ only,
 $n\geq 1$.
\Ei
Let $\fG^\checkmark$ be any one of these families of graphs. Then one can define
an $\fG^\checkmark$-algebra as in \S 2.1.1 by requiring that all the graphs
$G$, $H$ and $G/H$ involved in that
definition belong to the subset $\fG^\checkmark$ (cf. \cite{Me-Perm}). Then:
\Bi
\item[(a)]
an $\fG^\uparrow$- algebra is called a {\em prop} \cite{Mc};
\item[(b)]
an $\fG_c^\uparrow$-algebra is called a {\em properad} \cite{V};
\item[(c)]
an $\fG^\uparrow_{c,0}$-algebra is called a {\em dioperad} \cite{G};
\item[(d)]
an $\fG^\uparrow_{oper}$-algebra is called an {\em operad} \cite{May}.
\Ei

A {\em quadratic}\ $\fG^\checkmark$-algebra is defined (in all the above cases)
 as a quotient of a free $\fG^\checkmark$-algebra,
$\cF^\checkmark\hspace{-0.6mm}\langle E\rangle$, by the ideal  generated by a subspace
$R\subset \fG_{gen}^\checkmark\langle E\rangle$, where $\fG_{gen}^\checkmark$ is the
{\em minimal}\,
 subset of $\fG^\checkmark$ whose elements generate all possible compositions, $\mu_G$, via iteration
 (cf.\ \S 2.1.4). We apply the same minimality principle for the definition of a {\em quadratic}\, wheeled properad
 in \S 2.6 below.


\sip
\no{\bf 2.2. Morphisms of wheeled props.} One can make dg wheeled prop(erad)s into a category
by defining a morphism, $f: \cP_1\rar \cP_2$, as a morphism of the underlying dg $\bS$-bimodules,
$\{f: \cP_1(m,n)\rar \cP_2(m,n)\}_{m,n\geq 0}$,
 such that, for any graph $G\in \fG^\circlearrowright$,
 one has $f\circ \mu_G= \mu_G\circ (f^{\ot G})$,
where $f^{\ot G}$ means a map, $G\langle \cP_1\rangle \rar G\langle \cP_2\rangle$, which changes decorations of each
vertex
 in $G$ in accordance with $f$.

\sip

\no{\bf 2.2.1. Definition.}
A morphism of wheeled prop(erad)s, $\cP\rar \cE nd_V$, is called a {\em representation}\, of the
wheeled prop(erad)  $\cP$ in a graded
vector space $V$.

\sip

\no{\bf 2.2.2. Definition.}
A morphism of dg wheeled prop(erad)s, $\cP_1\rar \cP_2$, is called a {\em quasi-isomorphism}, if the induced
morphism of cohomology prop(erad)s, $H(\cP_1)\rar H(\cP_2)$, is an isomorphism.

\sip

\no{\bf 2.2.3. A useful fact.}
If $\cP_2$ is an arbitrary wheeled prop(erad) and  $\cP_1$ is a free wheeled prop(erad),
$\cF^\circlearrowright\hspace{-0.3mm}
 \langle E\rangle$, generated by some $\bS$-bimodule $E$,
then the set of morphisms of wheeled prop(erad)s, $\{f: \cP_1\rar \cP_2\}$, is in
one-to-one correspondence with the
vector space of degree zero morphisms of $\bS$-bimodules, $\{f|_E: E\rar \cP_2\}$, i.e.\ $f$ is uniquely
determined by its values on the generators. In particular, the set of morphisms,
$\cF^\circlearrowright\hspace{-0.3mm}
 \langle E\rangle\rar \cP_2$, has a graded vector space structure for any $\cP_2$.

\sip

\no{\bf 2.2.4. Definition.}
 A {\em free resolution}\, of a dg wheeled prop(erad)
$\cP$ is, by definition, a dg free wheeled prop(erad), $(\cF^\circlearrowright \hspace{-0.3mm}\langle E \rangle, \delta)$,
generated by some $\bS$-bimodule $E$ together with an epimorphism,
$\pi: (\cF^\circlearrowright \hspace{-0.2mm}\langle E \rangle, \delta) \rar \cP$, which is a quasi-isomorphism.
If the differential $\delta$ in $\cF^\circlearrowright \hspace{-0.2mm}\langle \cE \rangle$ is
decomposable with respect to the compositions $\mu_G$, then
$\pi: (\cF^\circlearrowright \hspace{-0.2mm}\langle E \rangle, \delta) \rar
\cP$ is called a {\em minimal model}\, of $\cP$.

\mip

\no{\bf 2.3. Coprop(erad)s.}
A {\em wheeled coproperad}\,  is an $\bS$-bimodule $\cP=\{\cP(m,n)\}$ together with
a family of linear $\bS_m\times \bS_n$-equivariant maps,
$$
\left\{\Delta_G: \cP(m,n) \rar G\langle \cP\rangle \right\}_{G\in \fG_c^\circlearrowright(m,n), m,n\geq 0},
$$
parameterized by elements $G\in \fG_c^\circlearrowright$, which
satisfy the  condition
\Beq\label{graph-coassociativity}
\Delta_G=\Delta_H' \circ \Delta_{G/H}
\Eeq
 for any connected
 subgraph $H\subset G$. Here $\Delta_H':  (G/H)\langle E\rangle \rar G\langle E \rangle$
 is the map
which equals $\Delta_H$ on the distinguished vertex of $G/H$  and which is identity on
all other vertices of $G$. {\em Wheeled coprops}\, are defined analogously.

\sip

If the $\bS$-bimodule $\cP$ underlying a wheeled coprop(erad) has a differential $\delta$ satisfying,
for any $G\in \fG^\circlearrowright$, the condition $\Delta_G\circ \delta=
\delta_G\circ \Delta_G$, then the
wheeled coprop(erad) $\cP$
is called {\em differential}.

\sip

For any $\bS$-bimodule,
 $E=\{E(m,n)\}$, the associated  $\bS$-bimodule,
$\cF^\circlearrowright\hspace{-0.2mm} \langle E\rangle$,
has  a natural coproperad
structure  with the co-contraction map
$$
\Delta:=\sum_{G\in \fG^\circlearrowright(m,n)}
\Delta_G: \cF^\circlearrowright\hspace{-0.3mm} \langle E\rangle
\lon \sum_{G\in \fG^\circlearrowright(m,n)}
 G\langle\cF^\circlearrowright\hspace{-0.2mm} \langle E\rangle\rangle =
\cF^\circlearrowright\hspace{-0.2mm} \langle\cF^\circlearrowright\hspace{-0.2mm}
\langle E\rangle\rangle
$$
given, on an arbitrary element $g\in G\langle E\rangle\subset  \cF^\circlearrowright\hspace{-0.3mm} \langle E\rangle$,
  by \cite{MMS},
$$
\Delta g = \sum_{f: Edg(G)\rar \{0,1\}} g_f
$$
where the sums runs over markings, $f: Edg(G)\rar \{0,1\}$, of the set,  $Edg(G)$,
 of internal edges of $G$ by numbers 0 and 1, and
 $g_f$ is an element of
 $\cF^\circlearrowright\hspace{-0.2mm} \langle\cF^\circlearrowright\hspace{-0.2mm}
\langle E\rangle\rangle$  obtained from $g$ by the following recipe:
\Bi
\item[(i)] cut every internal edge of the graph $G$
 marked by $0$ in the middle; let $G_1, \ldots, G_k$, for some $k\geq 1$,
be the resulting connected components of $G$; the vertices of the latter graphs inherit
$E$-decorations, and hence the marking $f$ defines elements $g_1\in G_1\langle E\rangle, \ldots,
g_k\in G_k\langle E\rangle$;
\item[(ii)] let $G'$ be the graph with $k$-vertices obtained from $G$ by shrinking each subgraph
$G_1, \ldots, G_k$ into a single vertex; then $g_f$ is, by definition, the decorated graph $g$
viewed as an element of $G'\langle \cF^\circlearrowright\hspace{-0.2mm}
\langle E\rangle\rangle   $, i.e.\ it equals $G'$ with vertices decorated
by elements $g_1, \ldots g_k \in \cF^\circlearrowright\hspace{-0.2mm}
\langle E\rangle$.
\Ei
The wheeled coprop $(\cF^\circlearrowright\hspace{-0.2mm}
\langle E\rangle, \Delta)$ is called the {\em free}\,  coprop generated by the
$\bS$-module $E$.

\sip

One can show analogously
that  $\cF_c^\circlearrowright\hspace{-0.2mm}
\langle E\rangle$ has a natural coproperad structure $\Delta$; the data
$(\cF_c^\circlearrowright\hspace{-0.2mm}
\langle E\rangle, \Delta)$ is
called the {\em free coproperad}\, generated by the
$\bS$-module $E$. We denote it by $\cF_{co}^\circlearrowright\hspace{-0.2mm}
\langle E\rangle$ (to avoid confusion with the natural properad structure
in $\cF_c^\circlearrowright\hspace{-0.2mm}
\langle E\rangle$).

\mip
\no{\bf 2.4.  Bar construction.}
With an $\bS$-module $E=\{E(m,n)\}$
one can associate two other $\bS$-bimodules,
$$
{\mathsf w} E=:\left\{ E(m,n)\ot sgn_n [-n] \right\}, \ \ \ \
{\mathsf w}^{-1} E:=\left\{ E(m,n)\ot sgn_n [n] \right\},
$$
where $sgn_m$ stands for the 1-dimensional sign representation of $\bS_m$.
We shall show in this subsection that for any properad $\cP$ the associated free coproperad,
$$
B(\cP):=\cF_{co}^\circlearrowright\hspace{-0.2mm}
\langle {\mathsf w}^{-1} \cP \rangle,
$$
comes canonically equipped  with a differential, $\delta_\cP$, encoding
all the generating
properadic compositions $\{\mu_G: G\langle\cP \rangle\rar \cP\}_{G\in \fG_{gen}^\circlearrowright}$.
For this purpose let us consider a family of graphs,
 $\fG^\circlearrowright_\bullet$, obtained
from the family of directed connected graphs $\fG_c^\circlearrowright$ by
inserting into each input leg and each internal edge of a graph
$G\in \fG_c^\circlearrowright$
a black $(1,1)$-corolla, $\xy
 <0mm,-2mm>*{};<0mm,2mm>*{}**@{-},
<0mm,0mm>*{\bullet};
\endxy$, and denoting the resulting graph  by $G_\bullet$. For example,
$$
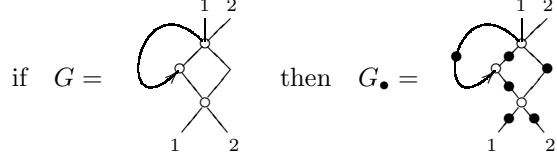
\mbox{if}\ \ \
G=
\Ba{c}
\begin{xy}
 <0mm,0mm>*{\circ};
<0.39mm,0.39mm>*{};<3.4mm,3.4mm>*{}**@{-},
<0mm,0.39mm>*{};<0mm,3.4mm>*{}**@{-},
 <0.39mm,-0.39mm>*{};<3.4mm,-3.4mm>*{}**@{-},
 <-0.35mm,-0.35mm>*{};<-2.9mm,-2.9mm>*{}**@{-},
 <-3.4mm,-3.4mm>*{\circ};
 <-3.0mm,-3.8mm>*{};<-0.4mm,-7.1mm>*{}**@{-},
 <3.4mm,-3.4mm>*{};<0.4mm,-7.1mm>*{}**@{-},
  <0mm,-7.8mm>*{\circ};
  <0.39mm,-8.5mm>*{};<3mm,-11.7mm>*{}**@{-},
<-0.39mm,-8.5mm>*{};<-3mm,-11.7mm>*{}**@{-},
  <0.4mm,-8.5mm>*{};<3.9mm,-14mm>*{^2}**@{},
 <0.4mm,-6.5mm>*{};<-3.9mm,-14mm>*{^1}**@{},
 <0mm,0mm>*{};<0mm,4.5mm>*{^1}**@{},
 <0mm,0mm>*{};<3.6mm,4.5mm>*{^2}**@{},
%
%
(-0.39,0.39)*{}
   \ar@{->}@(ul,dl) (-3.6,-3.6)*{}
 \end{xy}
\Ea
\ \ \ \mbox{then}\ \ \
G_\bullet=
\Ba{c}
\begin{xy}
 <0mm,0mm>*{\circ};
<0.39mm,0.39mm>*{};<3.4mm,3.4mm>*{}**@{-},
<0mm,0.39mm>*{};<0mm,3.4mm>*{}**@{-},
 <0.39mm,-0.39mm>*{};<3.4mm,-3.4mm>*{}**@{-},
 <-0.35mm,-0.35mm>*{};<-2.9mm,-2.9mm>*{}**@{-},
 <-3.4mm,-3.4mm>*{\circ};
 <-3.0mm,-3.8mm>*{};<-0.4mm,-7.1mm>*{}**@{-},
 <3.4mm,-3.4mm>*{};<0.4mm,-7.1mm>*{}**@{-},
  <0mm,-7.8mm>*{\circ};
  <0.39mm,-8.5mm>*{};<3mm,-11.7mm>*{}**@{-},
<-0.39mm,-8.5mm>*{};<-3mm,-11.7mm>*{}**@{-},
  <0.4mm,-8.5mm>*{};<3.9mm,-14mm>*{^2}**@{},
 <0.4mm,-6.5mm>*{};<-3.9mm,-14mm>*{^1}**@{},
 <0mm,0mm>*{};<0mm,4.5mm>*{^1}**@{},
 <0mm,0mm>*{};<3.6mm,4.5mm>*{^2}**@{},
<-1.7mm,-1.7mm>*{\bullet};
<-8.7mm,-1.7mm>*{\bullet};
<-1.7mm,-5.7mm>*{\bullet};
<3.4mm,-3.4mm>*{\bullet};
<-1.7mm,-10mm>*{\bullet};
<1.7mm,-10mm>*{\bullet};
(-0.39,0.39)*{}
   \ar@{->}@(ul,dl) (-3.6,-3.6)*{}
 \end{xy}
\Ea
$$
The automorphism group of such a graph $G_\bullet$
 is defined as in \S 2.1 with an extra assumption that the colour is preserved.
Then, obviously, $Aut(G)= Aut(G_\bullet)$.

\sip

Let $1$ stand for the unit in the field $\K$, and $\bar{1}$ for its image under the isomorphism $\K\rar \K[1]$.
The vector $\bar{1}$ has degree $-1$. For an arbitrary $\bS$-bimodule $E$ and an arbitrary graph
$G\in \fG_\bullet^\circlearrowright$ we denote
by $G_\bullet \langle E\rangle$ the vector space spanned by the graph $G_\bullet$ whose white vertices
are decorated by elements of $E$  and the special black $(1,1)$-vertices
are decorated by $\bar{1}$.

\sip
\no{\bf 2.4.1. Lemma.} {\em
For any $\bS$-module $E$ there is a canonical isomorphism of $\bS$-modules,}
$$
 \cF_c^\circlearrowright\hspace{-0.2mm}
\langle {\mathsf w}^{-1} E\rangle =
\bigoplus_{G_\bullet\in \fG_\bullet^\circlearrowright(m,n)} G_\bullet\langle E \rangle
$$

\sip

\Proof It is enough to show a canonical isomorphism $\bS_m\times \bS_n$-modules,
$G\langle {\mathsf w}^{-1} E\rangle= G_\bullet\langle E\rangle$ for an arbitrary graph
 $G\in \fG_c^\circlearrowright(m,n)$.
The graph $G_\bullet$ is obtained from  $G$ be replacing each constituting $(m,n)$-corolla of $G$
as follows,
$$
{\mathfrak C}_{m,n}=
\begin{xy}
 <0mm,0mm>*{\circ};
 <-0.5mm,0.2mm>*{};<-8mm,3mm>*{}**@{-},
 <-0.4mm,0.3mm>*{};<-4.5mm,3mm>*{}**@{-},
 <0mm,0mm>*{};<0mm,2.6mm>*{\ldots}**@{},
 <0.4mm,0.3mm>*{};<4.5mm,3mm>*{}**@{-},
 <0.5mm,0.2mm>*{};<8mm,3mm>*{}**@{-},
<-0.4mm,-0.2mm>*{};<-8mm,-3mm>*{}**@{-},
 <-0.5mm,-0.3mm>*{};<-4.5mm,-3mm>*{}**@{-},
 <0mm,0mm>*{};<0mm,-2.6mm>*{\ldots}**@{},
 <0.5mm,-0.3mm>*{};<4.5mm,-3mm>*{}**@{-},
 <0.4mm,-0.2mm>*{};<8mm,-3mm>*{}**@{-};
<0mm,5mm>*{\overbrace{\ \ \ \ \ \ \ \ \ \ \ \ \ \  }};
<0mm,-5mm>*{\underbrace{\ \ \ \ \ \ \ \ \ \ \ \ \ \ }};
<0mm,7mm>*{^{m\ \ output\ legs}};
<0mm,-7mm>*{_{n\ \ input\ legs}};
 \end{xy} \lon
{\mathfrak C}_{m,n\, \bullet}=
\begin{xy}
 <0mm,0mm>*{\circ};
 <-0.5mm,0.2mm>*{};<-8mm,3mm>*{}**@{-},
 <-0.4mm,0.3mm>*{};<-4.5mm,3mm>*{}**@{-},
 <0mm,0mm>*{};<0mm,2.6mm>*{\ldots}**@{},
 <0.4mm,0.3mm>*{};<4.5mm,3mm>*{}**@{-},
 <0.5mm,0.2mm>*{};<8mm,3mm>*{}**@{-},
<-0.4mm,-0.2mm>*{};<-8mm,-3mm>*{}**@{-},
 <-0.5mm,-0.3mm>*{};<-4.5mm,-3mm>*{}**@{-},
 <0mm,0mm>*{};<0mm,-2.6mm>*{\ldots}**@{},
 <0.5mm,-0.3mm>*{};<4.5mm,-3mm>*{}**@{-},
 <0.4mm,-0.2mm>*{};<8mm,-3mm>*{}**@{-};
<0mm,5mm>*{\overbrace{\ \ \ \ \ \ \ \ \ \ \ \ \ \  }};
<0mm,-5mm>*{\underbrace{\ \ \ \ \ \ \ \ \ \ \ \ \ \ }};
<0mm,7mm>*{^{m\ \ output\ legs}};
<0mm,-7mm>*{_{n\ \ input\ legs}};
 <-6mm,-2.2mm>*{\bullet};
<-2.8mm,-2.2mm>*{\bullet};
<6mm,-2.2mm>*{\bullet};
<2.8mm,-2.2mm>*{\bullet};
 \end{xy}
$$
It is obvious that ${\mathfrak C}_{m,n\, \bullet}\langle E \rangle=
{\mathfrak C}_{m,n}\langle {\mathsf w}^{-1} E\rangle$
as $\bS_m\times \bS_n$-bimodules. If we set $E(Out_v,In_v):= \bar{1}$ for every black vertex $v$ in
$G_\bullet$,
then
$\bigotimes_{v\in V(G_\bullet)} E(Out_v, In_v)= \bigotimes_{v\in V(G)}
{ \mathsf w}^{-1} E(Out_v, In_v)$ and the claim follows
finally from the isomorphism $Aut(G_\bullet)=Aut(G)$. \hfill $\Box$

\sip

{\bf 2.4.2. Corollary.}
{\em For any wheeled properad $\cP$ there is a canonical isomorphism
of $\bS$-modules, }
\Beq\label{Bar_constr_iso}
B(\cP)= \bigoplus_{G_\bullet\in \fG_\bullet^\circlearrowright} G_\bullet\langle \cP \rangle
\Eeq

\sip

The r.h.s of
 (\ref{Bar_constr_iso}) is denoted sometimes by $B_\bullet(\cP)$.

\mip

{\bf 2.4.3. Fact.} Let $\cP$ be an arbitrary wheeled properad. The $\bS$-module
 $\bigoplus_{G_\bullet\in \fG_\bullet^\circlearrowright} G_\bullet\langle \cP \rangle$
can be made naturally  into a complex with the differential,
$$
\delta_\cP=``\frac{\p}{\p \bullet_{edge}}",
$$
which is equal to zero an all white vertices and all black vertices attached to {\em legs},
and which
deletes a black vertex lying on every {\em internal edge} and contracts the associated
internal edge with the help of the corresponding  composition in $\cP$;
equation $\delta_\cP^2=0$
follows then from associativity conditions (\ref{graph-associativity}). More precisely,
one defines
 $\delta_\cP g$ for some $g\in G_\bullet\langle E\rangle=(\bigotimes_{v\in V(G_\bullet)}
 E(Out_v, In_v))_{Aut(G_\bullet)}$ as follows:
  choose first a representative,
 $\tilde{g}\in E(Out_{v_1}, In_{v_2})\ot \ldots
E(Out_{v_p}, In_{v_p})$,
of the equivalence class $g$ associated with some  ordering of all vertices in $G$, apply then
$\delta_\cP$ to the vertices of
$\tilde{g}$ in the chosen order, and finally set $\delta_\cP g= \pi(\delta_\cP \tilde{g})$,
where $\pi$ is the natural surjection
$$
\pi:  E(Out_{v_1}, In_{v_2})\ot \ldots\ot
E(Out_{v_p}, In_{v_p})\lon (\bigotimes_{v\in V(G_\bullet)}
 E(Out_v, In_v))_{Aut(G_\bullet)}.
$$
The result does not depend on the choice of a section, $g\rar \tilde{g}$, of $\pi$ used in the definition.
For example,
If
$$
g=
\Ba{c}
\begin{xy}
 <0mm,0mm>*{\circ};
<0.39mm,0.39mm>*{};<3.4mm,3.4mm>*{}**@{-},
<0mm,0.39mm>*{};<0mm,3.4mm>*{}**@{-},
 <0.39mm,-0.39mm>*{};<3.4mm,-3.7mm>*{}**@{-},
 <-0.35mm,-0.35mm>*{};<-2.9mm,-2.9mm>*{}**@{-},
 <-3.4mm,-3.4mm>*{\circ};
 <-3.0mm,-3.8mm>*{};<-0.4mm,-7.1mm>*{}**@{-},
 <3.4mm,-3.7mm>*{};<0.4mm,-7.1mm>*{}**@{-},
  <0mm,-7.8mm>*{\circ};
  <0.39mm,-8.5mm>*{};<3mm,-11.7mm>*{}**@{-},
<-0.39mm,-8.5mm>*{};<-3mm,-11.7mm>*{}**@{-},
  <0.4mm,-8.5mm>*{};<3.9mm,-14mm>*{^2}**@{},
 <0.4mm,-6.5mm>*{};<-3.9mm,-14mm>*{^1}**@{},
 <0mm,0mm>*{};<0mm,4.5mm>*{^1}**@{},
 <0mm,0mm>*{};<3.6mm,4.5mm>*{^2}**@{},
<-1.7mm,-1.7mm>*{\bullet};
<-6.7mm,1.7mm>*{\bullet};
<-1.7mm,-5.7mm>*{\bullet};
<3.4mm,-4.2mm>*{\bullet};
<-1.7mm,-10.2mm>*{\bullet};
<1.5mm,-9.5mm>*{\bullet};
<3.6mm,0.2mm>*{c};
<-5.5mm,-2.4mm>*{b};
<3.2mm,-7.4mm>*{a};
(-0.39,0.39)*{}
   \ar@{->}@(ul,dl) (-3.6,-3.6)*{}
 \end{xy}
\Ea
\ \ \ \mbox{for some}\ \ a\in E(2,2), b\in E(1,2), c\in E(3,2),
$$
then, ordering  the vertices from the bottom to the top, we obtain that
$\delta_\cP g$ is the equivalence class (in the unordered tensor product) of the following graph,
$$
(-1)^{a+b}\hspace{-3mm}
\Ba{c}
\begin{xy}
 <0mm,0mm>*{\circ};
<0.39mm,0.39mm>*{};<3.4mm,3.4mm>*{}**@{-},
<0mm,0.39mm>*{};<0mm,3.4mm>*{}**@{-},
 <0.39mm,-0.39mm>*{};<3.4mm,-3.4mm>*{}**@{-},
 <-0mm,-0.4mm>*{};<0mm,-7.1mm>*{}**@{-},
 <3.4mm,-3.4mm>*{};<0.4mm,-7.1mm>*{}**@{-},
  <0mm,-7.8mm>*{\circ};
  <0.39mm,-8.5mm>*{};<3mm,-11.7mm>*{}**@{-},
<-0.39mm,-8.5mm>*{};<-3mm,-11.7mm>*{}**@{-},
  <0.4mm,-8.5mm>*{};<3.9mm,-14mm>*{^2}**@{},
 <0.4mm,-6.5mm>*{};<-3.9mm,-14mm>*{^1}**@{},
 <0mm,0mm>*{};<0mm,4.5mm>*{^1}**@{},
 <0mm,0mm>*{};<3.6mm,4.5mm>*{^2}**@{},
<0mm,-2mm>*{\bullet};
<-4.5mm,1.9mm>*{\bullet};
<3.4mm,-3.4mm>*{\bullet};
<-1.7mm,-10.2mm>*{\bullet};
<1.5mm,-9.5mm>*{\bullet};
<4mm,0mm>*{c};
<3.6mm,-7.1mm>*{e_1};
(-0.39,0.39)*{}
   \ar@{->}@(ul,dl) (-0.6,-7.5)*{}
 \end{xy}
\Ea
\ - \ (-1)^{a+bc}\hspace{-3mm}\ \
\Ba{c}
\begin{xy}
 <0mm,0mm>*{\circ};
<-0.3mm,-7.6mm>*{};<-3.4mm,-5.4mm>*{}**@{-},
<-0.35mm,-7.6mm>*{};<-6mm,-5.4mm>*{}**@{-},
 <0.39mm,-0.39mm>*{};<3.4mm,-3.4mm>*{}**@{-},
 <-0mm,-0.4mm>*{};<0mm,-7.1mm>*{}**@{-},
 <3.4mm,-3.4mm>*{};<0.4mm,-7.1mm>*{}**@{-},
  <0mm,-7.8mm>*{\circ};
  <0.39mm,-8.5mm>*{};<3mm,-11.7mm>*{}**@{-},
<-0.39mm,-8.5mm>*{};<-3mm,-11.7mm>*{}**@{-},
  <0.4mm,-8.5mm>*{};<3.9mm,-14mm>*{^2}**@{},
 <0.4mm,-6.5mm>*{};<-3.9mm,-14mm>*{^1}**@{},
 <0mm,0mm>*{};<-7mm,-4.8mm>*{^1}**@{},
 <0mm,0mm>*{};<-3.0mm,-4.8mm>*{^2}**@{},
<0mm,-5mm>*{\bullet};
<-4.5mm,2.4mm>*{\bullet};
<3.4mm,-3.4mm>*{\bullet};
<-1.7mm,-10.2mm>*{\bullet};
<1.5mm,-9.5mm>*{\bullet};
<3mm,0.5mm>*{b};
<3.4mm,-7.3mm>*{e_2};
(-0.39,0.39)*{}
   \ar@{->}@(ul,dl) (-0.39,-0.39)*{}
 \end{xy}
\Ea
\ + \ (-1)^{a+b}\hspace{-3mm}
\Ba{c}
\begin{xy}
 <0mm,0mm>*{\circ};
<0.39mm,0.39mm>*{};<3.4mm,3.4mm>*{}**@{-},
<0mm,0.39mm>*{};<0mm,3.4mm>*{}**@{-},
 <0.39mm,-0.39mm>*{};<3.4mm,-3.4mm>*{}**@{-},
 <-0mm,-0.4mm>*{};<0mm,-7.1mm>*{}**@{-},
 <3.4mm,-3.4mm>*{};<0.4mm,-7.1mm>*{}**@{-},
  <0mm,-7.8mm>*{\circ};
  <0.39mm,-8.5mm>*{};<3mm,-11.7mm>*{}**@{-},
<-0.39mm,-8.5mm>*{};<-3mm,-11.7mm>*{}**@{-},
  <0.4mm,-8.5mm>*{};<3.9mm,-14mm>*{^2}**@{},
 <0.4mm,-6.5mm>*{};<-3.9mm,-14mm>*{^1}**@{},
 <0mm,0mm>*{};<0mm,4.5mm>*{^1}**@{},
 <0mm,0mm>*{};<3.6mm,4.5mm>*{^2}**@{},
<0mm,-5mm>*{\bullet};
<-4.5mm,2.4mm>*{\bullet};
<3.4mm,-3.4mm>*{\bullet};
<-1.7mm,-10.2mm>*{\bullet};
<1.5mm,-9.5mm>*{\bullet};
<4mm,0mm>*{e_3};
<3.4mm,-7.3mm>*{a};
(-0.39,0.39)*{}
   \ar@{->}@(ul,dl) (-0.39,-0.39)*{}
 \end{xy}
\Ea
\ - \ (-1)^{a+c}\hspace{-3mm}
\Ba{c}
\begin{xy}
 <0mm,0mm>*{\circ};
<0.39mm,0.39mm>*{};<3.4mm,3.4mm>*{}**@{-},
<0mm,0.39mm>*{};<0mm,3.4mm>*{}**@{-},
 <0.39mm,-0.39mm>*{};<3.4mm,-3.4mm>*{}**@{-},
 <-0mm,-0.4mm>*{};<0mm,-7.1mm>*{}**@{-},
 <3.4mm,-3.4mm>*{};<0.4mm,-7.1mm>*{}**@{-},
  <0mm,-7.8mm>*{\circ};
  <0.39mm,-8.5mm>*{};<3mm,-11.7mm>*{}**@{-},
<-0.39mm,-8.5mm>*{};<-3mm,-11.7mm>*{}**@{-},
  <0.4mm,-8.5mm>*{};<3.9mm,-14mm>*{^2}**@{},
 <0.4mm,-6.5mm>*{};<-3.9mm,-14mm>*{^1}**@{},
 <0mm,0mm>*{};<0mm,4.5mm>*{^1}**@{},
 <0mm,0mm>*{};<3.6mm,4.5mm>*{^2}**@{},
<0mm,-5mm>*{\bullet};
<-4.5mm,2.4mm>*{\bullet};
<3.4mm,-3.4mm>*{\bullet};
<-1.7mm,-10.2mm>*{\bullet};
<1.5mm,-9.5mm>*{\bullet};
<4mm,0mm>*{e_4};
<3.4mm,-7.3mm>*{a};
(-0.39,0.39)*{}
   \ar@{->}@(ul,dl) (-0.39,-0.39)*{}
 \end{xy}
\Ea
,
$$
where
$
e_1:=\mu\left(
\Ba{c}
\begin{xy}
 <0mm,0mm>*{\circ};
 <0.39mm,0.39mm>*{};<2.4mm,2.4mm>*{}**@{-},
 <-0.35mm,0.35mm>*{};<-1.9mm,1.9mm>*{}**@{-},
<0.39mm,-0.39mm>*{};<2.4mm,-2.4mm>*{}**@{-},
<-0.39mm,-0.39mm>*{};<-2.4mm,-2.4mm>*{}**@{-},
 <-2.4mm,2.4mm>*{\circ};
 <-2.6mm,2.2mm>*{};<-4.6mm,0mm>*{}**@{-},
 <-2.4mm,2.8mm>*{};<-2.4mm,5.5mm>*{}**@{-},
 <-4.2mm,2.9mm>*{b};
 <3mm,0mm>*{a};
 \end{xy}
\Ea
\right)\in E(2,3)$, etc. Applying  $\delta_\cP$ again and using
associativity relations (\ref{graph-associativity}),
one easily concludes that  $\delta_\cP^2=0$.

\sip

Isomorphism (\ref{Bar_constr_iso}) induces a differential in the free coproperad
$B(\cP)$ which we denoted by the same symbol $\delta_\cP$. It obviously respects
the coproperad structure in $B(\cP)$. If $\cP$ is a {\em differential}\, operad with differential
$d$, then the sum $d+\delta_\cP$ is a differential in $B(\cP)$.

\mip

{\bf 2.4.4. Definition.} The dg coproperad $(B(\cP), d+\delta_\cP)$ is called the
{\em bar construction}\, of a dg properad $(\cP, d)$.

\sip

This notion was first introduced in \cite{MMS} but with a different  $\bS$-module
structure and $\Z$-grading on $B(\cP)$. We shall be most interested below in the situations
when $d=0$.

\mip

{\bf 2.5. Cobar construction}. If $(\cC, d)$ is a dg coproperad,
then its {\em cobar construction}\, is, by definition, a free wheeled properad,
 $B^c(\cC):= \cF_c^\circlearrowright\langle{\mathsf w}\cC \rangle$, equipped with a differential,
$d+\p_\cC$, where $\p_\cC$ is the differential encoding the co-composition maps $\Delta_G: \cC\rar G\langle\cC\rangle$
in a way dual to the definition of $\p_\cP$ in \S 2.4 (see \cite{MMS}).
Let
 $\fG^\circlearrowright_\diamond$ be a family of graphs obtained
 from $\fG_c^\circlearrowright$ by
inserting into each input leg and each internal edge of a graph
$G\in \fG_c^\circlearrowright$
a white rhombic $(1,1)$-corolla, $\xy
 <0mm,-0.7mm>*{};<0mm,-2mm>*{}**@{-},
<0mm,0.7mm>*{};<0mm,2mm>*{}**@{-},
<0mm,0mm>*{\diamond};
\endxy$, and let us denote the resulting graph  by $G_\diamond$. Then, by analogy to \S 2.4.2, we have a canonical
degree $0$ isomorphism
of $\bS$-modules,
\Beq\label{White rhombic isomorphism}
B^c(\cC)= \bigoplus_{G_\diamond\in \fG_\diamond^\circlearrowright} G_\diamond\langle \cC \rangle
\Eeq
where in the r.h.s. we used $s(1)$, $s$ being the isomorphism $\K\rar \K[-1]$, to decorate special
 $\xy
 <0mm,-0.7mm>*{};<0mm,-2mm>*{}**@{-},
<0mm,0.7mm>*{};<0mm,2mm>*{}**@{-},
<0mm,0mm>*{\diamond};
\endxy$-vertices. The differential $\p_\cC$ is, by definition, equal to zero on the special
white rhombic  corollas while on ordinary
(decorated by $\cC$) vertices it is equal to the map
$\sum_G\Delta_G: \cC \rar \sum_{G\in \fG^\circlearrowright_{gen}}
 G\langle\cC\rangle$ with the sum running
over all possible graphs of the form (\ref{2.1.4-generating compositions-i}); the unique internal edge in the image of the map
$\sum_G\Delta_G$ is then  decorated by  $\xy
 <0mm,-0.7mm>*{};<0mm,-2mm>*{}**@{-},
<0mm,0.7mm>*{};<0mm,2mm>*{}**@{-},
<0mm,0mm>*{\diamond};
\endxy$ so that $\p_\cC$ increases the number of rhombic white vertices by one.

\sip

In the case when $\cC $ is the bar construction,
$(\cC=B(\cP), d+\delta_\cP)$ on some dg properad $(\cP,d)$ one has a natural epimorphism of dg properads,
$$
\bar{\pi}:
\left(B^c(B(\cP)), \delta:=
d+ \delta_\cP + \p_{B(\cP)}\right) \lon \left({\mathsf w}({\mathsf w}^{-1}\cP)= \cP, d\right),
$$
which is a quasi-isomorphism \cite{MMS}. If we now apply constructions
(\ref{Bar_constr_iso}) and
(\ref{White rhombic isomorphism}) to $B^c(B(\cP)$, we shall get decorated graphs whose internal
edges are decorated by either black vertices or simultaneously by black and white rhombic vertices.
As white rhombic and black corollas placed on the same edge ``annihilate"
 each other with respect to the their total impact on graph, we have
 a degree $0$ isomorphism of  $\bS$-bimodules,
$$
B^c(B(\cP))=  \bigoplus_{G_\diamond\in \fG_{\bullet,-}^\circlearrowright} G_{\bullet,-}
\langle \cP \rangle
$$
where, by definition, $ \fG_{\bullet,-}$ is a family of graphs obtained
 from graphs in  $\fG_c^\circlearrowright$ by
inserting into {\em some}\, (possibly, none) internal edges  black $(1,1)$-corollas; thus every {\em internal}\,
 edge  of a graph $G$ from  $ \fG_{\bullet,-}$ is either straight or equipped with the black  $(1,1)$-corolla,
and every input
or output leg of $G$ is straight.
Then the differential $\p_{B(\cP)}$ gets a very simple interpretations ---
it eliminates, in accordance with the Leibnitz rule, each black corolla making the corresponding edge straight;
on the other hand, the differential $\delta_\cP$ contracts (again in accordance with the Leibnitz rule) each internal edge decorated by the black corolla and performs a corresponding to this contraction composition in the original
properad $\cP$.

\mip
{\bf 2.6. Quadratic wheeled (co)properads and  Koszul duality}.
Koszul duality for ordinary quadratic operads was introduced in \cite{GK}, for dioperads in \cite{G},
for ordinary properads in \cite{V} and for wheeled operads in \cite{MMS}. In this section we extend
the idea to arbitrary quadratic wheeled properads.

\sip

For a graph $G\in \fG^\circlearrowright$ with $p$ vertices and $q$ wheels (that is, closed paths
of directed internal edges) set $||G||:=p+q$
and call it the {\em weight}\, of $G$.   For an $\bS$-module $E$ set
$\cF_{(\lambda)}^\circlearrowright\langle E \rangle$
to be a submodule of the free properad $\cF_c^\circlearrowright\langle E \rangle$
spanned by decorated graphs of weight $\la$. Note that properadic compositions
$\{\mu_G: G\langle\cF_c^\circlearrowright\langle E \rangle\rangle\rar
\cF_c^\circlearrowright\langle E \rangle\}_{G\in \fG_c^\circlearrowright}$
are homogeneous with respect to this weight gradation.
 Note also that the quadratic subspace
$\cF^\circlearrowright_{(2)}\langle E \rangle\subset
\cF_c^\circlearrowright\langle E \rangle$ is distinguished as it is spanned,
$$
\cF^\circlearrowright_{(2)}\langle E \rangle = \sum_{G\in \fG_{gen}^\circlearrowright}
G\langle E\rangle,
$$
by the minimal set of graphs (\ref{2.1.4-generating compositions-i})
which generate {\em all}\, possible wheeled properadic compositions.

\mip

{\bf 2.6.1. Definition.}
 A wheeled properad $\cP$ is called {\em quadratic}\, if it is the
  quotient, $\cP:= \cF_c^\circlearrowright\langle E\rangle / I$, of a free wheeled properad
 (generated by some $\bS$-bimodule $E$) by the ideal, $I$,  generated
by some subspace $R\subset \cF^\circlearrowright_{(2)}\langle E \rangle$.

\mip

An obvious dualization of the above definition gives
the notion of a {\em quadratic}\, coproperad.

\mip

Any quadratic (co)properad, $\cP$, comes equipped with an induced weight gradation,
$\cP=\sum_{\lambda\geq 1} \cP_{(\la)}$, where $\cP_{(\la)}$ is the image of
$\cF_{(\lambda)}^\circlearrowright\langle E \rangle$ under the natural surjection
$\cF_{c}^\circlearrowright\langle E \rangle\rar \cP$. Note that $\cP_{(1)}=E$ and $\cP_{(2)}$
is given by an
exact sequence
\Beq\label{Relations quad}
0 \lon \cR \lon \cF_{(2)}^\circlearrowright\langle E \rangle \lon \cP_{(2)}\lon 0.
\Eeq
The subspace $\cB(\cP_{(1)})\subset \cB(\cP)$ is obviously a sub-coproperad, but, in general,  it is not preserved
by the bar differential $\p_\cP$. It is not hard to check that
an $\bS$-bimodule $\cP^{{\mbox {\scriptsize !`}}}$ defined by the exact sequence,
$$
0\lon \cP^{\mbox {\scriptsize !`}} \lon B(\cP_{(1)}) \stackrel{\p_\cP}{\lon}
 B(\cP)[1],
$$
is a sub-coproperad of $\cB(\cP_{(1)})$
so that the natural composition of inclusions,
\Beq\label{i from P bot}
\imath: \cP^{\mbox {\scriptsize !`}} \lon B(\cP_{(1)})  \lon B(\cP),
\Eeq
is a monomorphism of {\em dg}\, wheeled coproperads.

\mip

{\bf 2.6.2. Definition.} The coproperad $\cP^{\mbox {\scriptsize !`}}$
is called {\em Koszul
dual}\, to a quadratic  wheeled properad $\cP$.

\mip

{\bf 2.6.3. Definition }.
A quadratic wheeled properad $\cP$ is called {\em Koszul}, if the associated
morphism of dg coproperads, $\imath: \cP^{\mbox {\scriptsize !`}}  \lon B(\cP)$, is
 a quasi-isomorphism.

\mip

As the cobar construction functor $B^c$ is exact \cite{MMS},
 the composition
$$
\pi: \cP_\infty:=B^c(\cP^{\mbox {\scriptsize !`}}) \stackrel{B^c(i)}{\lon} B^c(B(\cP))
\stackrel{\bar{\pi}}{\lon} \cP
$$
is a quasi-isomorphism if $\cP$ is Koszul; then the dg free wheeled properad $\cP_\infty$ gives us
 a minimal resolution of $\cP$.

\mip

{\bf 2.6.4. Remark on notation}.
 In general (i.e.\ if $\cP$ is not Koszul), the  dg
  properad $B^c(\cP^{\mbox {\scriptsize !`}})$
 is only an approximation to the genuine minimal wheeled resolution of $\cP$ (if it exists
at all); it is, however,
 associated {\em canonically}\, to $\cP$, and, slightly abusing tradition, {\em we continue denoting it
in this paper
by $\cP_\infty$ even in the cases when $\cP$ is not Koszul}.
\sip

\mip

Note that $B(\cP_{(1)})$ is the free co-properad generated by the $\bS$-module
$\mathsf w^{-1}\cP_{(1)}={\mathsf w}^{-1}E$. By the definition of the bar differential $\p_\cP$,
 the image, $I^{co}$,  of the degree $0$ map $\p_\cP:B(\cP_{(1)}) \rar
 B(\cP)[1]$ is spanned by graphs with all (except one!) vertices decorated by the
 $\mathsf w^{-1}\cP_{(1)}$ and with the exceptional vertex decorated by
 $\mathsf w^{-1}\cP_{(2)}[1]$. Thus we have an exact sequence,
 $$
 0\lon \cP^{\mbox {\scriptsize !`}} \lon \cF_{co}^\circlearrowright \langle {\mathsf w}^{-1}E
 \rangle \lon I^{co} \lon 0,
 $$
As $\mathsf w^{-1} \cF_{(2)}^\circlearrowright\langle E \rangle =
\cF_{(2)}^\circlearrowright\langle \mathsf w^{-1}E \rangle[-1]$, we can rewrite
(\ref{Relations quad}) as follows,
\Beq\label{Relations quad w}
0 \lon  {\mathsf w}^{-1} \cR[1] \lon \cF_{(2)}^\circlearrowright\langle {\mathsf w}^{-1} E \rangle \lon
{\mathsf w}^{-1} \cP_{(2)}[1]\lon 0,
\Eeq
and conclude that $I^{co}$ is the co-ideal of $\cF_{cc}^\circlearrowright \langle {\mathsf w}^{-1}E
 \rangle $ cogenerated by {\em quadratic}\,  co-relations $\mathsf w^{-1}\cP_{(2)}[1]$. Hence
 we proved the following
 \mip

{\bf 2.6.5. Proposition}. {\em For any quadratic wheeled properad $\cP$ the associated Koszul dual
wheeled coproperad
 $\cP^{\mbox {\scriptsize !`}}$ is quadratic}.

 \mip

 {\bf 2.6.6. Remark.}
If the $\bS$-bimodule  $E=\{E(m,n)\}$ is of finite type (i.e.\ each $E(m,n)$ is finite-dimensional),
then it is often easier to work with the wheeled properad $\cP^!:= (\cP^{\mbox {\scriptsize !`}})^*$,
the ordinary dual of the coproperad $\cP^{\mbox {\scriptsize !`}}$. It is a quadratic
wheeled properad  generated by the $\bS$-bimodule,
$$
E^\vee:=\left\{E(m,n)^* \ot sgn_n[-n]\right\},
$$
with the relations, $\cR^\bot$, given by the exact sequence,
\Beq\label{Relations for P!}
0\lon \cR^\bot \lon \cF_{(2)}^\circlearrowright\langle E^\vee\rangle \lon
{\mathsf w}\cR^*[-1] \lon 0,
\Eeq
where $\cR$ are the quadratic relations for $\cP$.

\mip

{\bf 2.6.7. Remark.} Definition 2.6.1 implies that there exists a canonical {\em wheelification
functor},
\Beq\label{wheelification functor}
\Ba{rccc}
^\circlearrowright: &
\Ba{c}\mbox{\sf Category of quadratic}\\
\mbox{\sf (co)dioperads}
\Ea &
\lon &
\Ba{c}\mbox{\sf Category of quadratic wheeled}\\
\mbox{\sf (co)properads}
\Ea
\\
& \caD & \lon & \caD^\circlearrowright
\Ea
\Eeq
which is, by definition, identity on the (co)generators and the quadratic (co)relations of
the dioperad $\caD$. It is worth noting that, in general,
$(\caD^{\mbox {\scriptsize !`}})^\circlearrowright \neq
 (\caD^\circlearrowright)^{\mbox {\scriptsize !`}}$,
implying that $(\caD^\circlearrowright)_\infty$ may be substantially larger that
$(\caD_\infty)^\circlearrowright$, where $\caD_\infty$ stands for the cobar construction
on the Koszul dual co-dioperad $\caD^{\mbox {\scriptsize !`}}$ in the category of dioperads;
we refer to \cite{MMS} for explicit examples of this phenomenon
for the cases $\caD=\cA ss$ and $\caD=\cC omm$, the operads of associative and, respectively,
commutative algebras. In the case of the operad of Lie
algebras one actually has an equality, $(\caL ie^{\mbox {\scriptsize !`}})^\circlearrowright =
 (\caL ie^\circlearrowright)^{\mbox {\scriptsize !`}}$ (see \cite{Me-graphs}).

\sip

The wheelification functor does not, in general, preserve Koszulness: a Koszul dioperad, $\caD$, may have
a non-Koszul wheelification, $\caD^\circlearrowright$. We give an explicit example of
this phenomenon in \S 4. It is worth noting in this connection that the functor $^\circlearrowright$
applied to the three classical operads $\cA ss$, $\cC omm$, and $\caL ie$ does preserve  Koszulness
(see \cite{Me-graphs, MMS} for the proofs).

\mip

{\bf 2.7. Homotopy Transfer Formulae.} Let $(V,d)$ and $(W,d)$ be dg vector spaces
equipped with linear degree 0 maps of complexes, $i: W \rar V$ and $p: V\rar W$,
such that the composition  $i\circ p: V\rar V$ is homotopy equivalent to the identity map,
$\Id: V\rar V$,
\Beq\label{cohomogical splitting}
Id_V=i\circ p + d\circ h + h\circ d,
\Eeq
via a fixed homotopy $h: V\rar V[-1]$. Without loss of generality we may assume
that the data $(i,p,h)$ satisfies the so called {\em side conditions} \cite{LS},
$$
p\circ i=\Id_W, \ \ \ p\circ h=0,\ \ \ h\circ i=0, \ \ \ h\circ h=0.
$$
When $W$ is the cohomology of the complex $V$ the above data is often called
a {\em cohomological splitting}\ of $(V,d)$.

\mip

{\bf 2.7.1. Theorem}.
{\em Let $\cP$ be a quadratic
 wheeled properad, and  $\rho: \cP\rar \cE nd_V$  an arbitrary
  $\cP$-algebra structure on the complex $V$.
  For any element
$G\in {\mathsf w}\cP^{\mbox {\scriptsize !`}}(m,n)$\,
let
$$
G\langle i,h,p,\rho\rangle\in \cE nd_W(m,n),
$$
be a linear map $W^{\ot n}\rar W^{\ot n}$ defined
as follows:
\Bi
\item[(i)] consider the image,
$\imath(G)$,
of $G$
under the canonical inclusion
$\imath:  {\mathsf w}\cP^{\mbox {\scriptsize !`}}\rar {\mathsf w}B_\bullet(\cP)$;
\item[(ii)]
decorate the input legs of  each graph summand in the image
$\imath(G)$ with $i$, the output legs with $p$,
and the special vertices, $\bullet$, lying on the internal edges   with $h$,
\item[(iii)]  replace a decoration, $e$, of every non-special vertex in $\imath(G)$
by $\rho(e)$,
and finally
\item[(iv)] interpret the resulting decorated graph as a scheme for the composition of maps
$i$, $h$, $\rho(e)$ and $p$.
\Ei
Then the family of maps,
$$
\left\{G \lon
 G\langle i,h,p, \rho\rangle\in \cE nd_{W}\right\}_{G\in  {\mathsf w}\cP^{\mbox {\scriptsize !`}}},
$$
defines a representation of the  dg free wheeled properad $\cP_\infty$ in the
dg space $W$.
}

\mip

\Proof Any morphism, $\cP_\infty=\cF^\circlearrowright\langle
{\mathsf w}\cP^{\mbox {\scriptsize !`}}\rangle \rar \cE nd_{W}$, of wheeled properads is
uniquely determined by its values
on the generators, i.e.\ by a morphism,
${\mathsf w}^{-1}\cP^{\mbox {\scriptsize !`}} \lon  \cE nd_{W}$ of $\bS$-modules.
 Define such a morphism, $\rho_\infty: \cP_\infty \rar \cE nd_{W}$,  by setting
$$
\rho_\infty(G):= G\langle i,h,p,\rho\rangle.
$$
This morphism gives  a representation of the dg properad $\cP_\infty$ if and only if
it respects the differentials, i.e.
$$
\rho_\infty(\p_{B(\cP)} G)=d (G\langle i,h,p\rangle),
$$
where $\p_{B(\cP)}$ is the differential in $B^c(B(\cP))$ restricted to
the subcomplex $B^c(\cP^\Koz)$,  and $d$ is the differential
in $\cE nd_W$ (induced by the differential $d$ in $W$ and denoted by the same letter).
%

\sip

Let us assume, for an illustration,  that
$$G=
\Ba{c}
\begin{xy}
<-3mm,7.8mm>*{};<-6mm,10.8mm>*{}**@{-},
<3mm,7.8mm>*{};<6mm,10.8mm>*{}**@{-},
 <0mm,6mm>*{\mbox{$\xy *=<6mm,4mm>
\txt{{{$a$}}}*\frm{-}\endxy$}};
<-3mm,3.8mm>*{};<-6mm,0mm>*{}**@{-},
<3mm,3.8mm>*{};<6mm,0mm>*{}**@{-},
<-3mm,-3.8mm>*{};<-6mm,0mm>*{}**@{-},
<3mm,-3.8mm>*{};<6mm,0mm>*{}**@{-},
 <0mm,-6mm>*{\mbox{$\xy *=<6mm,4mm>
\txt{{{$b$}}}*\frm{-}\endxy$}};
 <0mm,-7.9mm>*{};<0mm,-10.9mm>*{}**@{-},
 \end{xy} \in {\mathsf w}\cP^{\mbox {\scriptsize !`}}
\Ea
$$
for some $a,b\in \cP$.
Then
$$
G\langle i,h,p,\rho\rangle=
\Ba{c}
\begin{xy}
<-3mm,7.8mm>*{};<-6mm,10.8mm>*{}**@{-},
<3mm,7.8mm>*{};<6mm,10.8mm>*{}**@{-},
 <0mm,6mm>*{\mbox{$\xy *=<6mm,4mm>
\txt{{{$_{\rho(a)}$}}}*\frm{-}\endxy$}};
<-3mm,3.8mm>*{};<-6mm,0mm>*{}**@{-},
<3mm,3.8mm>*{};<6mm,0mm>*{}**@{-},
<-3mm,-3.8mm>*{};<-6mm,0mm>*{}**@{-},
<3mm,-3.8mm>*{};<6mm,0mm>*{}**@{-},
 <0mm,-6mm>*{\mbox{$\xy *=<6mm,4mm>
\txt{{{$_{\rho(b)}$}}}*\frm{-}\endxy$}};
 <0mm,-8mm>*{};<0mm,-10.9mm>*{}**@{-},
  <6mm,0mm>*{\bullet};
   <-6mm,0mm>*{\bullet};
 <-7.5mm,11.5mm>*{^p};
<7.5mm,11.5mm>*{^p};
<0mm,-12.5mm>*{_i};
<-7.5mm,0mm>*{_h};
<8mm,0mm>*{_h};
 \end{xy}
\Ea
$$
and, using (\ref{cohomogical splitting}), we obtain
\Beqrn
d(G\langle i,h,p,\rho\rangle)&=&
(-1)^b\
\begin{xy}
<-3mm,7.8mm>*{};<-6mm,10.8mm>*{}**@{-},
<3mm,7.8mm>*{};<6mm,10.8mm>*{}**@{-},
 <0mm,6mm>*{\mbox{$\xy *=<6mm,4mm>
\txt{{{$_{\rho(a)}$}}}*\frm{-}\endxy$}};
<-3mm,3.8mm>*{};<-6mm,0mm>*{}**@{-},
<3mm,3.8mm>*{};<6mm,0mm>*{}**@{-},
<-3mm,-3.8mm>*{};<-6mm,0mm>*{}**@{-},
<3mm,-3.8mm>*{};<6mm,0mm>*{}**@{-},
 <0mm,-6mm>*{\mbox{$\xy *=<6mm,4mm>
\txt{{{$_{\rho(b)}$}}}*\frm{-}\endxy$}};
 <0mm,-8mm>*{};<0mm,-10.9mm>*{}**@{-},
   <6mm,0mm>*{\bullet};
 <-7.5mm,11.5mm>*{^p};
<7.5mm,11.5mm>*{^p};
<0mm,-12.5mm>*{_i};
<8mm,0mm>*{_h};
<-7.8mm,0mm>*{_\Id}
 \end{xy}
\ \ -\ (-1)^b\
\begin{xy}
<-3mm,7.8mm>*{};<-6mm,10.8mm>*{}**@{-},
<3mm,7.8mm>*{};<6mm,10.8mm>*{}**@{-},
 <0mm,6mm>*{\mbox{$\xy *=<6mm,4mm>
\txt{{{$_{\rho(a)}$}}}*\frm{-}\endxy$}};
<-3mm,3.8mm>*{};<-6mm,0mm>*{}**@{-},
<3mm,3.8mm>*{};<6mm,0mm>*{}**@{-},
<-3mm,-3.8mm>*{};<-6mm,0mm>*{}**@{-},
<3mm,-3.8mm>*{};<6mm,0mm>*{}**@{-},
 <0mm,-6mm>*{\mbox{$\xy *=<6mm,4mm>
\txt{{{$_{\rho(b)}$}}}*\frm{-}\endxy$}};
 <0mm,-8mm>*{};<0mm,-10.9mm>*{}**@{-},
   <-6mm,0mm>*{\bullet};
 <-7.5mm,11.5mm>*{^p};
<7.5mm,11.5mm>*{^p};
<0mm,-12.5mm>*{_i};
<-7.5mm,0mm>*{_h};
<8.2mm,0mm>*{_\Id}
 \end{xy}\\
 &&
 -(-1)^b
\begin{xy}
<-3mm,7.8mm>*{};<-6mm,10.8mm>*{}**@{-},
<3mm,7.8mm>*{};<6mm,10.8mm>*{}**@{-},
 <0mm,6mm>*{\mbox{$\xy *=<6mm,4mm>
\txt{{{$_{\rho(a)}$}}}*\frm{-}\endxy$}};
<-3mm,3.8mm>*{};<-6mm,0mm>*{}**@{-},
<3mm,3.8mm>*{};<6mm,0mm>*{}**@{-},
<-3mm,-3.8mm>*{};<-6mm,0mm>*{}**@{-},
<3mm,-3.8mm>*{};<6mm,0mm>*{}**@{-},
 <0mm,-6mm>*{\mbox{$\xy *=<6mm,4mm>
\txt{{{$_{\rho(b)}$}}}*\frm{-}\endxy$}};
 <0mm,-8mm>*{};<0mm,-10.9mm>*{}**@{-},
   <6mm,0mm>*{\bullet};
 <-7.5mm,11.5mm>*{^p};
<7.5mm,11.5mm>*{^p};
<0mm,-12.5mm>*{_i};
<8mm,0mm>*{_h};
<-6mm,-3mm>*{_p};
<-6mm,3mm>*{_i};
 \end{xy}
\ \ +\ (-1)^b\
\begin{xy}
<-3mm,7.8mm>*{};<-6mm,10.8mm>*{}**@{-},
<3mm,7.8mm>*{};<6mm,10.8mm>*{}**@{-},
 <0mm,6mm>*{\mbox{$\xy *=<6mm,4mm>
\txt{{{$_{\rho(a)}$}}}*\frm{-}\endxy$}};
<-3mm,3.8mm>*{};<-6mm,0mm>*{}**@{-},
<3mm,3.8mm>*{};<6mm,0mm>*{}**@{-},
<-3mm,-3.8mm>*{};<-6mm,0mm>*{}**@{-},
<3mm,-3.8mm>*{};<6mm,0mm>*{}**@{-},
 <0mm,-6mm>*{\mbox{$\xy *=<6mm,4mm>
\txt{{{$_{\rho(b)}$}}}*\frm{-}\endxy$}};
 <0mm,-8mm>*{};<0mm,-10.9mm>*{}**@{-},
   <-6mm,0mm>*{\bullet};
 <-7.5mm,11.5mm>*{^p};
<7.5mm,11.5mm>*{^p};
<0mm,-12.5mm>*{_i};
<-7.5mm,0mm>*{_h};
<6mm,-3mm>*{_p};
<6mm,3mm>*{_i};
 \end{xy}\\
 &=&\rho(\p_\cP G) +  \rho(\p_{B(\cP)} G)\\
 &=&\rho(\p_{B(\cP)} G).
\Eeqrn
In the above calculation we used
\Bi
\item[-]
the identification of $\p_\cP$ with a machine deleting the black vertices and
contracting the associated internal edge (so that the first two terms in the above sum of 4 graphs
are precisely $\rho(\p_\cP G\langle i,h,p,\rho\rangle)$),
\item[-]
 the identification of $\p_{B(\cP)}$ with a machine deleting the black vertices without
 subsequent contraction of the associated internal edge
  (so that the last two terms in the above sum
are precisely $\rho(\p_{B(\cP)} G\langle i,h,p,\rho\rangle)$), and
\item[-] the fact that, by definition of $\cP^\Koz$, one has
 $\rho(\p_{\cP} G)=0$ for any $G\in   {\mathsf w}\cP^{\mbox {\scriptsize !`}}$.
\Ei

The pattern exposed above is universal, i.e.\ it does not depend on the particularities of  $G$.
This simple calculation proves the claim. \hfill
$\Box$

\mip

{\bf 2.7.2. Remarks}.
(i) The above arguments work for {\em any}\, dg sub-coproperad of $B(\cP)$, not only for $\cP^\Koz$.

\sip

(ii) Theorem 2.7.1 gives a conceptual explanation of the well-known  ``experimental"
fact that homotopy transfer formulae are given by sums over certain families of
{\em decorated graphs}. Moreover,
it follows that these sums describe essentially
a {\em  morphism}\, of coproperads $\cP^\Koz\rar B(\cP)$. This fact prompts one to think about the following
two closely related problems:
\Bi
\item[$\Rightarrow$] Given a quadratic (wheeled) properad, construct a quantum field theory
whose Feynman's perturbation series for the effective action
gives precisely the homotopy transfer formulae, i.e.\
a morphism of (wheeled) coproperads. This idea was first proposed
by A. Losev.
\item[$\Leftarrow$] Given a quantum field theory, find dg (wheeled) props
such that  Feynman's perturbation series for certain expectation values
 can be interpreted as their morphism.
\Ei
A simple and beautiful example where both problems have been successfully addressed
was constructed by Mnev in \cite{Mn}. Another example is studied in \S 6 of this paper.
A much less trivial example is given by the works of Kontsevich \cite{Ko} and
Cattaneo and Felder \cite{CF1} which imply
that the quantum Poisson sigma model on the
2-disk describes  a morphism of certain dg wheeled props
(see \cite{Me-lec} for their explicit construction).

\sip

In the rest of the paper we apply the above theory to a rather simple quadratic wheeled
properad, $\ULB$, of unimodular Lie 1-bialgebras. Remarkably, representations of the associated
dg wheeled prop,  $\ULB_\infty$, are in one-to-one correspondence with  so called {\em
(quasi-classically split) quantum BV manifolds},
interesting structures which one encounters in the Batalin-Vilkovisky quantization of certain gauge systems.

\mip
\bip

\section{Geometry of quantum Batalin-Vilkovisky manifolds}

{\bf 3.1. $\Z$-graded formal manifolds.} Batalin-Vilkovisky (shortly, BV) formalism \cite{BV}
is
one of most effective and universal methods
 for perturbative quantization of field theories with gauge
symmetries.
The first attempt
to understand the BV formalism as a {\em geometric}\, theory was done by Schwarz
who introduced and studied in  \cite{Sc} a category of so called
{\em $SP$-manifolds}\, to understand BV structures.
We adopt, however, in this paper
 a slightly different picture
of BV geometry based on semidensities and Khudaverdian's laplacian \cite{Kh}.
When one works in a fixed background the difference between these two pictures is not principal,
but we are going to concentrate in this section on {\em morphisms}\, and equivalences of BV structures,
and in this case the difference becomes decisive.

\sip

First we note that
\Bi
\item[(i)] ``manifolds", i.e.\ spaces of fields,
 used in the BV quantization are often {\em pointed}; the distinguished
point is called a {\em vacuum}\, state;
\item[(ii)] to make sense of perturbation series around the vacuum state  one is only
interested in a {\em formal}\, neighborhood of that state in the space of fields, not
in the global structure of the latter.
\Ei
Of course, one can try to ignore the formal nature of the perturbation series
and accept a genuine smooth supermanifold as a toy model for a space of fields.
However,
the {\em formal}\, nature of the basic notions and operations used in the BV formalism
resurrects  again  when one attempts  to make sense of expressions of the type
$e^{\frac{\Gamma(x, \hbar)}{\hbar}}$, where the function in the exponent,
$$
\Gamma(x, \hbar)=\Gamma_0(x) + \Gamma_1(x)\hbar + \ldots + \Gamma_n(x)\hbar^n + \ldots
$$
is a formal power series in a parameter (``Planck constant") $\hbar$. One can try to ignore this
issue as well,
and set $\hbar=1$. This is what is often done in many papers on geometric aspects
of the BV formalism.
We, however, can not afford setting the formal parameter to $1$ in the present paper
as without $\hbar$ no link
 between BV manifolds and  the homotopy theory of
 unimodular Lie 1-bialgebras
  holds true.
  Therefore right from the beginning we shall be working
 in the category of {\em formal}\, $\Z$-graded manifolds in which one can easily make
 a  coordinate independent sense to functions of the type
 $e^{\frac{\Gamma(x,\hbar)}{\hbar}}$ by demanding, for example,
 that $\Gamma_0(x)\in \cI$, where $\cI$ is the maximal ideal
 of the distinguished point. Let us give precise definitions.

\mip

The category of {\em finite-dimensional  $\Z$-graded formal manifolds}\,
 over a field $\K$  is,
by definition, the  category opposite to the category whose
\Bi
\item
objects are
 (isomorphism classes of)
of completed finitely generated free $\Z$-graded commutative $\K$-algebras;
every such a $\K$-algebra $\cR$ has a natural translation invariant
{\em adic}\, topology defined by the condition that
the powers, $\{\cI^n\}_{n\geq 1}$, of the maximal
ideal $\cI\subset \cR$ form a basis of open neighborhoods of $0\in \cR$;
\item morphisms are (isomorphism  classes of) continuous morphisms of topological $\K$-algebras.
\Ei

Thus, every $\Z$-graded formal manifold $\cM$ corresponds to a certain
isomorphism class, $\f_\cM$, of completed free finitely generated algebras
of the form $\K[[x^1, \ldots, x^n]]$, where formal variables $x^a$
(called {\em coordinates})
 are assigned
some degrees $|x^a|\in \Z$. The isomorphism class, $\f_\cM$, of $\K$-algebras
 is called the {\em structure sheaf}\,\footnote{We apologize for
 using the term {\em sheaf}\, in the present formal context as all the {\em sheaves}\,
 the reader encounters in the present section are rather primitive ---
 they are skyscrapers consisting of a single {\em stalk}\, over the distinguished point;
 this terminology helps, however,
 the geometric intuition (cf.\ \cite{Me-DG}).}  of the manifold $\cM$.
 A representation of $\f_\cM$ in
 the form
$\K[[x_1, \ldots, x_n]]$ is called a {\em coordinate chart}\, on $\cM$.
Such a representation  is not  canonical: a coordinate chart is defined up to
an arbitrary  (preserving $\Z$-grading)
change of coordinates  of the form
\Beq\label{Formal change of coordinates}
x^a \lon \hat{x^a} = \phi^a(x):=\sum_{k=1}^\infty \phi^a_{b_1\ldots b_k} x^{b_1}\ldots x^{b_k},
\ \mbox{for some}\  \phi^a_{b_1\ldots b_k}\in \K
\Eeq
where $\phi_{b_1}^a$ form an invertible matrix. Such changes form a group
of  formal diffeomorphisms, $\mathit D\mathit i\mathit f\mathit f(\cM)$,
and the Constitution of (formal) Geometry
says that every construction on a $\Z$-graded formal manifold
$\cM$ must be invariant under this group. There are, unfortunately,  not that many
$\mathit D\mathit i\mathit f\mathit f(\cM)$-invariant constructions possible in nature,
and their study is the major theme of {\em geometry}\, rather than {\em algebra}.
This is why we use
geometric terminology and intuition  throughout this section.

\sip

A {\em smooth map}, $\phi: \cM\rar \cN$, of formal graded manifolds is the same as a morphism,
$\phi^*: \f_\cN\rar \f_\cM$, of their structure sheaves. It is given in local coordinates  by formulae
of the type (\ref{Formal change of coordinates}) with $\phi_{b_1}^a$ not
 necessarily forming an invertible matrix.
A smooth invertible map is called a {\em diffeomorphism}.

\sip

The {\em tangent sheaf}, $\cT_\cM$,
of a $\Z$-graded manifold $\cM$ is, by definition, the $\Z$-graded $\f_\cM$-module
of derivations
 of the structure sheaf, that is, the module
of linear maps
$X: \f_\cM \rar \f_\cM$ satisfying the Leibnitz condition,
$X(fg)= (X(f))g + (-1)^{|X||f|}fX(g)$.
It is a  free $\f_\cM$-module generated, in a coordinate chart $\{x^a\}$, by partial derivatives,
$\p/\p x^a$. Elements of $\cT_\cM$ are called {\em smooth vector fields}\, on $\cM$. Every vector
field $X\in \cT_\cM$ is given in a coordinate chart as a linear combination,
$$
X=\sum_{a} X^a(x)\frac{\p}{\p x^a}, \ \ \ X^a(x)\in \K[[x^a]].
$$
This representation is not canonical: if $\{\hat{x}^a\}$ is another set of generators of $\f_\cM$
related to $\{x^a\}$ via (\ref{Formal change of coordinates}), then
$$
X=\sum_{a} X^a(x)\frac{\p}{\p x^a}= \sum_{b} \hat{X}^b(\hat{x})\frac{\p}{\p \hat{x}^b}
$$
with
$$
 \hat{X}^b(\hat{x})\mid_{\hat{x}^a=\phi^a(x)}= \sum_a X^a(x)\frac{\p \phi^b(x)}{\p x^a} .
$$
The matrix $({\p \phi^b(x)}/{\p x^a})$ is called the {\em Jacobian}\, of the coordinate
transformation (\ref{Formal change of coordinates}). The $\f_\cM$-module $\cT_\cM$
has a natural graded Lie algebra structure with respect to the ordinary graded commutator
of derivations, $[X_1, X_2]=X_1\circ X_2 - (-1)^{|X_1||X_2|} X_2\circ X_1$.
 The rank of $\cT_\cM$ is equal to the number of generators of the algebra
 $\f_{\cM}$ and is called the {\em dimension}\, of the
 graded manifold $\cM$.

\sip

Let $V$ be a finite dimensional $\Z$-graded vector space. One can associate to $V$
a $\Z$-graded formal manifold $\cV$ by defining $\f_\cV$ to be the isomorphism class
of the $\K$-algebra $\widehat{\odot^\bullet}V^*$, where $V^*:=\Hom(V, \K)$. The manifold
$\cV$ is said to be {\em modeled}\, on the graded vector space $V$.
Every formal $\Z$-graded manifold $\cM$ is modeled by some uniquely defined graded vector
space
$\cT_{\star\in \cM}:= (\cI/\cI^2)^*$ called the {\em tangent vector space at the distinguished point
$\star$
in $\cM$}. Note that every morphism of graded manifolds, $\phi: \cM\rar \cN$, gives rise
to a well-defined map, $d\phi_*: \cT_{\star\in \cM} \rar \cT_{\star\in \cN}$, of tangent vector spaces,
but, in general, {\em not}\,
 to a morphism, $d\phi: \cT_\cM \rar \cT_\cN$, of tangent sheaves. The latter is well-defined
if, for example, $\phi$ is an isomorphism.

\sip

Let
$\cT_\cM^*:= \Hom_{\f_\cM}(\cT_\cM, \f_\cM)$ be the dual $\f_\cM$-module, and let
$\Omega^1_\cM:=\cT_\cM^*[1]$  be the same $\f_\cM$-module $\cT_\cM^*$ but with shifted  grading.
The latter is called
the sheaf of {\em differential 1-forms}\, on the graded manifold $\cM$.
Note that the natural pairing,
$$
\Ba{rccc}
\langle\ , \ \rangle : &  \cT_\cM \times \Omega^1_\cM & \lon & \f_\cM\\
& X \ot \tau & \lon & \langle X, \tau\rangle
 \Ea
 $$
 has degree $1$.
There is a
 canonical degree -1 $\K$-linear morphism,
 $$
 \Ba{rccc}
 d: & \f_\cM & \lon & \Omega^1_\cM\\
 &    f & \lon & df,
 \Ea
 $$
defined, for  arbitrary vector field $X\in \cT_\cM$ by the equality $\langle X , df\rangle
=  X(f)$. It is clear that $\Omega^1_\cM$ is a free $\f_\cM$ module with
 a basis given, in some coordinate chart $\{x^a\}$,  by 1-forms $dx^a$, i.e.\ every
 $1$-form $\tau$ can be represented in this chart as a linear combination,
 $\tau= \sum_a dx^a \tau_a(x)$, for some $\tau_a(x)\in \K[[x^a]]$. We also have
 $df=\sum_a dx^a \p f/\p x^a$.

 \sip

The sheaf of graded commutative algebras,
$\Omega^\bullet_\cM:= \odot^\bullet_{\f_\cM} \Omega^1_\cM$, generated by 1-forms
 is called the
{\em De Rham}\, sheaf on $\cM$. Elements of $\Omega^k_M:= \odot^k \Omega^1_\cM$
are called differential $k$-forms.
The morphism $d: \Omega^0_\cM\rar \Omega^1_\cM$ extends naturally to a
morphism $d: \Omega^k_\cM \rar \Omega^{k+1}_\cM$ for any $k$ making
thereby $(\Omega^\bullet_\cM, d)$
into a sheaf of {\em differential}\,  algebras, i.e.\ satisfying
$d^2=0$ and $d(\tau_1\tau_2)= (d\tau_1)\tau_2 + (-1)^{|\tau_1|}\tau_1 d\tau_2$
for any $\tau_1, \tau_2\in \Omega_\cM$. In a local coordinate chart $\{x^a\}$ on $\cM$
we have an  isomorphism
$$
\Omega^\bullet_\cM \simeq \K[[x^a, dx^a]], \ \ \  |dx^a|= |x^a|-1,
$$
with the de Rham differential given on generators by
$d(x^a):= dx^a$, $d(dx^a):=0$.


\mip

{\bf 3.2. Odd Poisson structure.}  Let $\cM$ be a formal $\Z$-graded manifold.
 A {\em odd Poisson structure}\,
on $\cM$  is a degree -1 linear map,
$$
\Ba{rccc}
\{\ \bullet \ \}: & \f_\cM \ot_\K \f_\cM & \lon & \f_\cM\\
& f\ot g & \lon & \ \{f\bullet g\},
\Ea
$$
such that $\{f\bullet g\}= (-1)^{fg + f + g}\{g\bullet f\}$ and
\Beqr
\label{Jacobi-bullet}
\{f\bullet \{g\bullet h\}\} &=& \{\{f\bullet g\}\bullet h\} + (-1)^{(|f|+1)(|g|+1}\{g\bullet \{f\bullet h\}\\
\{f\bullet g h\} &=&
\label{Poisson-bullet} \{f\bullet g\} h + (-1)^{fg +g}g\{f\bullet h\}.
\Eeqr
for any $f,g,h\in \f_\cM$. Thus  brackets $\{\ \bullet \ \}$ and the ordinary product of functions
make the structure sheaf
$\f_\cM$ into a sheaf
of so called {\em Gerstenhaber algebras}.

\sip

A $\Z$-graded formal manifold $\cM$ equipped with a degree -1 Poisson structure is
called an {\em odd Poisson
manifold}. A {\em Poisson map}, $\phi: (\cM_1, \{\ \bullet\ \}) \rar (\cM_2, \{\ \bullet\ \})$,
of odd Poisson manifolds
is a degree $0$ smooth map $\phi: \cM_1\rar\cM_2$ such that
$$
\{\phi^*(f)\bullet \phi^*(g)\}= \phi^*\{f\bullet g\}
$$
for any $f,g\in \f_{\cM_2}$.  If one translates
 brackets $\{\ \bullet \ \}$
from $\f_\cM$ to its ``shifted" version, $\f_\cM[1]$ via the natural isomorphisms
$\f_\cM\rightleftarrows \f_\cM[1]$, one obtains an ordinary $\Z$-graded Lie algebra structure
on $\f_\cM[1]$.
\sip

An important example of an odd Poisson structure
comes from the sheaf of polyvector fields defined next.

\mip

{\bf 3.3. Polyvector fields.} For any $\Z$-graded $n$-dimensional  formal manifold $M$ the
completed
graded commutative algebra $\f_{\cM}:=\widehat{\odot^\bullet} (\cT_M[-1])$ is free of
rank $2n$ and hence defines a $\Z$-graded formal manifold $\cM$ which is often called
{\em the total space, $\cM:= \Omega^1_M$, of the  bundle of 1-forms on $M$}.
Elements of its structure sheaf $\f_\cM$ are called {\em polyvector fields} on the manifold $M$
and the structure sheaf itself is often denoted by $\Poly(M)$. One sets
$\Poly^k(M):= \odot^k(\cT_M[-1])$ and call its elements $k$-{\em vector fields}.
 This terminology for $\cM$ and its structure sheaf
originates from the duality $\cT_M[-1]= \Hom_{\f_M}(\Omega^1_M, \f_M)$ and from
the natural inclusion of the degree shifted\footnote{
To avoid such a degree shifting the sheaf of polyvector fields
$\Poly(M)$ is defined by some authors  as $\f_{\cM}[-1]$.}
 tangent sheaf $\cT_M[-1] \subset \Poly(M)$.

\sip

A  coordinate
chart $\{x^a\}$ on $M$ induces a  coordinate chart
$$
\left\{x^a, \psi_a:= \Pi \frac{\p}{\p x^a},\  |\psi_a|= -|x^a| + 1\right\}
$$
on $\cM$, where $\Pi: \cT_M \rar \cT_M[-1]$ is the natural isomorphism.
A change of coordinates (\ref{Formal change of coordinates}) on $M$ induces a change
of coordinates,
\Beq\label{Total cotangent Glueing map}
\Ba{ccl}
x^a &\lon & \hat{x}^a= \phi^a(x)\\
\psi_a &  \lon & \hat{\psi}_a= \sum_b \frac{\p \phi^a(x)}{\p x^b} \psi_b
\Ea
\Eeq
on $\cM$. In these coordinates we have an isomorphism,
$$
\f_\cM\equiv \Poly(M)\simeq \K[[x^a, \psi_a]].
$$
It is not hard to check that the degree -1 brackets on $\f_\cM$ defined in such a coordinate
chart by
\Beq\label{Poisson_structure_on TM[1]}
\{f\bullet g\}:= \sum_a \left((-1)^{|f| |x^a|}\frac{\p f}{\p x^a} \frac{\p g}{\p \psi_a}
+  (-1)^{|f| (|x^a|+1)}\frac{\p f}{\p \psi_a} \frac{\p g}{\p x^a}.
\right)
\Eeq
satisfy the axioms (\ref{Jacobi-bullet}) and (\ref{Poisson-bullet}), and, moreover,
are invariant under transformations (\ref{Total cotangent Glueing map}). Hence they
define an odd Poisson structure on the manifold $\cM$. Brackets
(\ref{Poisson_structure_on TM[1]})
 on $\Poly(M)$ are often denoted by $[\ \bullet\ ]_S$
and called {\em Schouten}\, brackets.

\sip

Schouten
brackets $[\ \bullet\ ]_S$ restricted to the subsheaf  $\cT_M[-1]\subset \Poly(M)$ give,
modulo the degree shifting, the ordinary commutator of vector fields.

\mip

{\bf 3.4. Odd symplectic structures.}
Any odd Poisson structure on a graded formal manifold $\cM$ defines a homogeneous (of degree 1)
section, $\nu$, of the bundle $\Poly^2(\cM)$ by the
 formula,
$$
\langle \nu, df
dg\rangle= \{f\bullet g\}\ \ \ \forall f,g\in \f_\cM.
$$
where $\langle\ ,\  \rangle$ stands for the natural duality pairing between
$\Poly^2(\cM)=\Hom_{\f_\cM}(\Omega^2_\cM, \f_\cM)$
and $\Omega^2_\cM$.
An odd Poisson structure on $\cM$ is called {\em non-degenerate}\, or {\em odd symplectic}\,
if the associated 2-vector
field is non-degenerate in the sense that the
induced ``raising of indices" morphism of sheaves,
$$
\Omega^1_\cM\stackrel{\lrcorner\nu}{\lon}\cT_\cM ,
$$
is an isomorphism. The inverse map gives the rise to a degree $-1$ differential 2-form,
$\omega:= ``\nu^{-1}"$, on $\cM$ which satisfies, due to the Jacobi identity
(\ref{Jacobi-bullet}), the condition $d\omega=0$.

\mip

{\bf 3.4.1. Darboux lemma (see, e.g., \cite{Kh, Le, Sc})}. {\em Any $\Z$-graded manifold
 with a non-degenerate odd Poisson structure
is locally isomorphic to an odd Poisson manifold $\cM$ described in \S 3.3.}

\mip

Thus any odd symplectic manifold $\cM$ admits local coordinates,
$
\{(x^a, \psi_a), |\psi_a|= -|x^a| +1\}
$,
 in which
the odd Poisson brackets are given by (\ref{Poisson_structure_on TM[1]}).
The associated symplectic 2-form is then given by
$
\omega=\sum_a dx^a d\psi_a$.
These coordinates are called {\em Darboux coordinates}.

\mip


\
{\bf 3.4.2. Symplectomorphisms and canonical transformations.}
By Lemma~3.4.1, any odd symplectic manifold can be covered by a Darboux coordinate chart
$(x^a, \psi_a)$. For future reference we note that a generic change of coordinates
\Beq\label{Glueing map}
\Ba{ccc}
x^a &\lon & \hat{x}^a= \phi^a(x, \psi)\\
\psi_a &  \lon & \hat{\psi}_a= \phi_a(x, \psi)
\Ea
\Eeq
defines a new  Darboux coordinate
chart $(\hat{x}^a, \hat{\psi}_a)$ if and only if the equations
\Beq\label{Canonical transform in Darboux}
\Ba{c}
\displaystyle
\sum_a(-1)^{|x^a|(|x^c|+1)}\frac{\p \phi^a}{\p x^b}\frac{\p \phi_a}{\p x^c}=0, \ \
\sum_a(-1)^{|x^a||x^c|}\frac{\p \phi^a}{\p \psi_b}\frac{\p \phi_a}{\p \psi_c}=0, \\
\displaystyle
 \sum_a\left((-1)^{|x^a||x^c|}\frac{\p \phi^a}{\p x^b}\frac{\p \phi_a}{\p \psi_c} +
(-1)^{(|x^a|+1)(|x^b|+1)}\frac{\p \phi^a}{\p \psi_c}\frac{\p \phi_a}{\p x^b}\right)
=\delta_b^c:=\left\{\Ba{rr} 1 & \mbox{if}\ b=c,\\
0  & \mbox{if}\ b\neq c.\Ea\right.
\Ea
\Eeq
are satisfied.
A diffeomorphism (\ref{Glueing map}) satisfying equations (\ref{Canonical transform in Darboux})
is called a {\em canonical transformation}. It is easy to check that (\ref{Total cotangent Glueing map})
is a canonical transformation for arbitrary functions $\phi^a(x)$ which have the
associated Jacobi matrix $\p \phi^a(x)/\p x^b|_{x=0}$
invertible.

\sip

 A Poisson  diffeomorphism of odd symplectic manifolds,
 $\phi: (\cM, \om) \rar (\hat{\cM}, \hat{\om})$, is called a
 {\em symplectomorphism}.
 This is the same as a  diffeomorphism
$\phi: \cM\rar \hat{\cM}$ of smooth $\Z$-graded manifolds such that
$\phi^*(\hat{\omega})=\omega$. It is the assumption on the non-degeneracy of the
odd symplectic forms  which forces one to define symplectomorphisms as
 special cases of {\em diffeomorphisms}. In Darboux coordinate charts, $(x^a, \psi_a)$ and $(\hat{x}^a,
\hat{\psi}_a)$,
 on $\cM$ and, respectively,
 $\hat{\cM}$ any  symplectomorphism is given by functions (\ref{Glueing map})
 satisfying equations (\ref{Canonical transform in Darboux}); for that reason
 a  symplectomorphism is also often called a
 canonical transformation.

 \mip

 {\bf 3.4.3. Remark.}
 The collection of
  Poisson morphisms of odd symplectic manifolds is much richer than the collection
  of symplectomorphisms. For example,
 if $\cM_1$ and $\cM_2$ are odd symplectic manifolds, then $\cM_1\times \cM_2$
 is naturally an odd symplectic manifold and the natural projection
 $ \cM_1\times \cM_2\rar \cM_1$  is a well defined
 Poisson morphism, which is {\em not}\, a symplectomorphism.

\mip
{\bf 3.4.4. Hamiltonian vector fields}.
For any function $\Phi\in \f_\cM$ on an odd Poisson manifold $\cM$, the associated map
$$
\Ba{rccc}
H_\Phi: & \f_M & \lon & \f_M \\
      &   g  & \lon & \{\Phi\bullet g\}
\Ea
$$
is a derivation of the structure ring $\f_\cM$ and hence is a smooth vector field on
$\cM$ called the {\em
Hamiltonian vector field associated with a function $\Phi$}.
It is not hard to check that $[H_{\Phi_1}, H_{\Phi_2}]= H_{\{\Phi_1\bullet \Phi_2\}}$
 for any $\Phi_1,\Phi_2\in \f_M$.

\sip

If the Poisson structure is non-degenerate, then
in a local Darboux coordinate chart one has,
$$
H_\Phi=\sum_a \left((-1)^{|\Phi| |x^a|}\frac{\p \Phi}{\p x^a} \frac{\p }{\p \psi_a}
+  (-1)^{\Phi| (|x^a|+1)}\frac{\p \Phi}{\p \psi_a} \frac{\p }{\p x^a}
\right).
$$
Note that if the function $\Phi(x^a,\psi_a)$ has degree $1$, then the associated hamiltonian vector
field $H_\Phi$ has degree zero, and it makes sense to consider a system of
ordinary  differential equations,
\Beq\label{Integral curves}
\Ba{ccl}
\frac{d \phi^a(x,\psi, t)}{d t} &=& - (-1)^{|x^a|}
\frac{\p \Phi(\hat{x},\hat{\psi})}{\p \hat{\psi_a}}\mid_{
\hat{x}^a=\phi^a(x,\psi, t),\,  \hat{\psi}_a= \phi_a(x,\psi, t)} \\
\frac{d \phi_a(x,\psi, t)}{d t} &=&  (-1)^{|x^a|}
\frac{\p \Phi(\hat{x},\hat{\psi})}{\p \hat{x^a}}\mid_{
\hat{x}^a=\phi^a(x,\psi, t),\,  \hat{\psi}_a= \phi_a(x,\psi, t)} \\
\phi^a(x,\psi, t)|_{t=0}&=& x^a,\\
\phi_a(x,\psi, t)|_{t=0}&=& \psi_a,
\Ea
\Eeq
for the unknown functions $\phi^a(x,\psi, t)$ and $\phi_a(x,\psi, t)$ of degrees $|x^a|$
and, respectively, $|\psi_a|$. Moreover,
a classical theorem from the theory of systems of ordinary fifferential
equations guarantees that,
for a sufficiently small strictly positive $\var\in \R$,
its solution,
$$
\{\phi_t^a=\phi^a(x,\psi, t), \phi_{t\, a}=\phi_a(x,\psi, t)\},
$$
exists and is unique for all $t$ in the interval $[0, \var)$. Moreover, the solution
is real analytic with respect to the parameter $t$.
Using the above equations it is easy to check that
$$
\frac{d}{dt}
\sum_a(-1)^{|x^a|(|x^c|+1)}
\frac{\p \phi_t^a}{\p x^b}\frac{\p \phi_{t\, a}}{\p x^c}=0, \ \
\frac{d}{dt}\sum_a (-1)^{|x^a||x^c|}
\frac{\p \phi^a}{\p \psi_b}\frac{\p \phi_a}{\p \psi_c}=0,
$$
$$
\frac{d}{dt}\sum_a
\left((-1)^{|x^a||x^c|}\frac{\p \phi^a}{\p x^b}\frac{\p \phi_a}{\p \psi_c} +
(-1)^{(|x^a|+1)(|x^b|+1)}\frac{\p \phi^a}{\p \psi_c}\frac{\p \phi_a}{\p x^b}\right)
=0
$$
implying, in view of the boundary $t=0$ conditions on $\phi_t$, that, for any
$t\in [0, \var)$ the map
\Beqrn
x^a &\lon & \hat{x}^a = \phi^a(x, \psi, t)\\
\psi_a &\lon & \hat{\psi}_a = \phi_a(x, \psi, t)
\Eeqrn
satisfies equation (\ref{Canonical transform in Darboux}) and hence defines
a canonical transformations. Thus any degree 1 function $\Phi$ on $\cM$ gives naturally
 rise to a 1-parameter family of local symplectomorphisms $\phi_t:\cM \rar \cM$.

\mip

 {\bf 3.4.5. Lagrangian submanifolds.} Let $\cM$ be a $\Z$-graded manifold and $I\subset
 \f_\cM$ an ideal such that the quotient ring $\f_\cS\subset
 \f_\cM/I$ is free; this ring corresponds, therefore,
 to a $\Z$-graded manifold $\cS$ which is called a {\em submanifold}\, of $\cM$;
 the natural epimorphism $\f_\cM\rar \f_\cS$  is called  an {\em embedding}\,
 $\cS\hook \cM$.

 \sip

 Lemma 3.4.1 implies that any odd symplectic
 manifold $(\cM, \omega)$ is even dimensional, say $\dim \cM=2n$.
 An $n$-dimensional submanifold
 $\caL\hook \cM$ is called {\em Lagrangian}\, if $\omega|_{\caL}=0$, i.e.\
 the induced map $\omega|_\caL: \odot^2(\cT_\caL[-1]) \rar \f_\caL$ is zero.
 The {\em normal}
 sheaf, $N_{\caL|\cM}$, of the submanifold $\caL\hook \cM$ is defined by the
 short exact sequence of sheaves of $\f_\caL$-modules,
$$
 0\lon \cT_\caL \stackrel{i}{\lon}
  \cT_\cM\mid_\caL \lon N_{\caL|\cM} \lon 0
$$
 so that its dualization (and degree shifting) gives,
 $$
 0 \lon N_{\caL|\cM}^*[1] \lon  \Omega^1_\cM\mid_\caL \stackrel{p}{\lon}
  \Omega^1_\caL \lon 0.
 $$
The odd symplectic form $\om$ provides us with a {\em degree $0$}\,
 isomorphism
of the middle terms of the short exact sequences above,
$$
\cT_\cM\mid_{\caL} \stackrel{\lrcorner\omega}{\lon} \Omega^1_\cM.
$$
The condition $\omega|_{\caL}=0$ is equivalent to saying that the composition,
$$
  \cT_\caL\stackrel{i}{\lon}
  \cT_\cM\mid_\caL \stackrel{\lrcorner\omega}{\lon} \Omega^1_\cM  \stackrel{p}{\lon}
  \Omega^1_\caL,
  $$
 vanishes. Hence we get a canonical monomorphism of sheaves,
 $$
 \lrcorner \omega \circ i: \cT_\caL \lon N_{\caL|\cM}^*[1],
 $$
 which is an isomorphism because  both sheaves have the same rank as locally free
 $\f_\caL$-modules. Hence, for any Lagrangian submanifold $\caL\hook \cM$,
 $N_{\caL|\cM}= (\cT_\caL)^*[1]=\Omega^1_\caL$ and
 there is a canonically associated exact sequence,
 \Beq\label{Normal-bundle_Lagr}
 O\lon \cT_\caL \stackrel{i}{\lon}
  \cT_\cM\mid_\caL \lon \Omega^1_\caL \lon 0
 \Eeq
of sheaves.

\mip
{\bf 3.5. Densities and semidensities}. If $V$ is a $\Z$-graded  free  module over
a $\Z$-graded ring $R$, then $\cB er(V)$ is, by definition, a degree $0$ rank 1
free module over $R$ equipped with a distinguished family of bases $\{D_e\}$
defined as follows \cite{Ma} (see also \cite{Ca3}):
\Bi
\item[(i)]
for any base $e=\{e_\al\}$ of the module $V$, there is an associated
basis vector, $D_e$, of $\cB er(V)$;
\item[(ii)]
if $e=\{e_\al\}$ and $\hat{e}=\{\hat{e}_\be\}$ are
two bases of $V$ with the relation $\hat{e}_\be= \sum_\al e_\al A_\be^\al$ for some
non-degenerate matrix $A_\al^\be\in R$, then $D_{\hat{e}}= Ber(A)D_e$, where
$Ber(A)$ is the Berezinian  of the matrix $A$.
\Ei

If $\cM$ is a $\Z$-graded manifold, then
$$
\cB er(\cM):= \cB er(\cT_\cM^*)= \left(\cB er(\Omega^1_\cM)\right)^*
$$
is a rank 1 locally
free sheaf of $\f_\cM$-modules. Its elements which do not vanish
 at the distinguished point, are called
{\em densities}\, or {\em volume forms}\, on the manifold $\cM$.
Let, for concreteness, $\cM$ be an odd symplectic manifold. To every Darboux coordinate chart
$(x^b,\psi_b)$ on $\cM$ there corresponds, by definition,
a basis section, $D_{x,\psi}$, of $\cB er(\cM)$;
if $(x^b,\psi_b)$ and $(\hat{x}^a,\hat{\psi}_a)$ are two Darboux coordinate charts
related to each other by a canonical transformation (\ref{Glueing map}), then
$$
D_{\hat{x},\hat{\psi}}=
Ber\left(\frac{\p(\hat{x},\hat{\psi})}{\p (x,\psi)}\right)D_{x,\psi},
$$
where $\frac{\p(\hat{x},\hat{\psi})}{\p (x,\psi)}$ stands for the Jacobi matrix of
the natural
transformation (\ref{Glueing map}).

\sip

One can show (see, e.g., \cite{Sc,KV}) that for any odd symplectic manifold $\cM$
 the sheaf $\cB er(\cM)$ admits a square root, that is, there exists a sheaf $\cB er^{1/2}(\cM)$
of so called {\em semidensities}\, such that
\Beq\label{sqroot of berezinian sheaf}
\cB er(\cM)= \left(\cB er^{1/2}(\cM) \right)^{\ot 2}.
\Eeq
An element $\Theta$ of the $\f_\cM$-module
 $\cB er^{1/2}(\cM)$ which does not vanish at the distinguished point,
 is called a {\em semidensity}\,
on $\cM$. In a  Darboux coordinate chart on $\cM$
  a semidensity, $\Theta$,
can be represented in the form
$\Theta= \Theta_{x,\psi}\sqrt{D}_{x,\psi}$ for some degree zero formal power series
 $\Theta_{x,\psi}\in \K[[x^a, \psi_a]]$ with  $ \Theta_{x,\psi}|_{x=0, \psi=0}\in \K^*$.
Under a canonical transformation this representation changes as follows,
\Beq\label{Canonical_transform_of semidensities}
\Theta_{\hat{x},\hat{\psi}}=
\left(Ber\left(\frac{\p(\hat{x},\hat{\psi})}{\p (x,\psi)}\right)\right)^{-1/2}
 \Theta_{x,\psi}
\Eeq

\mip

{\bf 3.5.1. Odd Laplacian on semidensities.} Let $(\cM, \omega)$ be an odd symplectic manifold.
The odd symplectic structure on $\cM$ gives canonically rise to a differential
operator on semidensities,
$$
\Ba{rccc}
\Delta_\om: & \cB er^{1/2}(\cM) & \lon & \cB er^{1/2}(\cM)\\
 & \Theta & \lon & \Delta_\om\Theta,
\Ea
$$
defined in an arbitrary Darboux coordinate system  as follows \cite{Kh},
$$
\Delta_\om \Theta := \left(\sum_a\frac{\p^2 \Theta_{x,\psi}}{\p x^a \p\psi_a}\right)
\sqrt{D}_{x,\om}.
$$
A remarkable fact is that $\Delta_\om $ is well-defined, i.e.\ does not depend on
a particular choice of Darboux coordinates used in the definition, as under
arbitrary canonical transformations
(\ref{Glueing map}) one has \cite{Kh},
\Beq\label{Khud_formula}
\sum_a\frac{\p^2 \Theta_{\hat{x},\hat{\psi}}}{\p \hat{x}^a \p\hat{\psi}_a}
\mid_{\hat{x}^a=\phi^a(x,\psi) \atop \hat{\psi}_a=\phi_a(x,\psi)}
= \left(Ber\left(\frac{\p(\hat{x},\hat{\psi})}{\p (x,\psi)}\right)\right)^{-1/2}
\sum_a\frac{\p^2 \Theta_{{x},{\psi}}}{\p {x}^a \p{\psi}_a}.
\Eeq
The operator $\Delta_\om$ is uniquely determined by the underlying odd symplectic
structure and is called the {\em odd Laplacian}. This is an odd analogue
of the modular vector field in ordinary  Poisson geometry \cite{We}.
Its invariant definition
can be found in \cite{Se}; that definition is a bit tricky and involves a beautiful
Manin's
description of the Berezinian $\cB er(\cM)$  as a cohomology class in a
certain complex (see Chapter 3, \S 4.7 in \cite{Ma}).

\sip

For an arbitrary Darboux coordinate chart $(x^a, \psi_a)$
the second order operator  $\sum_a\frac{\p^2 }{\p {x}^a \p{\psi}_a}$
is denoted from now on by $\Delta_{x,\psi}$ or, when a particular choice of  Darboux
coordinates
is implicitly assumed, simply by $\Delta_0$. This operator has an invariant meaning
only when applied to (coordinate representatives of) semidensities, not to
ordinary functions.

\mip

{\bf 3.5.2. Lemma.} $\Delta_\om^2=0$.

\sip

Proof is evident when one uses  Darboux coordinates.

\mip

{\bf 3.6. Batalin-Vilkovisky manifolds.} A {\em Batalin-Vilkovisky structure}\,
(or, shortly, {\em BV-structure})
 on an odd symplectic manifold $(\cM, \omega)$ is a semidensity
 $\Theta\in \cB er^{1/2}(\cM)$   satisfying an equation,
 \Beq\label{Master equation invariant def}
 \Delta_{\om}\Theta =0.
\Eeq
Equation (\ref{Master equation invariant def}) is called a
{\em master equation}, while its solution $\Theta\in \cB er^{1/2}(\cM)$
a {\em master semidensity}.

\sip

Such structures first emerged in the powerful Batalin-Vilkovisky approach to the quantization
of field theories with gauge symmetries (see, e.g., \cite{BV, Sc,Ca1,Ca3,Ca2} and
references cited there). A concrete example of such a BV quantization
machine is considered in \S 6 below.

\sip

The automorphism group, $Aut(\cM,\om)$, of the odd symplectic manifold
(that is, the group of symplectomorphisms $(\cM, \om) \rar (\cM, \om)$) acts naturally on
the set of BV-structures on $(\cM,\om)$: if $\Theta$ is a master semidensity,
then for any $\phi\in  Aut(\cM,\om)$, its pullback $\phi^*(\Theta)$ defined by
(\ref{Canonical_transform_of semidensities}) is again a
master semidensity.

\mip

{\bf 3.6.1. Remark.} Following Schwarz \cite{Sc}, BV-structures on odd symplectic maniolds
are defined in many papers in a different way: one first fixes an extra structure,
a
volume form $\rho$ on $\cM$, and then one defines an odd Laplacian, $\Delta_{\om,\rho}$,
 on {\em functions} as a map
$$
\Ba{rccc}
\Delta_{\om, \rho}: & \f_\cM & \lon & \f_\cM\\
&                   f   & \lon & \frac{\caL_{H_f}\rho}{\rho},
\Ea
$$
where $\caL_{H_f}$ stands for the Lie derivative along the Hamiltonian vector field
$H_f$ associated with a function $f\in \f_\cM$. The data $(M,\om, \rho)$ is called
in \cite{Sc} an $SP$-{\em manifold}, and a BV-structure on an $SP$-manifold
is  defined as a function $f\in \f_\cM$ satisfying the equation $\Delta_{\om, \rho}f=0$.
In fact, the volume form $\rho$ can not be arbitrary but must satisfy
en extra condition \cite{Sc} which assures that in some Darboux coordinate
the equation for $f$ takes the form $\Delta_{x,\psi}f=0$ making it completely equivalent
to the above semidensity approach via an association
$$
 f \rightleftarrows \frac{\Theta}{\sqrt{\rho}}.
 $$
The $SP$-manifold approach to the BV-geometry does not seem to be a natural one for the following two reasons:
\Bi
\item[(i)] it is often accompanied with an unduly restriction of the gauge group of  the set of BV-structures on $\cM$
from arbitrary symplectomorphisms to {\em volume preserving}\, symplectomorphisms
(in contrast to ordinary symplectic geometry in {\em odd}\, symplectic geometry a generic
symplectomorphism is {\em not}\, necessarily volume preserving);
\item[(ii)] in applications of the BV formalism to
quantizations
one never integrates over the ``phase space" $\cM$ itself but rather over its Lagrangian
submanifolds $\caL\hook \cM$ (see \S 6 for a concrete example) depending on
a gauge fixing. Thus what one  needs in applications
 is not a volume form on
$\cM$ but rather a global object on that odd symplectic manifold
 which restricts to a volume form on its arbitrary Lagrangian submanifold
$\caL \hook \cM$.  Extension
(\ref{Normal-bundle_Lagr}) implies,
\Beq\label{Semidensity_on_Lagrangian}
\cB er(\cM)\mid_\caL= \cB er (\caL)\ot \cB er(\Omega^1_\caL)^*= \Ber (\caL)^{\ot 2},
\Eeq
which in turn implies that it is an appropriately chosen semidensity on $\cM$
(rather than a volume
form on $\cM$) which might restrict to a volume form on the Lagrangian submanifold.
\Ei

In fact we have no choice as to adopt a definition of BV structures via semidensities
rather than via $SP$-manifolds as in our approach Definition~3.6
(as well as definition of a morphism of $BV$-manifolds, see \S 3.9 below) {\em follows}\,
from the homotopy theory of Lie 1-bialgebras and the associated
homotopy transfer formulae.

 \mip

{\bf 3.6.2. Definition.}
A data $(\cM, \om, \Theta)$ consisting of an odd symplectic manifold $(\cM, \om)$
and a master semidensity $\Theta\in \cB er^{1/2}(\cM)$
is called a {\em Batalin-Vilkovisky manifold},
or simply a {\em BV-manifold}.

\mip

{\bf 3.6.3. Dilation group action and pointed manifolds.}
If $\Theta$ is a BV structure on an odd symplectic manifold
$\cM$ then $\lambda\Theta$ is again a BV structure for any non-zero constant
$\lambda\in \K$. From now one we identify such BV structures, i.e.\  we understand
a master semidensity $\Theta$
as an element of the projective space $\P\cB er^{1/2}(\cM)$.

\sip

Formal $\Z$-graded manifolds $\cM$ are always
{\em pointed}, i.e.\ have a distinguished point $* \in \cM$ corresponding to the unique maximal ideal
in $\f_\cM$ which is often called (in the quantization context)
{\em the vacuum state}. In a Darboux coordinate system $(x^a, \psi_a)$
centered at $*$ one can always normalize a master semidensity
$\Theta=\Theta_{x,\psi} D_{x,\psi}\in \P\cB er^{1/2}(\cM)$
in such a way that $\Theta_{x,\psi}|_{x=\psi=0}=1\in \K$, and this normalization is invariant under
formal canonical transformations. It is often suitable to represent
such a normalized semidensity in the form
$\Theta=e^{\Gamma(x,\psi)} \sqrt{D}_{x,\psi}$ for some
smooth formal function $\Gamma(x,\psi)$ {\em vanishing at zero}
(so that its exponent is well-defined as a formal power series);
the master equation takes then the form
$$
\Delta_\om \Theta= \left(\Delta_0\Gamma + \frac{1}{2}\{\Gamma\bullet \Gamma\}\right)
\Theta=0.
$$
As $\Theta$ is, by assumption, non-vanishing, the latter equation is equivalent to
\Beq\label{Master equation for Gamma}
\Delta_0\Gamma + \frac{1}{2}\{\Gamma\bullet \Gamma\}=0.
\Eeq
where $\{\ \bullet\ \}$ are the odd Poisson brackets on $\cM$.
The normalization $\Gamma|_*=0$ is assumed from now on.

\mip

{\bf 3.7. Sheaves of Gerstenhaber-Batalin-Vilkovisky (GBV) algebras.}
A $\Z$-graded commutative unital algebra $\cA$ equipped with a degree -1
 linear map
$\Delta: \cA \rar \cA$ satisfying
\Bi
\item[(i)] $\Delta^2=0$,
\item[(ii)] and, for any $a,b,c\in \cA$,
\Beqrn
\Delta(a b  c) & = & \Delta(a b)  c +
(-1)^{|b|(|a|+1)} b\Delta(a c) +
(-1)^{|a|}a \Delta(b c)\\
&& - \Delta(a) b  c - (-1)^{|a|} a
\Delta(b) c - (-1)^{|a|+|b|} a b
\Delta(c).
\Eeqrn
\Ei
is called a {\em GBV-algebra} (see, e.g., \cite{Ma2}).
 Note that $\Delta(1)=0$.
One can check  \cite{Ma2}  that the linear map
\Beq\label{GBV brackets}
\Ba{rccl}
[\ \bullet \ ] & : \cA\ot \cA & \lon & \cA \\
& a \ot b & \lon & [a\bullet b]:= (-1)^{|a|} \Delta(a b) -
(-1)^{|a|} \Delta(a)\circ b - a \Delta(b)
\Ea
\Eeq
makes $\cA$ into an odd Lie superalgebra, i.e.\ the  Jacobi identities
of type
(\ref{Jacobi-bullet}) are satisfied. Moreover, the odd Poisson identity
(\ref{Poisson-bullet}) also holds true,
$$
[a\bullet (b c)] = [a\bullet b]c +
(-1)^{|a|(|b|+1)}b [a\bullet c],
$$
 for any $a,b\in \cA$. The operator $\Delta$ is called a BV-operator of the
 GBV-algebra $\cA$.

\mip

{\bf 3.7.1. Lemma.} {\em  Let $(\cM, \om, \Theta)$ be a BV-manifold.
Then its structure sheaf is naturally a sheaf of GBV-algebras with
the BV operator given by}
$$
\Ba{rccl}
\Delta_{\om, \Theta}: & \f_\cM & \lon &   \f_\cM \\
& f & \lon & \Delta_{\om, \Theta}f:= \frac{ \Delta_\om(f\Theta)}{\Theta}.
\Ea
$$

\sip

\Proof Representing $\Theta$ in a local Darboux coordinate system as
$e^\Gamma \sqrt{D}_{x,\psi}$, we get
$$
\Delta_{\om, \Theta} f = \frac{\Delta_0(fe^\Gamma) D_{x,\psi}}{\Theta}=
\frac{
\left(\Delta_0 f + \{\Gamma\bullet f\} \right)\Theta}{\Theta}=
\Delta_0 f + \{\Gamma\bullet f\}
$$
Then, for any $f\in \f_\cM$,
\Beqrn
(\Delta_{\om, \Theta})^2f &=& \Delta_{\om, \Theta}\left(\Delta_0 f +
\{\Gamma\bullet f\}\right)\\
&=&\Delta_0\left(\Delta_0 f +
\{\Gamma\bullet f\}\right) +\left\{\Gamma\bullet \left(\Delta_0 f +
\{\Gamma\bullet f\}\right)\right\}\\
&=& \left\{\left(\Delta_0\Gamma + \frac{1}{2}\{\Gamma \bullet \Gamma\}\right) \bullet
f \right\}\\
&=& 0,
\Eeqrn
so that condition (i) in the definition of a GBV algebra is satisfied. Condition
(ii) can be checked analogously.
\hfill $\Box$

\mip

It is easy to see that the odd Lie brackets induced on the structure sheaf
$\f_\cM$ by formula (\ref{GBV brackets}) coincide precisely with the
Poisson brackets of the underlying odd symplectic structure.

 \mip

 {\bf 3.8. Quantum master equation.} Let $\hbar$ be a formal parameter of degree
 $2$ and let $\K[[\hbar]]:= \{\sum_{n\geq 0} a_n\hbar^n, a_n\in \R\}$
 be the associated graded commutative ring of formal power series. The latter defines a $\Z$-graded
 manifold which we denote by $\K[[\hbar]]^{\mathsf v}$.
 We are interested in considering $\hbar$-twisted formal smooth manifolds $\cM^\hbar$
 whose structure sheaves $\f_{\cM^\hbar}$ are non-canonically isomorphic to
 $\f_\cM[[\hbar]]:= \f_\cM\ot_\R \R[[\hbar]]$ where $\f_\cM$ is the structure sheaf of
 some $\Z$-graded smooth manifold $\cM$. Such an $\hbar$-twisted manifold $\cM^\hbar$
 is best understood as a formal {\em family}\, of manifolds, $\pi:\cM_\hbar\rar \K[[\hbar]]^{\mathsf v}$,
 over the 1-dimensional formal $\Z$-graded manifold $\K[[\hbar]]^{\mathsf v}$.
 The fiber, $\cM^0:=\pi^{-1}(\star)$, over the distinguished point $\star\in \K[[\hbar]]^{\mathsf v}$
 is a $\Z$-graded formal manifold called the {\em classical limit of}\, $\cM^\hbar$.

 \sip

 More precisely, let us consider a  category, $\cC_\hbar$,  whose objects are
 (isomorphism classes) of completed $\Z$-graded free
 $\K[[\hbar]]$-algebras; they are equipped with a natural monomorphism, $\pi^*: \K[[\hbar]] \rar \f_{\cM_\hbar}$,
 of $\K[[\hbar]]$-algebras; morphisms in this category are defined as {\em continuous}\,  morphisms
 of topological $\K[[\hbar]]$-algebras commuting with the monomorphism $\pi^*$. The quotient
 of an algebra $\f_{\cM^\hbar}$ by the ideal generated by $\hbar$ is denoted
 by $\f_{\cM^0}$; this is the structure sheaf of the classical limit  $\cM^0=\pi^{-1}(\star)$.

 \sip

 The opposite category $\cC_\hbar^{\mathsf v}$
 is called the category of $\hbar$-twisted manifolds.
 As in \S 3.1-3.4 on can define natural relative versions of all basic concepts ---
 tangent sheaves, De Rham sheaves, odd Poisson structures and odd symplectic structures.
 For example, an $\hbar$-twisted odd symplectic manifold can be defined as an equivalence
 class of Darboux coordinate charts
$(x^a, \psi_a)$ modulo canonical transformation of the form,
 \Beq\label{Glueing h-map}
\Ba{ccc}
x^a &\lon & \hat{x}^a= \phi^a(x, \psi, \hbar)\\
\psi_a &  \lon & \hat{\psi}_a= \phi_a(x, \psi, \hbar)
\Ea
\Eeq
where $\phi^a(x, \psi, \hbar)$ and $\phi_a(x, \psi, \hbar)$ are formal power series
from $\K[[x^a, \psi_a, \hbar]]$
such that equations (\ref{Canonical transform in Darboux}) hold
and
 the Jacobian $\frac{\p(\hat{x},\hat{\psi})}{\p (x,\psi)}$
gives an invertible matrix at the point $(x^a=0,\psi_a=0,\hbar=0)$.
\sip

Let $\cM_\hbar$ be an $\hbar$-twisted $\Z$-graded manifold. We need a singular
(with respect to $\hbar$) extension of its structure sheaf $\f_{\cM^\hbar}$.
Let us fix an isomorphism $i: \f_{\cM^\hbar}\simeq \f_{\cM^0}\ot_\K \K[[\hbar]]$, and
use it to extend $\f_{\cM^\hbar}$  as follows,
\Beq\label{h-ring}
\f_{\cM^{\hbar, \hbar^{-1}}}:=\left\{ \sum_{n= -\infty}^\infty f_n \hbar^{n}
\in \f_{\cM^0}\ot_\K \K[[\hbar,\hbar^{-1}]]:
f_{-n}\in
I^n\ \mbox{for}\ n \geq 1\right\}.
\Eeq
where $I$ is the maximal ideal in $\f_{\cM^0}$.
The resulting vector space has natural a $\K[[\hbar]]$-algebra structure extending that of  $\f_{\cM_\hbar}$;
moreover, the isomorphism class of this extension  is independent of a particular choice of a map $i$
used in the definition.
\sip

Let $\cM^\hbar$ be an $\hbar$-twisted odd symplectic manifold.
An invertible element $\Theta$ of the sheaf  $\cB er^{1/2}(\cM^\hbar)\ot_{\f_{\cM^\hbar}}
\f_{\cM^{\hbar, \hbar^{-1}}}$
is called {\em regular}\, if in some  Darboux coordinate charts, $(x^a, \psi_a)$ it can
be represented
in the form
$$
\Theta= e^{\frac{\Gamma}{\hbar}}\sqrt{D}_{x,\psi}
$$
for some  function $\Gamma(x,\psi,\hbar)\in \cI$ whose classical limit, $\Gamma|_{\hbar=0}$, lies in $I$.
(Here and elsewhere  $\cI$ stands for the maximal ideal
in the $\K$-algebra $\f_{\cM_\hbar}$, and $I$ for the maximal ideal of $\f_{\cM^0}$.) It clear that
this notion does not depend on the choice of a Darboux coordinate chart
used in the definition.

\mip

{\bf 3.8.1. Definition.} A {\em quantum Batalin-Vilkovisky structure}\,
 on an $\hbar$-twisted odd symplectic manifold $(\cM_\hbar, \omega)$ is a
 regular element  $\Theta \in \cB er^{1/2}(\cM^\hbar)\ot_{\f_{\cM^\hbar}}
\f_{\cM^{\hbar, \hbar^{-1}}}$
 satisfying an equation,
 \Beq\label{Quantum Master equation invariant def}
 \Delta_{\om}\Theta =0.
\Eeq
This equation is called a
{\em quantum master equation}, while its solution
$\Theta$
a {\em quantum master semidensity}.
In a  Darboux coordinate system the quantum master equation has the form,
\Beq\label{Quantum Master equation Darboux}
\hbar \Delta_0 \Gamma + \frac{1}{2}\{\Gamma\bullet \Gamma\}=0.
\Eeq
The structure sheaf $\f_{\cM^\hbar}$ can be made into a sheaf of
GBV algebras with respect to the operator,
\Beq
\Ba{rccl}\label{Ch3: odd laplacian on functions hbar}
\Delta_{\om, \Theta}: & \f_{\cM^\hbar} & \lon &   \f_{\cM^\hbar} \\
& f & \lon & \Delta_{\om, \Theta}f:= \hbar\frac{ \Delta_\om(f\Theta)}{\Theta}=
\hbar \Delta_0 f + \{\Gamma\bullet f\}.
\Ea
\Eeq

\mip

{\bf 3.8.2. Fact-definition.} Let  $\phi: (\cM^\hbar, \om) \rar (\hat{\cM}^\hbar, \hat{\om})$
 be a symplectomorphism
of odd symplectic manifolds, and $\Theta$ a quantum BV structure on $\cM^\hbar$. Then
$\hat{\Theta}:=(\phi^{-1})^*\Theta$ (with the pullback map $(\phi^{-1})^*$  given in local coordinates by (\ref{Canonical_transform_of semidensities}))  is a quantum BV structure on $\hat{\cM}^\hbar$; such a pair,
$(\cM^\hbar, \om, \Theta)$ and $(\hat{\cM}^\hbar, \hat{\om}, \hat{\Theta})$,
of quantum BV structures is called  {\em symplectomorphic}.

\mip

{\bf 3.9. Quantum BV manifolds.}
Let $V$ be a $\Z$-graded finite-dimensional vector space. Slightly abusing notations,
the formal $\hbar$-twisted $\Z$-graded
 manifold
corresponding to the isomorphism class of the $\K[[\hbar]]$-algebra
$\widehat{\odot^\bullet}(V\oplus V^*[-1])^*\ot_\K \K[[\hbar]]$ is denoted from now on by
$\cM^\hbar_{V}$ rather than by $\cM^\hbar_{V\oplus V^*[-1]}$;
it has a natural odd symplectic
structure induced by the pairing between $V$ and $V^*[-1]$.
Any $\hbar$-twisted odd symplectic manifold is isomorphic to $\cM^\hbar_V$
for some non-canonically defined vector space $V$. We shall  consider next
an extra structure --- a {\em ordered}\, pair of  transversal Lagrangian submanifolds
in the classical $\hbar\rar 0$ limit, $\cM_V^{0}$, of $\cM_V^\hbar$ --- which will make
the correspondence $\cM_V^\hbar \rightleftarrows V$ canonical.

\sip

The inclusions $V \subset V\oplus V^*[-1]$
and $V^*[-1]\subset V\oplus V^*[-1]$
correspond to two transversal  Lagrangian submanifolds in
$\cM_V^{0}$ which we denote by the symbols $\caL_V$ and, respectively,
 $\caL_{V}^\bot$ and consider from now as an extra  part
of the definition of an $\hbar$ twisted odd symplectic manifold
$\pi: \cM_V^\hbar\rar \K[[\hbar]]^{\mathsf v}$.
The automorphism group of $\cM_V^\hbar$ consists, therefore, of those symplectomorphisms
$\phi: \cM_V^\hbar\rar \cM_V^\hbar$ which leave Lagrangian submanifolds $\caL_V$ and $\caL_{V}^\bot$
in the fiber $\cM_V^0$ over $\hbar=0$
invariant. We can always find an {\em adopted}\, Darboux coordinate chart,
$(x^a, \psi_a)$ on $\cM_V^\hbar$ such that the Lagrangian submanifold $\caL_V\hook \cM_V^0$ is
given by the equations $\psi_a=0$ and  the Lagrangian submanifold $\caL_V^\bot\hook \cM_V^0$
 by the equations $x^a=0$. Then $Aut(\cM_V^\hbar)$ consists of canonical transformations
 (\ref{Glueing h-map}) satisfying the conditions,
 \Beq\label{Glueing h-map-Lagrang-V}
\phi^a(x,\psi, \hbar)\mid_{x=\hbar=0}=0, \ \ \ \phi_a(x,\psi, \hbar)\mid_{\psi=\hbar=0}=0.
\Eeq

Note that the vector space $V$ is canonically isomorphic
to the Lagrangian subspace $\cT_{\star\in \caL_{V}}$ of the tangent space
$\cT_{\star\in \cM_{V}}$ and hence has an invariant meaning.
Moreover, the tangent space $\cT_{\star\in \cM_{V}}$  is canonically decomposed
into a direct sum of Lagrangian subspaces,
$\cT_{\star\in \caL_{V}}\oplus  \cT_{\star\in \caL_{V}}^\bot$, so that the odd symplectic structure
on $\cM_V^0$ does indeed coincide with the one which is induced
from the natural parings  between $\cT_{\star\in \caL_{V}}$
and $ \cT_{\star\in \caL_{V}}^\bot=  \Hom_K(\cT_{\star\in \caL_{V}}, \K)[-1]$.
\mip

{\bf 3.9.1. Fact (cf.\ \cite{Sc})}. Any odd symplectic manifold $\cM^\hbar$ equipped
with an ordered pair of
transversal Lagrangian submanifolds $\caL_1$ and $\caL_2$ in $\cM^{0}$ is symplectomorphic
to the odd symplectic manifold $\cM_V^\hbar$ for some uniquely defined vector space
$V:= \cT_{\star\in \caL_{1}}$.

\mip

{\bf 3.9.2. Definition.} (i)  A {\em quantum BV manifold}\, is an $\hbar$-twisted odd symplectic
manifold $\cM_V^\hbar$ associated with some graded vector space $V$ equipped with a quantum
Batalin-Vilkovisky structure $\Theta$ such that in an adopted Darboux coordinate chart
 $(x^a, \psi_a)$ one has $\Theta= e^{\frac{\Gamma}{\hbar}}D_{x,\psi}$ with
 ``classical" and ``semiclassical" parts of
 $\Gamma(x, \psi, \hbar)=\sum_{k\geq 0}\Gamma_k(x,\psi)\hbar^k$ satisfying the boundary conditions,
\Beq\label{Boundary conditions for Gamma}
\Gamma_0(x, \psi)\in I_{V} I_{V^\bot}
 \ \mbox{and}\ \ \Gamma_1(x, \psi)\in I_{V} + I_{V^\bot},
\Eeq
where $I_V\simeq \psi \K[[x, \psi]]$ and $I_{V^\bot}\simeq x\K[[x, \psi]]$
are the ideals of the  Lagrangian
submanifolds $\caL_V$ and, respectively, $\caL_{V^\bot}$ in the classical limit $\cM_V^0$.

\mip

(ii) A {\em symplectomorphism}, $\phi: (\cM^\hbar_{V_1}, \Theta_1) \rar (\cM^\hbar_{V_2}, \Theta_2)$, of quantum
BV manifolds is a symplectomorphism $\phi$ of the associated quantum BV structures (see \S 3.8.2) which
 respects Lagrangian submanifolds in the fibre over $\hbar=0$, i.e.
 $\lim_{\hbar\rar 0}\phi(\caL_{V_1})\subset \caL_{V_2}$
and  $\lim_{\hbar\rar 0}\phi(\caL_{V_1}^\bot)\subset \caL_{V_2}^\bot$.

\mip

{\bf 3.9.3. Remarks.} {\bf (i)}
The second boundary condition, $\Gamma_1(x, \psi)\in I_{V} + I_{V^\bot}$,
 says only
that $\Gamma_1$ has no constant term, i.e.\ $\Gamma$ itself has no term proportional
to $\hbar$; as the quantum
master equation is invariant under translations $\Gamma\rar \Gamma + \K[[\hbar]]$,
the second boundary condition is only a partial normalization condition on the quantum master function.
The first boundary condition in Definition 3.9.2(i) is quite restrictive, but still allows many
interesting examples such as, e.g.,  $BF$ theory and its various generalizations
(see, e.g., \cite{CR, Mn} and also \S 6).

\sip

{\bf (ii)}
 In view of the presence of boundary conditions, the structure in \S 3.9.2 should be
more precisely  called  a {\em quantum BV manifold with split quasi-classical limit}.
We abbreviate it to simply a {\em quantum BV manifolds}\, in this paper.

\mip
{\bf 3.10. Homotopy classification of quantum BV manifolds.}
Let $(\cM_V^\hbar, \Theta)$ be a quantum BV manifold associated with a $\Z$-graded vector space $V$.
 In an adopted Darboux coordinate chart we have
$\Theta= e^{\frac{\Gamma}{\hbar}}\sqrt{D}_{x,\psi}$, where, in view of  boundary
conditions (\ref{Boundary conditions for Gamma}),  the formal power series must have the form
$$
\Gamma(x, \psi, \hbar)=
\underbrace{\sum_{a,b}\Gamma_{(0)\, b}^{\ \ \ a}x^b \psi_a}_{\Gamma_0}
+ \underbrace{
\sum_{n\geq 1, p+q+2n\geq 3\atop
p+n\geq 2,q+ n\geq 2}\frac{1}{p!q!}
\Gamma_{(n)\, a_1\ldots a_p}^{\ \ \ b_1\ldots b_q} x^{a_1}\ldots  x^{a_p}\psi_{b_1}\ldots
\psi_{b_q} \hbar^n}_{\bGa}
,
$$
for some  $\Gamma_{(n)\, a_1\ldots a_p}^{\ \ \ b_1\ldots b_q}\in \K$.
 Quantum master equation
 (\ref{Quantum Master equation Darboux}) immediately implies
 $$
 \{\Gamma_0 \bullet\Gamma_0\} =0,
 $$
or, equivalently,
$$
 \sum_c \Gamma_{(0)\, c}^{\ \ \ a} \Gamma_{(0)\, b}^{\ \ \ c}=0.
$$
The linear functions $x^a\bmod \cI^2_\hbar$, where $\cI_\hbar$ is the maximal ideal
in the $\K[[\hbar]]$-algebra $\K[[x^a, \psi_a, \hbar]]$, form a basis of the vector space
$V^*$; let $\{e_a\}$
be the associated dual basis of $V$, and define a degree 1 map
$$
\Ba{rccc}
d: & V & \lon & V \\
   & e_a & \lon & d(e_a):= \sum_c e_c\Gamma_{(0)\, a}^{\ \ c}
\Ea
$$
Clearly, $d^2=0$. Moreover, the map $d$ does not depend on the choice
of an adopted Darboux coordinate chart $(x^a, \psi_a)$ used in its definition as
the third equation in (\ref{Canonical transform in Darboux}) implies that under a generic
canonical transformation (\ref{Glueing h-map}),
$$
\Ba{ccl}
x^a &\lon & \hat{x}^a= \phi^a(x, \psi, \hbar)= \sum_b
\cA_b^a x^b + \mbox{other terms}, \ \ \ \cA_b^a\in \K,\\
\psi_a &  \lon & \hat{\psi}_a= \phi_a(x, \psi, \hbar)= \sum_b
\cB_a^b \psi_b + \mbox{other terms},
\ \ \ \cB_b^a\in \K,
\Ea
$$
the leading matrices $\cA$ and $\cB$ must be inverse to each other. Hence the differential
$d$ on the vector space $V$ is defined canonically
and is called {\em the differential induced by the master semidensity}\,
or simply {\em induced differential}. In this situation we say that the quantum BV manifold
$(\cM_V^\hbar, \Theta)$ {\em is modeled on a dg vector space}\,  $(V,d)$.
Setting $d:=\{\Ga_0\bullet\ldots \}$ we can rewrite the quantum master equation in the form,
\Beq\label{Quantum master eqn for Gamma}
d\bGa + \hbar\Delta_0\bGa + \frac{1}{2}\{\bGa\bullet \bGa\}=0.
\Eeq

\sip

{\bf 3.10.1. Definition.}  A quantum BV-manifold $(\cM_V^\hbar, \Theta)$ is called {\em minimal}\, if
in some (and hence any) adopted Darboux coordinate system $(x^a, \psi_a)$ one has
$\Theta= e^{\frac{\Gamma}{\hbar}}\sqrt{D}_{x,\psi}$  with the quadratic part,
$\Ga_0$, of $\Gamma$ vanishing, i.e. with $\Gamma= \bGa$.

\mip

{\bf 3.10.2. Definition.}
 A quantum BV-manifold $(\cM_V^\hbar, \Theta)$ is called {\em contractible}\, if
 the associated complex $(V,d)$ is acyclic and
there exists an adopted Darboux coordinate system $(x^a, \psi_a)$ in which
$\Theta= e^{\frac{\Gamma}{\hbar}}\sqrt{D}_{x,\psi}$  with $\Gamma=\Gamma_0$,
 i.e.\ with $\bGa=0$.

\mip

{\bf 3.10.3. Symplectomorphisms and morphisms of tangent complexes.}
A symplectomorphism, $\phi: (\cM_{V}^\hbar, \Theta) \rar (\cM_{\hat{V}}^\hbar, \hat{\Theta})$, of quantum BV
manifolds induces a linear map, $ \cT_{\star\in \cM_{V}^\hbar}\rar \cT_{\star\in \cM_{\hat{V}}^\hbar}$,
of tangent spaces at the distinguished points and, as $\phi(\caL_{V})\subset \caL_{\hat{V}}$
and $\phi(\caL_{V}^\bot)\subset \caL_{\hat{V}}^\bot$, the linear maps,
$$
d\phi_\star
:
\cT_{\star\in \caL_{V}}=V    \lon \cT_{\star\in \caL_{\hat{V}}}=\hat{V}.
$$
and
$$
d\phi_\star^\bot:
\cT_{\star\in \caL_{V}^\bot}=V^*[-1]    \lon \cT_{\star\in \caL_{\hat{V}}\bot}=\hat{V}^*[-1],
$$
of the associated  subspaces.  The differentials $d$  in $V$ and $\hat{d}$ in $\hat{V}$
 induce, respectively, dual differentials
$d^*$ in $V^*[-1]$ and  $\hat{d}^*$ in $\hat{V}^*[-1]$.
\mip

{\bf 3.10.4. Lemma}.
{\em The maps  $d\phi_\star: V\lon \hat{V}$  and
$d\phi_\star^\bot: V^*[-1]\rar \hat{V}^*[-1]$
respect the induced differentials.}

\sip

\Proof Let $(x^a, \psi_a)$ and $(\hat{x}^A, \hat{\psi}_A)$ be arbitrary adopted
Darboux coordinate charts on
odd symplectic manifolds $\cM$ and, respectively, $\hat{\cM}$. The map $\phi$ is given
in these coordinates by
$$
\Ba{l}
\hat{x}^A=  \sum_a
\cA_a^A x^a + \mbox{higher order terms}, \ \ \ \cA_a^A\in \K,\\
\hat{\psi}_A= \sum_a \cB_A^a \psi_a + \mbox{higher order terms},
\ \ \ \ \cB_A^a\in \K.
\Ea
$$
The invertible matrix $\cA_a^B$ (resp.\, $\cB_A^a$) is a coordinate representative of the map
$d\phi_\star$ (resp., $d\phi_\star^\bot$) in the associated
(to a choice of Darboux coordinates) bases of $V$ and $\hat{V}$ (resp., of $V^*[-1]$ and
$\hat{V}^*[-1]$). Equality (\ref{Canonical transform in Darboux}) implies that matrices $\cA$ and $\cB$ are inverse to each other; then equality (\ref{Canonical_transform_of semidensities}) in the limit $\hbar\rar 0$  implies
$$
 \sum_{b}\cA^A_b \Gamma_{(0)\, a}^{\ \ \ b}  =  \sum_{B} \Gamma_{(0)\, B}^{\ \ \ A} \cA^B_a,
  \ \ \
  \sum_{b} \Gamma_{(0)\, b}^{\ \ \ A} \cB^b_A =  \sum_{B} \cB^a_B\Gamma_{(0)\, A}^{\ \ \ B} ,
$$
which in turn implies the required claims.
\hfill $\Box$

\mip

{\bf 3.10.5. Main Theorem.} {\em Every quantum $BV$-manifold $(\cM_V^\hbar, \Theta)$ is symplectomorphic to
the product, $(\cM_{V_1}^\hbar, \Theta_1)\times (\cM_{V_2}^\hbar, \Theta_2)$, of a minimal
quantum BV manifold $(\cM_{V_1}^\hbar, \Theta_1)$ and a contractible one, $(\cM_{V_2}^\hbar, \Theta_2)$}.

\mip

\Proof We shall construct by induction an adopted Darboux coordinate chart $(x^a, \psi_a)_{a\in I}$
on $\cM_V$  in which the quantum density $\Theta$  is represented by
$e^{\frac{\Gamma(x,\psi, \hbar)}{\hbar}}\sqrt{D}_{x,\psi}$
with
\Beq\label{Gamma separated variables}
{\Gamma}({x}, {\psi},\hbar)=
\underbrace{\sum_{A, B\in I'}
\Gamma_A^B x^A \psi_B}_{{\Ga}_0}
\,
+ \,
\underbrace{\sum_{N= 3}^\infty
\sum_{N=p+q+2n\atop p,q\geq 1, n\geq 0}\sum_{\tb_\bullet,\tc_\bullet \in I''} \frac{1}{p!q!}
{\Gamma}^{\ \tc_1\ldots \tc_q}_{(n)\tb_1\ldots \tb_p}x^{\tb_1}\dots x^{\tb_p}\psi_{\tc_1}\dots
\psi_{\tc_q}\hbar^n}_{{\bGa}},
\Eeq
for some partition of the labeling set $I=\{1,2,\ldots, \dim_\K V\}$ into {\em disjoint}\, subsets
$I=I'\coprod I''$. Then the data,
$$
\left(\K[[{x}^A, {\psi}_A, \hbar]], \ \ \Theta_2:=
e^{\frac{{\Ga}_0}{\hbar}}\sqrt{D}_{{x}^A,{\psi}_A}\right)_{A\in I'},
$$
 defines a contractible BV manifold $(\cM_{V_2}^\hbar, \Theta_2)$
while the data,
$$
\left(\K[[{x}^\ta, {\psi}_\ta, \hbar]],\ \ \Theta_1:=
e^{\frac{{{\bGa}}}{\hbar}}\sqrt{D}_{{x}^\ta,{\psi}_\ta}\right)_{\ta\in I''},
$$
defines a minimal BV manifold $(\cM_{V_1}^\hbar, \Theta_1)$ proving thereby the Main Theorem.
The complete separation of variables in (\ref{Gamma separated variables})
assures that $\Theta_1$ and $\Theta_2$ satisfy the corresponding quantum
master equations.

\mip

The required separation of variables
(\ref{Gamma separated variables})
 will be achieved
by induction on an integer valued  parameter $N$ starting with $N=2$. From now on  the {\em order}\, of a monomial
$$
 x^{a_1}\ldots  x^{a_p}\psi_{b_1}\ldots
\psi_{b_q} \hbar^n \in \K[[x,\psi,\hbar]]
$$
is assumed to be $p+q+2n$ and, for an natural number $N$,
 we denote by $\f(N)$ the subset of $\K[[x,\psi,\hbar]] $
consisting of formal power series spanned by monomials of order $\geq N $.
For formal power series $f, g\in \K[[x,\psi,\hbar]]$ the
equality $f=g \mod \f(N)$ means
equality of their polynomial parts of order strictly less than $N$.

\sip

As we already know,
the lowest second order polynomial part,
$$
\Gamma_0:= \sum_{a,b}\Gamma_{(0)\, b}^{\ \ \ a}x^b \psi_a 
$$
of  the master function  $\Gamma(x,\psi, \hbar)$ defines a differential
$d$ in the vector space $V$ (and hence in $V[-1]$).
A choice of an adopted Darboux coordinate system $(x^a, \psi_a)$
determines
the asoociated basis, $\{\psi_a \bmod \cI_\hbar^2\}$, of $V[-1]$ and the (dual)
basis, $\{x^a \bmod \cI_\hbar^2\}$, of $V^*$.
As we are working
over a field of characteristic zero, it is always possible to (non-canonically) represent
the complex $(V[-1],d)$ as a direct sum,
$$
V[-1]= H(V,d)[-1] \oplus B \oplus B[-1]
$$
with the differential $d$ given by $d(a\oplus b \oplus c)= b[-1]$. Let
$\{{\psi}_{\ta}\}_{\ta\in I'}$, be a basis of $H(V,d)[-1]$,
$\{\psi_\al\}_{\al\in J}$ a basis of $B$ and $\{\psi_\zal:=d\psi_\al\}_{\al\in J}$
the associated basis of $B[-1]$. In the basis
$\{{\psi}_\ta, {\psi}_\al, {\psi}_\zal\}$ of
$V[-1]$ the differential
$d$ is given by the block-matrix
\Beq\label{Differential matrix}
d=\left(\Ba{ccc}
0 & 0 & 0\\
0 & 0 & 0\\
0 & \Id & 0
\Ea
\right).
\Eeq
The above splitting induces an associated splitting of $V^*$, and hence an associated
dual base, $\{{x}^\ta, {x}^\al, {x}^\zal\}$ of $V^*$. Thus we can always
find an adopted Darboux coordinate chart
\Beq\label{Daboux chart- Main Theorem}
\left(x^a, \psi_a\right)= \left(
\underbrace{(x^\al, x^\zal, \psi_\al, \psi_\zal)}_{(x^A, \psi_A)},
(x^\ta, \psi_\ta)\right)
\Eeq
in which the master semidensity is given by
\Beqr
\Gamma(x, \psi,\hbar) &=&
\sum_{a,b\in I}\Gamma_{(0)\, b}^{\ \ \ a}x^b \psi_a \mod \f(3)\nonumber \\
&=&
\sum_{\al\in J} x^\al {\psi}_\zal \mod\f(3) \label{Gamma_2_chapter3}.
\Eeqr
Assume now that we have constructed an adopted Darboux chart
 (\ref{Daboux chart- Main Theorem}) in which $\Gamma(x, \psi,\hbar)$ is given by
 (\ref{Gamma separated variables}) modulo terms of order $N+1\geq 3$, i.e.,
$$
\Gamma({x}, {\psi},\hbar)=\sum_{\al\in J} x^\al {\psi}_\zal
\,
+ \, \sum_{k= 3}^N \bGa_k(x^\ta, \psi_\ta, \hbar) \mod \f(N+1)
$$
holds true for some $N\geq 2$. Here $\Ga_k$ stands for a sum of monomials of degree $k$.
It follows from quantum master equation (\ref{Quantum Master equation Darboux})
that the next term, $\Gamma_{N+1}(x^\ta, x^A, \psi_\ta, \psi_A,\hbar)$,  in the Taylor expansion of
$\Gamma$ must satisfy an equation,
\Beq\label{Quantum master eqn N+1}
\{\sum_{\al\in J} x^\zal {\psi}_\al\bullet \Gamma_{N+1}\} + \hbar \Delta_0 \Gamma_{N+1}
+ \frac{1}{2}\sum_{p+q=N+3\atop p,q\geq 3}\left\{\bGa_p\bullet \bGa_q\right\}=0.
\Eeq
The map
$$
\Ba{rccc}
d: & \K[[ x^\ta, x^A, \psi_\ta, \psi_A,\hbar]] & \lon & \K[[x^\ta, x^A, \psi_\ta, \psi_A,\hbar]]\\
\displaystyle
&        f & \lon &
\displaystyle
df:= \{\sum_{\al\in J} x^\al {\psi}_\zal\bullet f\}
\Ea
$$
is a differential which can be equivalently represented as
$$
d=\sum_{\al\in J} \left(\psi_\zal\frac{\p}{\p \psi_\al}+ (-1)^{|x^\al|} x^\al\frac{\p}{\p x^\zal}
   \right), \ \ \ \  |\psi_\zal|= |\psi_\al|+1,\  |x^\al|=|x^\zal|+1.
$$
As $d\Delta_0 + \Delta_0d=0$ and $\Delta_0^2=0$,  the map
 $d+\hbar\Delta_0$ is also a differential in  $\K[[ x^\ta, x^A, \psi_\ta, \psi_A,\hbar]]$.
Thus we can rewrite master equation (\ref{Quantum master eqn N+1}) in the form,
\Beq\label{Induction - main theorem}
(d+\hbar \Delta_0) \Gamma_{N+1}=
{\mathbf F}_{N+1},
\Eeq
where ${\mathbf F}_{N+1}:=  -\frac{1}{2}\sum_{p+q=N+3\atop p,q\geq 3}\{\bGa_p\bullet \bGa_q\}$
does {\em not}\, depend (by the induction assumption) on the variables $(x^A, \psi_A)$.

\mip

{\bf Lemma A}. {\em The vector subspace
$$
i: \K[[x^\ta, \psi_\ta,\hbar]]\subset \K[[x^\ta, x^A, \psi_\ta, \psi_A,\hbar]]
$$
is a subcomplex of the complex
$(\K[[x^\ta, x^A, \psi_\ta, \psi_A,\hbar]], d+ \hbar\Delta_0)$ with the induced differential
$\delta$ being equal to
$\hbar \sum_{\ta\in I''} \frac{\p^2}{\p x^\ta\p\psi_\ta}$.
The inclusion $i$ is a quasi-isomorphism of complexes}.\footnote{It is worth pointing
out that the homology of the complex $(\K[[x^\ta, \psi_\ta]], \frac{\p^2}{\p x^\ta\p\psi_\ta})$
is a one dimensional vector space spanned over $\K$
 by the product, $\eta=x^{\ta'}\cdots \psi_{\ta''}$, of all those
elements of the set, $(x^\ta, \psi_\ta)$, of generators  which have degrees in $2\Z +1$.
Hence the cohomology of the complex $(\K[[x^\ta, \psi_\ta,\hbar]], \delta)$ is equal to
the direct sum $A\oplus \hbar \K[[\hbar]]\ot\eta$, we $A$ is the kernel of
the operator $\frac{\p^2}{\p x^\ta\p\psi_\ta}$ in $\K[[x^\ta, \psi_\ta]]$.}

\mip

{\sc Proof of Lemma A}.
The inclusion $i$ respects the  filtrations,
$$
\K[[x^\ta, \psi_\ta,\hbar]]\ \supset \ \hbar \K[[x^\ta, \psi_\ta,\hbar]]\supset
\hbar^2\K[[x^\ta, \psi_\ta,\hbar]]\ \supset\  \ldots \hspace{31mm}
$$
$$
\K[[x^\ta, x^A, \psi_\ta, \psi_A,\hbar]]\supset \hbar
\K[[x^\ta, x^A, \psi_\ta, \psi_A,\hbar]]\supset \hbar^2
\K[[x^\ta, x^A, \psi_\ta, \psi_A,\hbar]]\supset \ldots
$$
and hence induces maps,
$$
i_r:  (E_r, \delta_r) \lon (\cE_r, D_r)\ \ r\geq 0,
$$
of the associated spectral sequences. The differential $\delta_0$ vanishes while
the differential $D_0$ is equal to $d$.
The Poincar\'{e} Lemma (see, e.g., \S 3.4.5 in \cite{Ma})
says that the cohomology of
the complex $(\cE_0=\K[[x^\ta, x^A, \psi_\ta, \psi_A]], D_0=d)$ is equal
$\K[[x^\ta, \psi_\ta]]=:\cE_1$. Hence the map $i_1: (E_1, \delta_1)\rar (\cE_1, D_1)$
is obviously an isomorphism. Both spectral sequences are regular (terminating at $r=2$),
and the filtrations are
complete and exhaustive. Hence by  classical Complete Convergence Theorem 5.5.10
(see p.139 in \cite{Weib}) they both converge. Then Comparison Theorem 5.2.12
in   \cite{Weib} says that the inclusion $i$ is a quasi-isomorphism completing the proof
of Lemma A. \hfill $\Box$

\mip

{\bf Lemma B}. {\em Every solution, $A=\sum_{k=0}^\infty A_k\hbar^k$, $
A_k \in  \K[[ x^\ta, x^A, \psi_\ta, \psi_A]]$,  of the equation,
$$
 (d+\hbar \Delta_0)A =0, \ \ \ \
$$
can be represented in the form
$$
A= {\mathbf B} + (d+\hbar \Delta_0)C
$$
 for some
${\mathbf B}\in \K[[ x^\ta, \psi_\ta,\hbar]]$
    and $C\in \K[[ x^\ta, x^A, \psi_\ta, \psi_A,\hbar]]$. Moreover, if $A$
    satisfies the boundary conditions (\ref{Boundary conditions for Gamma}), that is,
$$
A_0|_{x=0}=0, \ \     A_0|_{\psi=0}=0, \ \ A_1|_{x=0}\cdot A_1|_{\psi=0}=0,
$$
and has polynomial order $N+1\geq 3$.
then  ${\mathbf B}$ can be chosen to satisfy (\ref{Boundary conditions for Gamma})
and have order $N+1\geq 3$ as well.
    }
\mip

{\sc Proof of Lemma B}. The first part of this Lemma follows, of course, from Lemma A, but
we show another explicit proof which makes the second part of the Lemma immediate.  We have,
$$
dA_0=0, \ \ dA_1=- \Delta_0 A_0,\ \  \ldots\ \ , dA_i=-\Delta_0 A_{i+1},\ \ \ldots
$$
The Poincar\'{e} Lemma
says that the cohomology of
the complex $(\K[[x^\ta, x^A, \psi_\ta, \psi_A]], d)$ is equal $\K[[x^\ta, \psi_\ta]]$.
Hence we get,
$$
A_0= {\mathbf B}_0+ dC_0, \ \ A_1- \Delta_0 C_0={\mathbf B}_1 +  dC_1, \ \
\ldots\ \ , \ A_i - \Delta_0 C_{i-1}={\mathbf B}_i+
dC_i, \ \ \ldots
$$
for some ${\mathbf B}_i \in \K[[x^\ta, \psi_\ta]]$
and $C_i \in \K[[x^\ta, x^A, \psi_\ta, \psi_A]]$, $i=0,1,2,\ldots$.
Thus
$$
A=\underbrace{\sum_{k=0}^\infty {\mathbf B}_k\hbar^k}_{\mathbf B} +
(d+ \hbar \Delta_0)\underbrace{\sum_{k=0}^\infty C_k \hbar^k}_C.
$$
proving the first half of Lemma B. The differentials $d+\hbar\Delta_0$ and $\delta$
preserve the polynomial order, and the splitting homotopy in the proof
of the Poincare Lemma (see p.\ 171 in \cite{Ma}) can also be chosen to be degree preserving.
Thus if $A$ has order $N+1$, then ${\mathbf B}$ and $C$ can also be chosen to have
order $N+1$ (or be zero).
 The boundary conditions for $A$ imply
${\mathbf B}_0|_{x^\ta=0}+ d(C_0|_{x^\ta})=0$ which, by the Poincare Lemma,
in turn implies   ${\mathbf B}_0|_{x^\ta=0}$. 
 Analogously,
 ${\mathbf B}_0|_{\psi_\ta=0}$. 
Note that equation $A_1|_{x=0}\cdot A_1|_{\psi=0}=0$ says that the formal power series
$A_1\in \K[[x, \psi]]$ has no constant (i.e.\ belonging to $\K$) term.
As $C_0$ is of order $N+1\geq 3$ in $x$ and $\psi$,
the power series $\Delta_0 C_0$ has not constant term as well. Then ${\mathbf B}_1$, being the
 part of $ A_1- \Delta_0 C_0$ which is independent of $(x^A, \psi_A)$, has no constant term either
 and hence ${\mathbf B}_1|_{x=0}\cdot {\mathbf B}_1|_{\psi=0}=0$. The proof is completed.
\hfill $\Box$

\mip

{\bf Lemma C}. {\em Every solution of equation (\ref{Induction - main theorem})
can be represented in the form,
$$
\Gamma_{N+1}= \bGa_{N+1} + (d+\hbar \Delta_0)\Psi_{N+1}
$$
 for some
$\bGa_{N+1}\in \K[[ x^\ta, \psi_\ta,\hbar]]$
    and $\Psi_{N+1}\in \K[[ x^\ta, x^A, \psi_\ta, \psi_A,\hbar]]$}.
    Moreover, if $\Gamma_{N+1}$ satisfies the boundary conditions
    (\ref{Boundary conditions for Gamma}) and $N+1\geq 3$, then
  $\bGa_{N+1}$ also satisfies boundary conditions (\ref{Boundary conditions for Gamma}).
\mip

{\sc Proof of Lemma C}. If an element ${\mathbf F}_{N+1}\in \K[[ x^\ta, \psi_\ta,\hbar]]$
is $(d+\hbar \Delta_0)$-exact in $\K[[ x^\ta, x^A, \psi_\ta, \psi_A,\hbar]]$,
then, by Lemma A, it is $\delta$-exact, i.e.\
$$
{\mathbf F}_{N+1}= \delta {\mathbf G}_{N+1}
$$
for some ${\mathbf G}_{N+1}\in \K[[ x^\ta, \psi_\ta,\hbar]]$, and we can rewrite
(\ref{Induction - main theorem}) in the form,
$$
(d+\hbar \Delta_0)(\Gamma_{N+1} - {\mathbf G}_{N+1})=0.
$$
Then the claim follows from Lemma B. \hfill $\Box$

\mip

We continue with an inductive proof of the Main Theorem. Our task now is to show
that one can further adjust a Darboux coordinate chart (\ref{Daboux chart- Main Theorem})
in such a way that decomposition (\ref{Gamma separated variables}) holds true $\bmod \f(N+2)$,
i.e.
\Beq\label{Gamma separated variables-3}
\Gamma({x}, {\psi},\hbar)=\sum_{\al\in J} x^\al {\psi}_\zal
\,
+ \, \sum_{k= 3}^{N+1} \bGa_k(x^\ta, \psi_\ta, \hbar) \mod \f(N+2).
\Eeq

  Let $\Phi_{N+1}\in
\f_{\cM_V}\simeq  \K[[ x^\ta, x^A, \psi_\ta, \psi_A,\hbar]]$
have degree $1$ and order $N+1$. The associated
 degree $0$ Hamiltonian vector field $H_{\Phi_{N+1}}$
on $\cM_V$ generates a one parameter family of canonical transformations
(see \S 3.4.4) which makes sense at $t=1$,
$$
x^a\lon \hat{x}^a=\phi^a(x, \psi, \hbar), \ \ \ \psi_a\lon \hat{\psi}_a=\phi_a(x, \psi, \hbar)
$$
and induces the following change of the coordinate representation of the master function,
\Beq\label{Semidensity_Gamma_trasnformation}
e^{\frac{\hat{\Gamma}(\hat{x}, \hat{\psi}, \hbar)}{\hbar}}
|_{\hat{x}= \phi(x, \psi, \hbar)\atop \hat{\psi}= \phi(x, \psi, \hbar)}
=
\left(Ber\left(\frac{\p(\hat{x},\hat{\psi})}{\p (x,\psi)}\right)\right)^{-1/2}
e^{\frac{\Gamma({x}, {\psi}, \hbar)}{\hbar}},
\Eeq
Equations (\ref{Integral curves}) for the symplectomorphism generated by $H_{\Phi_{N+1}}$
imply,
\Beqrn
\hat{x}^a &=& x^a
- (-1)^{|x^a|}
\frac{\p \Phi_{N+1}}{\p {\psi_a}} \mod \f(N+1),\\
\hat{\psi}_a &=& \psi_a +
 (-1)^{|x^a|}
\frac{\p \Phi_{N+1}}{\p {x^a}} \mod \f(N+1)\\
\Eeqrn
so that
\Beqrn
\hat{\Gamma}(\hat{x}, \hat{\psi}, \hbar)&=&\hat{\Gamma}(x, \psi, \hbar) +
\{\Phi_{N+1}\bullet \hat{\Ga}\}
\mod \f(N+2)\\
&=&\hat{\Gamma}(x, \psi, \hbar) - \{\sum_{\al\in J} x^\al {\psi}_\zal \bullet \Phi_{N+1}\}
\mod \f(N+2)\\
&=&\hat{\Gamma}(x, \psi, \hbar) - d\Phi \mod \f(N+2),
\Eeqrn
and
\Beqrn
Ber\left(\frac{\p(\hat{x},\hat{\psi})}{\p (x,\psi)}\right) &=&
1 + \sum_a\left((-1)^{|x^a|}\frac{\p}{\p x^a} \left( - (-1)^{|x^a|}
\frac{\p \Phi}{\p {\psi_a}}   \right) + (-1)^{|\psi_a|}\frac{\p}{\p \psi_a} \left((-1)^{|x^a|}
\frac{\p \Phi}{\p {x^a}}   \right) \right)\\
&=&1 - 2\Delta_0 \Phi,
\Eeqrn
where we used a well-known fact that $Ber(1+X)=1+ Str(X)$ modulo higher
order polynomials in entries of $X$.
Thus equation (\ref{Semidensity_Gamma_trasnformation}) says
that $\hat{\Gamma}(x,\psi,\hbar)= \Ga(x, \psi, \hbar)\bmod \f(N+1)$ and
$$
\hat{\Gamma}_{N+1} - d \Phi_{N+1}
=\Ga_{N+1} + \hbar \Delta_0\Phi.
$$
Representing $\Gamma_{N+1}$ as in Lemma C, we obtain,
$$
\hat{\Ga}_{N+1}= \bGa_{N+1} + (d+\hbar\Delta_0)(\Phi_{N+1}+ \Psi_{N+1}),
$$
and  conclude that by choosing $\Phi_{N+1}=-\Psi_{N+1}$ we can always adjust
the adopted Darboux coordinate system in such a way that separation of variables
(\ref{Gamma separated variables}) holds true $\bmod \f(N+2)$. The induction completes
proof of the Main Theorem. \hfill $\Box$

\mip

{\bf 3.11. Quantum morphisms of BV manifolds.} A {\em quantum morphism},
$$
\phi_\hbar: \left(\cM_{V}^\hbar, \om, \Theta\right)\lon \left(\cM_{\hat{V}}^\hbar, \hat{\om}, \hat{\Theta}\right)
$$
of quantum BV manifolds is, by definition,
 a morphism of dg $\K[[\hbar]]$-modules (see (\ref{Ch3: odd laplacian on functions hbar})),
$$
\phi^*_\hbar: \left(\f_{\cM_{\hat{V}}^\hbar}, \Delta_{\hat{\om}, \hat{\Theta}}\right)\lon
 \left(\f_{\cM_V^\hbar}, \Delta_{{\om},{\Theta}}\right)
$$
inducing in the classical limit $\hbar\rar 0$ a morphism of algebras, $\phi^*_{0}:
\f_{\cM_{\hat{V}}^0} \rar \f_{\cM_{{V}}^0}$ which preserves the ideals
of the distinguished Lagrangian submanifolds in $\cM_V^0$ and $\cM_{\hat{V}}^0$.

\sip

It is easy to see that any quantum morphism $\phi_\hbar: \left(\cM_{V}^\hbar, \om, \Theta\right)\lon \left(\cM_{\hat{V}}^\hbar, \hat{\om}, \hat{\Theta}\right)$
induces a morphism, $d\phi_0: (V,d)\rar (\hat{V},\hat{d})$, of the associated tangent complexes;
such a morphism is called  a {\em quasi-isomorphism}\, if the  map $d\phi_0$ induces an isomorphism of the associated cohomology groups.

\sip

Note that a quantum morphism is a morphism of algebras only in the classical limit; therefore, in general, it is
not a morphism of smooth manifolds and can not be characterized in local coordinates (i.e.\ in terms
of generators of the structure
sheaves). Let us denote by $\widehat{Cat}(BV)$ the associated
to the above definition of quantum morphisms the category of quantum BV manifolds.

\mip

{\bf 3.11.1. Examples.} (i) {\em Symplectomorphisms}\, of quantum BV manifolds are obviously quantum morphisms.

\sip

(ii) Natural {\em projections},
$$
\phi_\hbar: \cM_{V}^\hbar \times \cM_{\hat{V}}^\hbar \lon \cM_{V}^\hbar
$$
are obviously quantum morphisms.

\mip

The above two examples are  special in the sense that the associated maps of dg $\K[[\hbar]]$-modules,
$\phi^*_\hbar: (\f_{\cM_{\hat{V}}^\hbar}, \Delta_{\hat{\om}, \hat{\Theta}})\rar
 (\f_{\cM_V^\hbar}, \Delta_{{\om},{\Theta}})$, are maps of $\K[[\hbar]]$-algebras. The next example
 does {\em not}\,  have this property in general.

 \mip

(iii) Let
$\left(\cM_{V}^\hbar, \om, \Theta\right)$ and $\left(\cM_{\hat{V}}^\hbar, \hat{\om}, \hat{\Theta}\right)$
be quantum BV manifolds. It is a well-known and very useful fact  \cite{Sc} that, for a Lagrangian submanifold $\caL^\hbar\subset \cM_{\hat{V}}^\hbar$,
the associated integration map
$$
\Ba{rccc}
\phi_\hbar^*: & \f_{\cM_{V}^\hbar\times \cM_{\hat{V}}^\hbar} & \lon &  \f_{\cM_{V}^\hbar}\\
              &      f   & \lon &        {\int_{\caL^\hbar} f \hat{\Theta}}.
\Ea
$$
satisfies,
\Beqrn
\phi_\hbar^*\left((\Delta_{{\om},{\Theta}}+ \Delta_{\hat{\om}, \hat{\Theta}})f\right) &=&
 {\int_{\caL^\hbar} \left((\Delta_{{\om},{\Theta}}+ \Delta_{\hat{\om}, \hat{\Theta}})f\right) \hat{\Theta}}\\
 &=& \hbar \int_{\caL^\hbar} \frac{\Delta_{{\om}}\left(f {\Theta}\right)}{\Theta}\hat{\Theta} + \hbar
 \int_{\caL^\hbar}\Delta_{{\hat{\om}}}\left(f \hat{\Theta}\right)\\
 &=& \hbar \frac{\Delta_{{\om}}\left(  \int_{\caL^\hbar} f \hat{\Theta}\right)}{\Theta}\\
 &=& \Delta_{{\om},{\Theta}} \phi_h^*(f),
\Eeqrn
and is, therefore,
a quantum morphism provided the integral exists as a perturbative series in $\hbar$.
We shall see \S 6 that
such a {\em quantum embedding}\,
$$
\phi_h: \cM_{V}^\hbar \lon \cM_{V}^\hbar\times \cM_{\hat{V}}^\hbar
$$
can always be constructed (as a well-defined formal power series in $\hbar$  satisfying the algebra morphism
condition in the limit $\hbar\rar 0$) in the case when
 the quantum manifold $\cM_{\hat{V}}^\hbar$ is contractible; in the latter case the  quantum embedding is also called {\em contractible}; such quantum embeddings are often given  by Feynman type sums over decorated graphs.

\mip
{\bf 3.11.2. Proposition}. {\em For any quantum BV manifold $\cM_V^\hbar$ and its any decomposition,
$\cM_V^\hbar\simeq \cM_{min}^\hbar \times \cM_{cntr}^\hbar$, into a product of a minimal quantum BV manifold
and a contractible one, there exists a contractible quantum embedding,
$$
\phi_h:  \cM_{min}^\hbar \lon \cM_V^\hbar,
$$
such that $\pi_\hbar \circ \phi_\hbar=\Id$, where $\pi_\hbar$ is the composition
$\cM_V^\hbar\stackrel{\simeq}{\rar} \cM_{min}^\hbar \times \cM_{cntr}^\hbar\rar  \cM_{min}^\hbar$.}

\mip

We shall prove this statement in \S 6 below by giving an explicit formula for $\phi_\hbar$. This fact has an important corollary which we discuss next.

\sip

Let $Cat(BV)$ be the full subcategory of $\widehat{Cat}(BV)$ whose class of morphisms consists, by definition,
of all possible compositions of symplectomorphisms, projections and contractible quantum embeddings. Then Theorem 3.10.5 and Proposition 3.11.2
imply that quasi-isomorphisms in this category are equivalence relations. Therefore, in the homotopy theory sense, the category
 $Cat(BV)$ is as good as, for example, the famous category of strong homotopy Lie algebras \cite{Ko,St}.

\bip

\mip

\section{From unimodular Lie 1-bialgebras to  quantum BV manifolds}
\bip

{\bf 4.1. Lie $n$-bialgebras} \cite{Me1,Me-lec}. A  {\em Lie n-bialgebra}\, is a graded vector space $V$,
equipped with linear maps,
$$
\Delta: V\rightarrow V\wedge V \ \ \ \mbox{and}\ \ \  [\ \bullet\ ]: \wedge^2 (V[-n]) \rightarrow V[-n],
$$
such that
\Bi
\item the data $(V,\delta)$ is a Lie coalgebra;
\item the data $(V[-n], [\ \bullet\ ])$ is a Lie algebra;
\item the compatibility condition,
$$
\Delta[a\bullet b] = \sum a_1\otimes [a_2\bullet b] +  [a\bullet
b_1]\otimes b_2 + (-1)^{|a||b|+n|a|+n|b|}( [b\bullet a_1]\otimes a_2
+ b_1\otimes [b_2\bullet a]),
$$
holds for any $a,b\in V$. Here $\Delta a=:\sum a_1\otimes a_2$, $\Delta b=:\sum
b_1\otimes b_2$.
\Ei
The case  $n=0$  gives us the ordinary definition of Lie bialgebra \cite{D1}.
The case $n=1$  is of most interest to us in this paper as it controls Poisson geometry \cite{Me1}
and, with unimodularity conditions added, controls the category of quantum BV manifolds (see \S 4.3 below).
Note that in this case one has
$\wedge^2 (V[-1])= (\odot^2V)[-2]$ so that the brackets $[\ \bullet\ ]$  describe a degree 1 linear map $\odot^2V
\rar V$.

\mip

\noindent{\bf 4.1.1.  Wheeled prop(erad) of Lie 1-bialgebras.} This is a wheeled prop(erad),
$\LB^\circlearrowright:= \cF^\circlearrowright\langle E\rangle/ <\cR>$, defined as the quotient
of the free wheeled prop(erad)
 generated by an  $\bS$-bimodule
\Beq\label{Generators of LieB}
E(m,n):=\left\{
\Ba{rr}
sgn_2\ot \id_1\equiv\mbox{span}\left\langle
\begin{xy}
 <0mm,-0.55mm>*{};<0mm,-2.5mm>*{}**@{-},
 <0.5mm,0.5mm>*{};<2.2mm,2.2mm>*{}**@{-},
 <-0.48mm,0.48mm>*{};<-2.2mm,2.2mm>*{}**@{-},
 <0mm,0mm>*{\circ};<0mm,0mm>*{}**@{},
 <0mm,-0.55mm>*{};<0mm,-3.8mm>*{_1}**@{},
 <0.5mm,0.5mm>*{};<2.7mm,2.8mm>*{^2}**@{},
 <-0.48mm,0.48mm>*{};<-2.7mm,2.8mm>*{^1}**@{},
 \end{xy}
=-
\begin{xy}
 <0mm,-0.55mm>*{};<0mm,-2.5mm>*{}**@{-},
 <0.5mm,0.5mm>*{};<2.2mm,2.2mm>*{}**@{-},
 <-0.48mm,0.48mm>*{};<-2.2mm,2.2mm>*{}**@{-},
 <0mm,0mm>*{\circ};<0mm,0mm>*{}**@{},
 <0mm,-0.55mm>*{};<0mm,-3.8mm>*{_1}**@{},
 <0.5mm,0.5mm>*{};<2.7mm,2.8mm>*{^1}**@{},
 <-0.48mm,0.48mm>*{};<-2.7mm,2.8mm>*{^2}**@{},
 \end{xy}
   \right\rangle  & \mbox{if}\ m=2, n=1,\vspace{3mm}\\
\id_1\ot \id_2[-1]\equiv
\mbox{span}\left\langle
\begin{xy}
 <0mm,0.66mm>*{};<0mm,3mm>*{}**@{-},
 <0.39mm,-0.39mm>*{};<2.2mm,-2.2mm>*{}**@{-},
 <-0.35mm,-0.35mm>*{};<-2.2mm,-2.2mm>*{}**@{-},
 <0mm,0mm>*{\circ};<0mm,0mm>*{}**@{},
   <0mm,0.66mm>*{};<0mm,3.4mm>*{^1}**@{},
   <0.39mm,-0.39mm>*{};<2.9mm,-4mm>*{^2}**@{},
   <-0.35mm,-0.35mm>*{};<-2.8mm,-4mm>*{^1}**@{},
\end{xy}=
\begin{xy}
 <0mm,0.66mm>*{};<0mm,3mm>*{}**@{-},
 <0.39mm,-0.39mm>*{};<2.2mm,-2.2mm>*{}**@{-},
 <-0.35mm,-0.35mm>*{};<-2.2mm,-2.2mm>*{}**@{-},
 <0mm,0mm>*{\circ};<0mm,0mm>*{}**@{},
   <0mm,0.66mm>*{};<0mm,3.4mm>*{^1}**@{},
   <0.39mm,-0.39mm>*{};<2.9mm,-4mm>*{^1}**@{},
   <-0.35mm,-0.35mm>*{};<-2.8mm,-4mm>*{^2}**@{},
\end{xy}
\right\rangle
\ & \mbox{if}\ m=1, n=2, \vspace{3mm}\\
0 & \mbox{otherwise}
\Ea
\right.
\Eeq
by the ideal generated by the relations
\Beq\label{R for LieB}
\cR:\left\{
\Ba{r}
\begin{xy}
 <0mm,0mm>*{\circ};<0mm,0mm>*{}**@{},
 <0mm,-0.49mm>*{};<0mm,-3.0mm>*{}**@{-},
 <0.49mm,0.49mm>*{};<1.9mm,1.9mm>*{}**@{-},
 <-0.5mm,0.5mm>*{};<-1.9mm,1.9mm>*{}**@{-},
 <-2.3mm,2.3mm>*{\circ};<-2.3mm,2.3mm>*{}**@{},
 <-1.8mm,2.8mm>*{};<0mm,4.9mm>*{}**@{-},
 <-2.8mm,2.9mm>*{};<-4.6mm,4.9mm>*{}**@{-},
   <0.49mm,0.49mm>*{};<2.7mm,2.3mm>*{^3}**@{},
   <-1.8mm,2.8mm>*{};<0.4mm,5.3mm>*{^2}**@{},
   <-2.8mm,2.9mm>*{};<-5.1mm,5.3mm>*{^1}**@{},
 \end{xy}
\ + \
\begin{xy}
 <0mm,0mm>*{\circ};<0mm,0mm>*{}**@{},
 <0mm,-0.49mm>*{};<0mm,-3.0mm>*{}**@{-},
 <0.49mm,0.49mm>*{};<1.9mm,1.9mm>*{}**@{-},
 <-0.5mm,0.5mm>*{};<-1.9mm,1.9mm>*{}**@{-},
 <-2.3mm,2.3mm>*{\circ};<-2.3mm,2.3mm>*{}**@{},
 <-1.8mm,2.8mm>*{};<0mm,4.9mm>*{}**@{-},
 <-2.8mm,2.9mm>*{};<-4.6mm,4.9mm>*{}**@{-},
   <0.49mm,0.49mm>*{};<2.7mm,2.3mm>*{^2}**@{},
   <-1.8mm,2.8mm>*{};<0.4mm,5.3mm>*{^1}**@{},
   <-2.8mm,2.9mm>*{};<-5.1mm,5.3mm>*{^3}**@{},
 \end{xy}
\ + \
\begin{xy}
 <0mm,0mm>*{\circ};<0mm,0mm>*{}**@{},
 <0mm,-0.49mm>*{};<0mm,-3.0mm>*{}**@{-},
 <0.49mm,0.49mm>*{};<1.9mm,1.9mm>*{}**@{-},
 <-0.5mm,0.5mm>*{};<-1.9mm,1.9mm>*{}**@{-},
 <-2.3mm,2.3mm>*{\circ};<-2.3mm,2.3mm>*{}**@{},
 <-1.8mm,2.8mm>*{};<0mm,4.9mm>*{}**@{-},
 <-2.8mm,2.9mm>*{};<-4.6mm,4.9mm>*{}**@{-},
   <0.49mm,0.49mm>*{};<2.7mm,2.3mm>*{^1}**@{},
   <-1.8mm,2.8mm>*{};<0.4mm,5.3mm>*{^3}**@{},
   <-2.8mm,2.9mm>*{};<-5.1mm,5.3mm>*{^2}**@{},
 \end{xy}
\ \ \ \in \cF_{(2)}^\circlearrowright \langle E\rangle(3,1) \vspace{3mm}\\
 \begin{xy}
 <0mm,0mm>*{\circ};<0mm,0mm>*{}**@{},
 <0mm,0.69mm>*{};<0mm,3.0mm>*{}**@{-},
 <0.39mm,-0.39mm>*{};<2.4mm,-2.4mm>*{}**@{-},
 <-0.35mm,-0.35mm>*{};<-1.9mm,-1.9mm>*{}**@{-},
 <-2.4mm,-2.4mm>*{\circ};<-2.4mm,-2.4mm>*{}**@{},
 <-2.0mm,-2.8mm>*{};<0mm,-4.9mm>*{}**@{-},
 <-2.8mm,-2.9mm>*{};<-4.7mm,-4.9mm>*{}**@{-},
    <0.39mm,-0.39mm>*{};<3.3mm,-4.0mm>*{^3}**@{},
    <-2.0mm,-2.8mm>*{};<0.5mm,-6.7mm>*{^2}**@{},
    <-2.8mm,-2.9mm>*{};<-5.2mm,-6.7mm>*{^1}**@{},
 \end{xy}
\ + \
 \begin{xy}
 <0mm,0mm>*{\circ};<0mm,0mm>*{}**@{},
 <0mm,0.69mm>*{};<0mm,3.0mm>*{}**@{-},
 <0.39mm,-0.39mm>*{};<2.4mm,-2.4mm>*{}**@{-},
 <-0.35mm,-0.35mm>*{};<-1.9mm,-1.9mm>*{}**@{-},
 <-2.4mm,-2.4mm>*{\circ};<-2.4mm,-2.4mm>*{}**@{},
 <-2.0mm,-2.8mm>*{};<0mm,-4.9mm>*{}**@{-},
 <-2.8mm,-2.9mm>*{};<-4.7mm,-4.9mm>*{}**@{-},
    <0.39mm,-0.39mm>*{};<3.3mm,-4.0mm>*{^2}**@{},
    <-2.0mm,-2.8mm>*{};<0.5mm,-6.7mm>*{^1}**@{},
    <-2.8mm,-2.9mm>*{};<-5.2mm,-6.7mm>*{^3}**@{},
 \end{xy}
\ + \
 \begin{xy}
 <0mm,0mm>*{\circ};<0mm,0mm>*{}**@{},
 <0mm,0.69mm>*{};<0mm,3.0mm>*{}**@{-},
 <0.39mm,-0.39mm>*{};<2.4mm,-2.4mm>*{}**@{-},
 <-0.35mm,-0.35mm>*{};<-1.9mm,-1.9mm>*{}**@{-},
 <-2.4mm,-2.4mm>*{\circ};<-2.4mm,-2.4mm>*{}**@{},
 <-2.0mm,-2.8mm>*{};<0mm,-4.9mm>*{}**@{-},
 <-2.8mm,-2.9mm>*{};<-4.7mm,-4.9mm>*{}**@{-},
    <0.39mm,-0.39mm>*{};<3.3mm,-4.0mm>*{^1}**@{},
    <-2.0mm,-2.8mm>*{};<0.5mm,-6.7mm>*{^3}**@{},
    <-2.8mm,-2.9mm>*{};<-5.2mm,-6.7mm>*{^2}**@{},
 \end{xy}
\  \ \ \in   \cF_{(2)}^\circlearrowright \langle E\rangle(1,3) \vspace{2mm} \\
 \begin{xy}
 <0mm,2.47mm>*{};<0mm,0.12mm>*{}**@{-},
 <0.5mm,3.5mm>*{};<2.2mm,5.2mm>*{}**@{-},
 <-0.48mm,3.48mm>*{};<-2.2mm,5.2mm>*{}**@{-},
 <0mm,3mm>*{\circ};<0mm,3mm>*{}**@{},
  <0mm,-0.8mm>*{\circ};<0mm,-0.8mm>*{}**@{},
<-0.39mm,-1.2mm>*{};<-2.2mm,-3.5mm>*{}**@{-},
 <0.39mm,-1.2mm>*{};<2.2mm,-3.5mm>*{}**@{-},
     <0.5mm,3.5mm>*{};<2.8mm,5.7mm>*{^2}**@{},
     <-0.48mm,3.48mm>*{};<-2.8mm,5.7mm>*{^1}**@{},
   <0mm,-0.8mm>*{};<-2.7mm,-5.2mm>*{^1}**@{},
   <0mm,-0.8mm>*{};<2.7mm,-5.2mm>*{^2}**@{},
\end{xy}
\  - \
\begin{xy}
 <0mm,-1.3mm>*{};<0mm,-3.5mm>*{}**@{-},
 <0.38mm,-0.2mm>*{};<2.0mm,2.0mm>*{}**@{-},
 <-0.38mm,-0.2mm>*{};<-2.2mm,2.2mm>*{}**@{-},
<0mm,-0.8mm>*{\circ};<0mm,0.8mm>*{}**@{},
 <2.4mm,2.4mm>*{\circ};<2.4mm,2.4mm>*{}**@{},
 <2.77mm,2.0mm>*{};<4.4mm,-0.8mm>*{}**@{-},
 <2.4mm,3mm>*{};<2.4mm,5.2mm>*{}**@{-},
     <0mm,-1.3mm>*{};<0mm,-5.3mm>*{^1}**@{},
     <2.5mm,2.3mm>*{};<5.1mm,-2.6mm>*{^2}**@{},
    <2.4mm,2.5mm>*{};<2.4mm,5.7mm>*{^2}**@{},
    <-0.38mm,-0.2mm>*{};<-2.8mm,2.5mm>*{^1}**@{},
    \end{xy}
\  - \
\begin{xy}
 <0mm,-1.3mm>*{};<0mm,-3.5mm>*{}**@{-},
 <0.38mm,-0.2mm>*{};<2.0mm,2.0mm>*{}**@{-},
 <-0.38mm,-0.2mm>*{};<-2.2mm,2.2mm>*{}**@{-},
<0mm,-0.8mm>*{\circ};<0mm,0.8mm>*{}**@{},
 <2.4mm,2.4mm>*{\circ};<2.4mm,2.4mm>*{}**@{},
 <2.77mm,2.0mm>*{};<4.4mm,-0.8mm>*{}**@{-},
 <2.4mm,3mm>*{};<2.4mm,5.2mm>*{}**@{-},
     <0mm,-1.3mm>*{};<0mm,-5.3mm>*{^2}**@{},
     <2.5mm,2.3mm>*{};<5.1mm,-2.6mm>*{^1}**@{},
    <2.4mm,2.5mm>*{};<2.4mm,5.7mm>*{^2}**@{},
    <-0.38mm,-0.2mm>*{};<-2.8mm,2.5mm>*{^1}**@{},
    \end{xy}
\  + \
\begin{xy}
 <0mm,-1.3mm>*{};<0mm,-3.5mm>*{}**@{-},
 <0.38mm,-0.2mm>*{};<2.0mm,2.0mm>*{}**@{-},
 <-0.38mm,-0.2mm>*{};<-2.2mm,2.2mm>*{}**@{-},
<0mm,-0.8mm>*{\circ};<0mm,0.8mm>*{}**@{},
 <2.4mm,2.4mm>*{\circ};<2.4mm,2.4mm>*{}**@{},
 <2.77mm,2.0mm>*{};<4.4mm,-0.8mm>*{}**@{-},
 <2.4mm,3mm>*{};<2.4mm,5.2mm>*{}**@{-},
     <0mm,-1.3mm>*{};<0mm,-5.3mm>*{^2}**@{},
     <2.5mm,2.3mm>*{};<5.1mm,-2.6mm>*{^1}**@{},
    <2.4mm,2.5mm>*{};<2.4mm,5.7mm>*{^1}**@{},
    <-0.38mm,-0.2mm>*{};<-2.8mm,2.5mm>*{^2}**@{},
    \end{xy}
\ + \
\begin{xy}
 <0mm,-1.3mm>*{};<0mm,-3.5mm>*{}**@{-},
 <0.38mm,-0.2mm>*{};<2.0mm,2.0mm>*{}**@{-},
 <-0.38mm,-0.2mm>*{};<-2.2mm,2.2mm>*{}**@{-},
<0mm,-0.8mm>*{\circ};<0mm,0.8mm>*{}**@{},
 <2.4mm,2.4mm>*{\circ};<2.4mm,2.4mm>*{}**@{},
 <2.77mm,2.0mm>*{};<4.4mm,-0.8mm>*{}**@{-},
 <2.4mm,3mm>*{};<2.4mm,5.2mm>*{}**@{-},
     <0mm,-1.3mm>*{};<0mm,-5.3mm>*{^1}**@{},
     <2.5mm,2.3mm>*{};<5.1mm,-2.6mm>*{^2}**@{},
    <2.4mm,2.5mm>*{};<2.4mm,5.7mm>*{^1}**@{},
    <-0.38mm,-0.2mm>*{};<-2.8mm,2.5mm>*{^2}**@{},
    \end{xy}
\ \ \ \in   \cF_{(2)}^\circlearrowright \langle E\rangle(2,2)
\Ea
\right.
\Eeq
It is clear from the association
$$
\vartriangle \leftrightarrow
 \begin{xy}
 <0mm,-0.55mm>*{};<0mm,-2.5mm>*{}**@{-},
 <0.5mm,0.5mm>*{};<2.2mm,2.2mm>*{}**@{-},
 <-0.48mm,0.48mm>*{};<-2.2mm,2.2mm>*{}**@{-},
 <0mm,0mm>*{\circ};<0mm,0mm>*{}**@{},
 \end{xy}
 \ \ \  , \ \ \
[\, \bullet \, ] \leftrightarrow
 \begin{xy}
 <0mm,0.66mm>*{};<0mm,3mm>*{}**@{-},
 <0.39mm,-0.39mm>*{};<2.2mm,-2.2mm>*{}**@{-},
 <-0.35mm,-0.35mm>*{};<-2.2mm,-2.2mm>*{}**@{-},
 <0mm,0mm>*{\circ};<0mm,0mm>*{}**@{},
 \end{xy}
$$
that there is a one-to-one correspondence between representations of $\LB^\circlearrowright$ in a finite
dimensional space $V$ and Lie 1-bialgebra structures in $V$.

\mip

{\bf 4.1.2. Cobar construction on the Koszul dual coproperad $(\LB^\circlearrowright)^\Koz$}. It follows from the exact sequence
(\ref{Relations for P!})  that the Koszul dual wheeled properad $(\LB^\circlearrowright)^!$ is the quotient,
$\cF^\circlearrowright\langle E^\vee\rangle/ <\cR^\bot>$, of the free wheeled prop(erad) generated
by the $\bS$-bimodule,
\Beq\label{Generators ULieB-vee}
E^\vee(m,n):=\left\{
\Ba{rr}
sgn_2\ot \id_1[1]\equiv\mbox{span}\left\langle
\begin{xy}
 <0mm,-0.55mm>*{};<0mm,-2.5mm>*{}**@{-},
 <0.5mm,0.5mm>*{};<2.2mm,2.2mm>*{}**@{-},
 <-0.48mm,0.48mm>*{};<-2.2mm,2.2mm>*{}**@{-},
 <0mm,0mm>*{\circ};<0mm,0mm>*{}**@{},
 <0mm,-0.55mm>*{};<0mm,-3.8mm>*{_1}**@{},
 <0.5mm,0.5mm>*{};<2.7mm,2.8mm>*{^2}**@{},
 <-0.48mm,0.48mm>*{};<-2.7mm,2.8mm>*{^1}**@{},
 \end{xy}
=-
\begin{xy}
 <0mm,-0.55mm>*{};<0mm,-2.5mm>*{}**@{-},
 <0.5mm,0.5mm>*{};<2.2mm,2.2mm>*{}**@{-},
 <-0.48mm,0.48mm>*{};<-2.2mm,2.2mm>*{}**@{-},
 <0mm,0mm>*{\circ};<0mm,0mm>*{}**@{},
 <0mm,-0.55mm>*{};<0mm,-3.8mm>*{_1}**@{},
 <0.5mm,0.5mm>*{};<2.7mm,2.8mm>*{^1}**@{},
 <-0.48mm,0.48mm>*{};<-2.7mm,2.8mm>*{^2}**@{},
 \end{xy}
   \right\rangle  & \mbox{if}\ m=2, n=1,\vspace{3mm}\\
\id_1\ot sgn_2[1]\equiv
\mbox{span}\left\langle
\begin{xy}
 <0mm,0.66mm>*{};<0mm,3mm>*{}**@{-},
 <0.39mm,-0.39mm>*{};<2.2mm,-2.2mm>*{}**@{-},
 <-0.35mm,-0.35mm>*{};<-2.2mm,-2.2mm>*{}**@{-},
 <0mm,0mm>*{\circ};<0mm,0mm>*{}**@{},
   <0mm,0.66mm>*{};<0mm,3.4mm>*{^1}**@{},
   <0.39mm,-0.39mm>*{};<2.9mm,-4mm>*{^2}**@{},
   <-0.35mm,-0.35mm>*{};<-2.8mm,-4mm>*{^1}**@{},
\end{xy}=-
\begin{xy}
 <0mm,0.66mm>*{};<0mm,3mm>*{}**@{-},
 <0.39mm,-0.39mm>*{};<2.2mm,-2.2mm>*{}**@{-},
 <-0.35mm,-0.35mm>*{};<-2.2mm,-2.2mm>*{}**@{-},
 <0mm,0mm>*{\circ};<0mm,0mm>*{}**@{},
   <0mm,0.66mm>*{};<0mm,3.4mm>*{^1}**@{},
   <0.39mm,-0.39mm>*{};<2.9mm,-4mm>*{^1}**@{},
   <-0.35mm,-0.35mm>*{};<-2.8mm,-4mm>*{^2}**@{},
\end{xy}
\right\rangle
\ & \mbox{if}\ m=1, n=2, \vspace{3mm}\\
0 & \mbox{otherwise},
\Ea
\right.
\Eeq
by the ideal generated by relations
$$
\cR^\bot:\left\{
\Ba{c}
\begin{xy}
 <0mm,-0.55mm>*{};<0mm,-2.5mm>*{}**@{-},
 <0.5mm,0.5mm>*{};<2.2mm,2.2mm>*{}**@{-},
 <-0.48mm,0.48mm>*{};<-2.2mm,2.2mm>*{}**@{-},
 <0mm,0mm>*{\circ};<0mm,0mm>*{}**@{},
 <-0.48mm,0.48mm>*{};<-2.7mm,2.8mm>*{^1}**@{},
(2.2,2.2)*{}
   \ar@{->}@(ur,dr) (0,-2.5)*{}
 \end{xy}=0\ , \ \ \
\begin{xy}
 <0mm,0.66mm>*{};<0mm,3mm>*{}**@{-},
 <0.39mm,-0.39mm>*{};<2.2mm,-2.2mm>*{}**@{-},
 <-0.35mm,-0.35mm>*{};<-2.2mm,-2.2mm>*{}**@{-},
 <0mm,0mm>*{\circ};<0mm,0mm>*{}**@{},
   <-0.35mm,-0.35mm>*{};<-2.8mm,-4mm>*{^1}**@{},
(0.0,3.0)*{}
   \ar@{->}@(ur,dr) (2.2,-2.2)*{}
\end{xy}=0
\\
\Ba{c}
\begin{xy}
 <0mm,0mm>*{\circ};<0mm,0mm>*{}**@{},
 <0mm,-0.49mm>*{};<0mm,-3.0mm>*{}**@{-},
 <0.49mm,0.49mm>*{};<1.9mm,1.9mm>*{}**@{-},
 <-0.5mm,0.5mm>*{};<-1.9mm,1.9mm>*{}**@{-},
 <-2.3mm,2.3mm>*{\circ};<-2.3mm,2.3mm>*{}**@{},
 <-1.8mm,2.8mm>*{};<0mm,4.9mm>*{}**@{-},
 <-2.8mm,2.9mm>*{};<-4.6mm,4.9mm>*{}**@{-},
   <0.49mm,0.49mm>*{};<2.7mm,2.3mm>*{^3}**@{},
   <-1.8mm,2.8mm>*{};<0.4mm,5.3mm>*{^2}**@{},
   <-2.8mm,2.9mm>*{};<-5.1mm,5.3mm>*{^1}**@{},
 \end{xy}
\ - \
\begin{xy}
 <0mm,0mm>*{\circ};<0mm,0mm>*{}**@{},
 <0mm,-0.49mm>*{};<0mm,-3.0mm>*{}**@{-},
 <0.49mm,0.49mm>*{};<1.9mm,1.9mm>*{}**@{-},
 <-0.5mm,0.5mm>*{};<-1.9mm,1.9mm>*{}**@{-},
 <-2.3mm,2.3mm>*{\circ};<-2.3mm,2.3mm>*{}**@{},
 <-1.8mm,2.8mm>*{};<0mm,4.9mm>*{}**@{-},
 <-2.8mm,2.9mm>*{};<-4.6mm,4.9mm>*{}**@{-},
   <0.49mm,0.49mm>*{};<2.7mm,2.3mm>*{^1}**@{},
   <-1.8mm,2.8mm>*{};<0.4mm,5.3mm>*{^3}**@{},
   <-2.8mm,2.9mm>*{};<-5.1mm,5.3mm>*{^2}**@{},
 \end{xy}
\Ea
=0, \ \
\Ba{c}
 \begin{xy}
 <0mm,0mm>*{\circ};<0mm,0mm>*{}**@{},
 <0mm,0.69mm>*{};<0mm,3.0mm>*{}**@{-},
 <0.39mm,-0.39mm>*{};<2.4mm,-2.4mm>*{}**@{-},
 <-0.35mm,-0.35mm>*{};<-1.9mm,-1.9mm>*{}**@{-},
 <-2.4mm,-2.4mm>*{\circ};<-2.4mm,-2.4mm>*{}**@{},
 <-2.0mm,-2.8mm>*{};<0mm,-4.9mm>*{}**@{-},
 <-2.8mm,-2.9mm>*{};<-4.7mm,-4.9mm>*{}**@{-},
    <0.39mm,-0.39mm>*{};<3.3mm,-4.0mm>*{^3}**@{},
    <-2.0mm,-2.8mm>*{};<0.5mm,-6.7mm>*{^2}**@{},
    <-2.8mm,-2.9mm>*{};<-5.2mm,-6.7mm>*{^1}**@{},
 \end{xy}
\ - \
 \begin{xy}
 <0mm,0mm>*{\circ};<0mm,0mm>*{}**@{},
 <0mm,0.69mm>*{};<0mm,3.0mm>*{}**@{-},
 <0.39mm,-0.39mm>*{};<2.4mm,-2.4mm>*{}**@{-},
 <-0.35mm,-0.35mm>*{};<-1.9mm,-1.9mm>*{}**@{-},
 <-2.4mm,-2.4mm>*{\circ};<-2.4mm,-2.4mm>*{}**@{},
 <-2.0mm,-2.8mm>*{};<0mm,-4.9mm>*{}**@{-},
 <-2.8mm,-2.9mm>*{};<-4.7mm,-4.9mm>*{}**@{-},
    <0.39mm,-0.39mm>*{};<3.3mm,-4.0mm>*{^1}**@{},
    <-2.0mm,-2.8mm>*{};<0.5mm,-6.7mm>*{^3}**@{},
    <-2.8mm,-2.9mm>*{};<-5.2mm,-6.7mm>*{^2}**@{},
 \end{xy}
\Ea=0, \ \
\Ba{c}
 \begin{xy}
 <0mm,2.47mm>*{};<0mm,0.12mm>*{}**@{-},
 <0.5mm,3.5mm>*{};<2.2mm,5.2mm>*{}**@{-},
 <-0.48mm,3.48mm>*{};<-2.2mm,5.2mm>*{}**@{-},
 <0mm,3mm>*{\circ};<0mm,3mm>*{}**@{},
  <0mm,-0.8mm>*{\circ};<0mm,-0.8mm>*{}**@{},
<-0.39mm,-1.2mm>*{};<-2.2mm,-3.5mm>*{}**@{-},
 <0.39mm,-1.2mm>*{};<2.2mm,-3.5mm>*{}**@{-},
     <0.5mm,3.5mm>*{};<2.8mm,5.7mm>*{^2}**@{},
     <-0.48mm,3.48mm>*{};<-2.8mm,5.7mm>*{^1}**@{},
   <0mm,-0.8mm>*{};<-2.7mm,-5.2mm>*{^1}**@{},
   <0mm,-0.8mm>*{};<2.7mm,-5.2mm>*{^2}**@{},
\end{xy}
\  + \
\begin{xy}
 <0mm,-1.3mm>*{};<0mm,-3.5mm>*{}**@{-},
 <0.38mm,-0.2mm>*{};<2.0mm,2.0mm>*{}**@{-},
 <-0.38mm,-0.2mm>*{};<-2.2mm,2.2mm>*{}**@{-},
<0mm,-0.8mm>*{\circ};<0mm,0.8mm>*{}**@{},
 <2.4mm,2.4mm>*{\circ};<2.4mm,2.4mm>*{}**@{},
 <2.77mm,2.0mm>*{};<4.4mm,-0.8mm>*{}**@{-},
 <2.4mm,3mm>*{};<2.4mm,5.2mm>*{}**@{-},
     <0mm,-1.3mm>*{};<0mm,-5.3mm>*{^1}**@{},
     <2.5mm,2.3mm>*{};<5.1mm,-2.6mm>*{^2}**@{},
    <2.4mm,2.5mm>*{};<2.4mm,5.7mm>*{^2}**@{},
    <-0.38mm,-0.2mm>*{};<-2.8mm,2.5mm>*{^1}**@{},
    \end{xy}
\Ea=0.
\Ea
\right.
$$
Thus
$$
(\LB^\circlearrowright)^\Koz(m,n)\simeq
(\LB^\circlearrowright)^!(m,n)= sgn_m\ot sgn_n [m+n-2]=\mbox{span}
\left\langle
\begin{xy}
<0mm,2mm>*{\circ};
<-0.26mm,2.26mm>*{};<-2.76mm,4.76mm>*{}**@{-},
<0.26mm,2.26mm>*{};<3mm,5mm>*{}**@{-},
 <-3mm,5mm>*{\circ};
<-3.26mm,5.26mm>*{};<-5.76mm,7.76mm>*{}**@{-},
<-2.74mm,5.26mm>*{};<0mm,8mm>*{}**@{-},
<-6.4mm,8.4mm>*{\cdot};
<-7.4mm,9.4mm>*{\cdot};
<-8.7mm,10.7mm>*{\circ};
<-9mm,11mm>*{};<-11.7mm,13.7mm>*{}**@{-},
<-8.4mm,11mm>*{};<-5.7mm,13.7mm>*{}**@{-},
<0mm,-1.6mm>*{};<0mm,1.6mm>*{}**@{-},
<0mm,-2mm>*{\circ};
<-0.26mm,-2.26mm>*{};<-2.76mm,-4.76mm>*{}**@{-},
<0.26mm,-2.26mm>*{};<3mm,-5mm>*{}**@{-},
 <-3mm,-5mm>*{\circ};
<-3.26mm,-5.26mm>*{};<-5.76mm,-7.76mm>*{}**@{-},
<-2.74mm,-5.26mm>*{};<0mm,-8mm>*{}**@{-},
<-6.4mm,-8.4mm>*{\cdot};
<-7.4mm,-9.4mm>*{\cdot};
<-8.7mm,-10.7mm>*{\circ};
<-9mm,-11mm>*{};<-11.7mm,-13.7mm>*{}**@{-},
<-8.4mm,-11mm>*{};<-5.7mm,-13.7mm>*{}**@{-},
\end{xy}
\ \ \
\right\rangle,
$$
and, in accordance with  \S 2.6,  the dg free wheeled prop $\LB_\infty^\circlearrowright:=B^c((\LB^\circlearrowright)^\Koz)$
is generated  by the $\bS$-bimodule,
\Beq\label{Generators for LieB_infty}
\mathsf w (\LB^\circlearrowright)^\Koz(m,n)= sgn_m\ot \id_n[m-2]=\mbox{span}\left\langle
\begin{xy}
 <0mm,0mm>*{\circ};<0mm,0mm>*{}**@{},
 <-0.6mm,0.44mm>*{};<-8mm,5mm>*{}**@{-},
 <-0.4mm,0.7mm>*{};<-4.5mm,5mm>*{}**@{-},
 <0mm,0mm>*{};<-1mm,5mm>*{\ldots}**@{},
 <0.4mm,0.7mm>*{};<4.5mm,5mm>*{}**@{-},
 <0.6mm,0.44mm>*{};<8mm,5mm>*{}**@{-},
   <0mm,0mm>*{};<-8.5mm,5.5mm>*{^1}**@{},
   <0mm,0mm>*{};<-5mm,5.5mm>*{^2}**@{},
   <0mm,0mm>*{};<4.5mm,5.5mm>*{^{m\hspace{-0.5mm}-\hspace{-0.5mm}1}}**@{},
   <0mm,0mm>*{};<9.0mm,5.5mm>*{^m}**@{},
 <-0.6mm,-0.44mm>*{};<-8mm,-5mm>*{}**@{-},
 <-0.4mm,-0.7mm>*{};<-4.5mm,-5mm>*{}**@{-},
 <0mm,0mm>*{};<-1mm,-5mm>*{\ldots}**@{},
 <0.4mm,-0.7mm>*{};<4.5mm,-5mm>*{}**@{-},
 <0.6mm,-0.44mm>*{};<8mm,-5mm>*{}**@{-},
   <0mm,0mm>*{};<-8.5mm,-6.9mm>*{^1}**@{},
   <0mm,0mm>*{};<-5mm,-6.9mm>*{^2}**@{},
   <0mm,0mm>*{};<4.5mm,-6.9mm>*{^{n\hspace{-0.5mm}-\hspace{-0.5mm}1}}**@{},
   <0mm,0mm>*{};<9.0mm,-6.9mm>*{^n}**@{},
 \end{xy}
\right\rangle,\ \ \ m,n\geq 1, m+n\geq 3,
\Eeq
and its differential is given on the generating corollas by (cf.\ \cite{Me1})
\Beq\label{LB_infty}
\delta
\begin{xy}
 <0mm,0mm>*{\circ};<0mm,0mm>*{}**@{},
 <-0.6mm,0.44mm>*{};<-8mm,5mm>*{}**@{-},
 <-0.4mm,0.7mm>*{};<-4.5mm,5mm>*{}**@{-},
 <0mm,0mm>*{};<-1mm,5mm>*{\ldots}**@{},
 <0.4mm,0.7mm>*{};<4.5mm,5mm>*{}**@{-},
 <0.6mm,0.44mm>*{};<8mm,5mm>*{}**@{-},
   <0mm,0mm>*{};<-8.5mm,5.5mm>*{^1}**@{},
   <0mm,0mm>*{};<-5mm,5.5mm>*{^2}**@{},
   <0mm,0mm>*{};<4.5mm,5.5mm>*{^{m\hspace{-0.5mm}-\hspace{-0.5mm}1}}**@{},
   <0mm,0mm>*{};<9.0mm,5.5mm>*{^m}**@{},
 <-0.6mm,-0.44mm>*{};<-8mm,-5mm>*{}**@{-},
 <-0.4mm,-0.7mm>*{};<-4.5mm,-5mm>*{}**@{-},
 <0mm,0mm>*{};<-1mm,-5mm>*{\ldots}**@{},
 <0.4mm,-0.7mm>*{};<4.5mm,-5mm>*{}**@{-},
 <0.6mm,-0.44mm>*{};<8mm,-5mm>*{}**@{-},
   <0mm,0mm>*{};<-8.5mm,-6.9mm>*{^1}**@{},
   <0mm,0mm>*{};<-5mm,-6.9mm>*{^2}**@{},
   <0mm,0mm>*{};<4.5mm,-6.9mm>*{^{n\hspace{-0.5mm}-\hspace{-0.5mm}1}}**@{},
   <0mm,0mm>*{};<9.0mm,-6.9mm>*{^n}**@{},
 \end{xy}
\ \ = \ \
 \sum_{[1,\ldots,m]=I_1\sqcup I_2\atop
 {|I_1|\geq 0, |I_2|\geq 1}}
 \sum_{[1,\ldots,n]=J_1\sqcup J_2\atop
 {|J_1|\geq 1, |J_2|\geq 1}
}\hspace{0mm}
(-1)^{\sigma(I_1\sqcup I_2) + |I_1|(|I_2|+1)}
 \begin{xy}
 <0mm,0mm>*{\circ};<0mm,0mm>*{}**@{},
 <-0.6mm,0.44mm>*{};<-8mm,5mm>*{}**@{-},
 <-0.4mm,0.7mm>*{};<-4.5mm,5mm>*{}**@{-},
 <0mm,0mm>*{};<0mm,5mm>*{\ldots}**@{},
 <0.4mm,0.7mm>*{};<4.5mm,5mm>*{}**@{-},
 <0.6mm,0.44mm>*{};<12.4mm,4.8mm>*{}**@{-},
     <0mm,0mm>*{};<-2mm,7mm>*{\overbrace{\ \ \ \ \ \ \ \ \ \ \ \ }}**@{},
     <0mm,0mm>*{};<-2mm,9mm>*{^{I_1}}**@{},
 <-0.6mm,-0.44mm>*{};<-8mm,-5mm>*{}**@{-},
 <-0.4mm,-0.7mm>*{};<-4.5mm,-5mm>*{}**@{-},
 <0mm,0mm>*{};<-1mm,-5mm>*{\ldots}**@{},
 <0.4mm,-0.7mm>*{};<4.5mm,-5mm>*{}**@{-},
 <0.6mm,-0.44mm>*{};<8mm,-5mm>*{}**@{-},
      <0mm,0mm>*{};<0mm,-7mm>*{\underbrace{\ \ \ \ \ \ \ \ \ \ \ \ \ \ \
      }}**@{},
      <0mm,0mm>*{};<0mm,-10.6mm>*{_{J_1}}**@{},
 <13mm,5mm>*{};<13mm,5mm>*{\circ}**@{},
 <12.6mm,5.44mm>*{};<5mm,10mm>*{}**@{-},
 <12.6mm,5.7mm>*{};<8.5mm,10mm>*{}**@{-},
 <13mm,5mm>*{};<13mm,10mm>*{\ldots}**@{},
 <13.4mm,5.7mm>*{};<16.5mm,10mm>*{}**@{-},
 <13.6mm,5.44mm>*{};<20mm,10mm>*{}**@{-},
      <13mm,5mm>*{};<13mm,12mm>*{\overbrace{\ \ \ \ \ \ \ \ \ \ \ \ \ \ }}**@{},
      <13mm,5mm>*{};<13mm,14mm>*{^{I_2}}**@{},
 <12.4mm,4.3mm>*{};<8mm,0mm>*{}**@{-},
 <12.6mm,4.3mm>*{};<12mm,0mm>*{\ldots}**@{},
 <13.4mm,4.5mm>*{};<16.5mm,0mm>*{}**@{-},
 <13.6mm,4.8mm>*{};<20mm,0mm>*{}**@{-},
     <13mm,5mm>*{};<14.3mm,-2mm>*{\underbrace{\ \ \ \ \ \ \ \ \ \ \ }}**@{},
     <13mm,5mm>*{};<14.3mm,-4.5mm>*{_{J_2}}**@{},
 \end{xy}
\Eeq
where $\sigma(I_1\sqcup I_2)$ is the signs of the shuffle
$[1,\ldots,m]=I_1\sqcup I_2$. It is easy to see
that representations
of  $\LB_\infty^\circlearrowright$-algebras in a finite-dimensional vector space $V$ are in one-to-one
correspondence with graded pointed formal Poisson structures on $V$, that is,
total degree degree 2 polyvector fields, $\pi\in \wedge^{\bullet\geq 1} \cT_V$,
which satisfy the Schouten equations $[\pi, \pi]_S=0$ and vanish at the distinguished point $0\in V$
(cf.\ \cite{Me1,Me-lec} and \S 4.3 below).

\mip

{\bf 4.1.3. Non-Koszulnes of $\LB^\circlearrowright$}. Let $\LB_\infty$ be a subcomplex
of the complex $\LB^\circlearrowright_\infty$ spanned by graphs with no closed directed paths,
i.e with no wheels. This subset has an obvious structure of an ordinary prop and, in fact, is a minimal
resolution of the ordinary  prop, $\LB$, of Lie 1-bialgebras (which is defined by the same generators
(\ref{Generators of LieB}) and relations (\ref{R for LieB}) as $\LB^\circlearrowright$ but in the category
of ordinary props). The natural epimorphism,
$$
\pi:(\LB_\infty, \delta) \lon (\LB, 0)
$$
which sends to zero all generating $(m,n)$-corollas (\ref{Generators for LieB_infty}) except
those with $m+n=3$, is a quasi-isomorphism \cite{Me1, Me-graphs}. This means that the prop $\LB$
is Koszul in the category of {\em ordinary}\, props. The wheelification functor
from the category of ordinary props to the category of wheeled props \cite{MMS} sends
these two props into precisely
$\LB^\circlearrowright$ and $\LB_\infty^\circlearrowright$, and the above morphism $\pi$
into the associated
morphism of dg wheeled props,
$$
\pi^\circlearrowright: (\LB^\circlearrowright_\infty, \delta) \lon (\LB^\circlearrowright, 0).
$$
The morphism $\pi^\circlearrowright$ is {\em not}, however, a
quasi-isomorphism: the following element \cite{Me-graphs}
\Beq\label{Three_wheels}
 \begin{xy}
<-5mm,5mm>*{\circ};
<-5mm,5mm>*{};<-5mm,8mm>*{}**@{-},
<-5mm,5mm>*{};<-7mm,4mm>*{}**@{-},
 <0mm,0mm>*{\circ};
<0mm,0mm>*{};<-5mm,-5mm>*{}**@{-},
 <0mm,0mm>*{};<-5mm,5mm>*{}**@{-},
<0mm,0mm>*{};<1.5mm,1.5mm>*{}**@{-},
 <0mm,0mm>*{};<1.5mm,-1.5mm>*{}**@{-};
<-5mm,-5mm>*{\circ};
<-5mm,-5mm>*{};<-7mm,-2mm>*{}**@{-},
<-5mm,-5mm>*{};<-5mm,-8mm>*{}**@{-},
   \ar@{->}@(ul,dl) (-5.0,8.0)*{};(-7.0,4.0)*{},
   \ar@{->}@(ur,dr) (1.5,1.5)*{};(1.5,-1.5)*{},
   \ar@{->}@(ul,dl) (-7.0,-2.0)*{};(-5.0,-8.0)*{},
\end{xy}
\ - \
 \begin{xy}
<0mm,4mm>*{\circ};
<0mm,4mm>*{};<0mm,8mm>*{}**@{-},
<0mm,4mm>*{};<3mm,2mm>*{}**@{-},
 <0mm,-5mm>*{\circ};
<0mm,-5mm>*{};<0mm,4mm>*{}**@{-},
<0mm,-5mm>*{};<2mm,-3mm>*{}**@{-},
<0mm,-5mm>*{};<0mm,-8mm>*{}**@{-},
<-5mm,-5mm>*{};<0mm,4mm>*{}**@{-};
<-5mm,-5mm>*{\circ};
<-5mm,-5mm>*{};<-7mm,-2mm>*{}**@{-},
<-5mm,-5mm>*{};<-5mm,-8mm>*{}**@{-},
   \ar@{->}@(ur,dr) (0,8.0)*{};(3.0,2.0)*{},
   \ar@{->}@(ur,dr) (2.0,-3.0)*{};(0.0,-8.0)*{},
   \ar@{->}@(ul,dl) (-7.0,-2.0)*{};(-5.0,-8.0)*{},
\end{xy}
\ + \
 \begin{xy}
<0mm,-4mm>*{\circ};
<0mm,-4mm>*{};<0mm,-8mm>*{}**@{-},
<0mm,-4mm>*{};<3mm,-2mm>*{}**@{-},
 <0mm,5mm>*{\circ};
<0mm,5mm>*{};<0mm,-4mm>*{}**@{-},
<0mm,5mm>*{};<2mm,3mm>*{}**@{-},
<0mm,5mm>*{};<0mm,8mm>*{}**@{-},
<-5mm,5mm>*{};<0mm,-4mm>*{}**@{-};
<-5mm,5mm>*{\circ};
<-5mm,5mm>*{};<-7mm,2mm>*{}**@{-},
<-5mm,5mm>*{};<-5mm,8mm>*{}**@{-},
   \ar@{->}@(ur,dr) (3.0,-2.0)*{};(0.0,-8.0)*{},
   \ar@{->}@(ur,dr) (0.0,8.0)*{};(2.0,3.0)*{},
   \ar@{->}@(ul,dl) (-5.0,8.0)*{};(-7.0,2.0)*{},
\end{xy}
\in \LB_\infty^\circlearrowright
\Eeq
gives a non-trivial cohomology class in $H(\LB_\infty^\circlearrowright, \delta)$ which is,
however, sent to zero under $\pi^\circlearrowright$. This means that the wheeled prop
of
Lie 1-bialgebras is {\em not}\ Koszul, and its minimal resolution, $(\LB^\circlearrowright)_\infty$ is
larger than $\LB_\infty^\circlearrowright$. Representations of $(\LB^\circlearrowright)_\infty$
in a vector space $V$ are called formal {\em wheeled Poisson structures}\,; these (at present mysterious)
structures are Maurer-Cartan elements
of a certain $L_\infty$ algebra\footnote{Graph (\ref{Three_wheels}) gives, in fact, an explicit
formula for particular $\mu_3$ composition in that $L_\infty$ algebra.}
which, in accordance with the general theory of \cite{MV}, is canonically associated to
$(\LB^\circlearrowright)_\infty$ and which involve not only Schouten brackets
but also divergence operators;  it was proven in \cite{Me-lec}
that wheeled Poisson structures can be deformation quantized over $\Q$.

\mip

{\bf 4.2. Wheeled prop, $\ULB$, of unimodular Lie 1-bialgebras}. A finite dimensional Lie 1-bialgebra
$V$ is called unimodular if, for any $e\in V$ and $e^*\in V^*$, the supertraces of linear maps,
$$
\Ba{rccc}
Ad_e: & V & \lon & V \\
&        v & \lon & [e\bullet v]
\Ea
\ \ \ \mbox{and}\ \ \
\Ba{rccc}
Ad_{e^*}: & V^* & \lon & V^* \\
&        v^* & \lon & [e^*, v^*]
\Ea,
$$
are zero. Here $[\ ,\ ]$ are the Lie brackets on $V^*$ induced by Lie coalgebra structure on $V$.
The wheeled prop(erad), $\ULB$
of unimodular Lie 1-bialgebras is is a quotient of the free wheeled
prop(erda)  generated by the $\bS$-bimodule (\ref{Generators of LieB}) by the ideal
generated by relations (\ref{R for LieB}) and the following ones,
$$
\begin{xy}
 <0mm,-0.55mm>*{};<0mm,-2.5mm>*{}**@{-},
 <0.5mm,0.5mm>*{};<2.2mm,2.2mm>*{}**@{-},
 <-0.48mm,0.48mm>*{};<-2.2mm,2.2mm>*{}**@{-},
 <0mm,0mm>*{\circ};<0mm,0mm>*{}**@{},
 <-0.48mm,0.48mm>*{};<-2.7mm,2.8mm>*{^1}**@{},
(2.2,2.2)*{}
   \ar@{->}@(ur,dr) (0,-2.5)*{}
 \end{xy}=0\ , \ \ \
\begin{xy}
 <0mm,0.66mm>*{};<0mm,3mm>*{}**@{-},
 <0.39mm,-0.39mm>*{};<2.2mm,-2.2mm>*{}**@{-},
 <-0.35mm,-0.35mm>*{};<-2.2mm,-2.2mm>*{}**@{-},
 <0mm,0mm>*{\circ};<0mm,0mm>*{}**@{},
   <-0.35mm,-0.35mm>*{};<-2.8mm,-4mm>*{^1}**@{},
(0.0,3.0)*{}
   \ar@{->}@(ur,dr) (2.2,-2.2)*{}
\end{xy}=0.
$$
Hence the Koszul dual properad, $(\ULB)^!$, is a quadratic wheeled properad
generated by the $\bS$-bimodule (\ref{Generators ULieB-vee}) modulo the relations,
$$
\Ba{c}
\begin{xy}
 <0mm,0mm>*{\circ};<0mm,0mm>*{}**@{},
 <0mm,-0.49mm>*{};<0mm,-3.0mm>*{}**@{-},
 <0.49mm,0.49mm>*{};<1.9mm,1.9mm>*{}**@{-},
 <-0.5mm,0.5mm>*{};<-1.9mm,1.9mm>*{}**@{-},
 <-2.3mm,2.3mm>*{\circ};<-2.3mm,2.3mm>*{}**@{},
 <-1.8mm,2.8mm>*{};<0mm,4.9mm>*{}**@{-},
 <-2.8mm,2.9mm>*{};<-4.6mm,4.9mm>*{}**@{-},
   <0.49mm,0.49mm>*{};<2.7mm,2.3mm>*{^3}**@{},
   <-1.8mm,2.8mm>*{};<0.4mm,5.3mm>*{^2}**@{},
   <-2.8mm,2.9mm>*{};<-5.1mm,5.3mm>*{^1}**@{},
 \end{xy}
\ - \
\begin{xy}
 <0mm,0mm>*{\circ};<0mm,0mm>*{}**@{},
 <0mm,-0.49mm>*{};<0mm,-3.0mm>*{}**@{-},
 <0.49mm,0.49mm>*{};<1.9mm,1.9mm>*{}**@{-},
 <-0.5mm,0.5mm>*{};<-1.9mm,1.9mm>*{}**@{-},
 <-2.3mm,2.3mm>*{\circ};<-2.3mm,2.3mm>*{}**@{},
 <-1.8mm,2.8mm>*{};<0mm,4.9mm>*{}**@{-},
 <-2.8mm,2.9mm>*{};<-4.6mm,4.9mm>*{}**@{-},
   <0.49mm,0.49mm>*{};<2.7mm,2.3mm>*{^1}**@{},
   <-1.8mm,2.8mm>*{};<0.4mm,5.3mm>*{^3}**@{},
   <-2.8mm,2.9mm>*{};<-5.1mm,5.3mm>*{^2}**@{},
 \end{xy}
\Ea
=0, \ \
\Ba{c}
 \begin{xy}
 <0mm,0mm>*{\circ};<0mm,0mm>*{}**@{},
 <0mm,0.69mm>*{};<0mm,3.0mm>*{}**@{-},
 <0.39mm,-0.39mm>*{};<2.4mm,-2.4mm>*{}**@{-},
 <-0.35mm,-0.35mm>*{};<-1.9mm,-1.9mm>*{}**@{-},
 <-2.4mm,-2.4mm>*{\circ};<-2.4mm,-2.4mm>*{}**@{},
 <-2.0mm,-2.8mm>*{};<0mm,-4.9mm>*{}**@{-},
 <-2.8mm,-2.9mm>*{};<-4.7mm,-4.9mm>*{}**@{-},
    <0.39mm,-0.39mm>*{};<3.3mm,-4.0mm>*{^3}**@{},
    <-2.0mm,-2.8mm>*{};<0.5mm,-6.7mm>*{^2}**@{},
    <-2.8mm,-2.9mm>*{};<-5.2mm,-6.7mm>*{^1}**@{},
 \end{xy}
\ - \
 \begin{xy}
 <0mm,0mm>*{\circ};<0mm,0mm>*{}**@{},
 <0mm,0.69mm>*{};<0mm,3.0mm>*{}**@{-},
 <0.39mm,-0.39mm>*{};<2.4mm,-2.4mm>*{}**@{-},
 <-0.35mm,-0.35mm>*{};<-1.9mm,-1.9mm>*{}**@{-},
 <-2.4mm,-2.4mm>*{\circ};<-2.4mm,-2.4mm>*{}**@{},
 <-2.0mm,-2.8mm>*{};<0mm,-4.9mm>*{}**@{-},
 <-2.8mm,-2.9mm>*{};<-4.7mm,-4.9mm>*{}**@{-},
    <0.39mm,-0.39mm>*{};<3.3mm,-4.0mm>*{^1}**@{},
    <-2.0mm,-2.8mm>*{};<0.5mm,-6.7mm>*{^3}**@{},
    <-2.8mm,-2.9mm>*{};<-5.2mm,-6.7mm>*{^2}**@{},
 \end{xy}
\Ea=0, \ \
\Ba{c}
 \begin{xy}
 <0mm,2.47mm>*{};<0mm,0.12mm>*{}**@{-},
 <0.5mm,3.5mm>*{};<2.2mm,5.2mm>*{}**@{-},
 <-0.48mm,3.48mm>*{};<-2.2mm,5.2mm>*{}**@{-},
 <0mm,3mm>*{\circ};<0mm,3mm>*{}**@{},
  <0mm,-0.8mm>*{\circ};<0mm,-0.8mm>*{}**@{},
<-0.39mm,-1.2mm>*{};<-2.2mm,-3.5mm>*{}**@{-},
 <0.39mm,-1.2mm>*{};<2.2mm,-3.5mm>*{}**@{-},
     <0.5mm,3.5mm>*{};<2.8mm,5.7mm>*{^2}**@{},
     <-0.48mm,3.48mm>*{};<-2.8mm,5.7mm>*{^1}**@{},
   <0mm,-0.8mm>*{};<-2.7mm,-5.2mm>*{^1}**@{},
   <0mm,-0.8mm>*{};<2.7mm,-5.2mm>*{^2}**@{},
\end{xy}
\  + \
\begin{xy}
 <0mm,-1.3mm>*{};<0mm,-3.5mm>*{}**@{-},
 <0.38mm,-0.2mm>*{};<2.0mm,2.0mm>*{}**@{-},
 <-0.38mm,-0.2mm>*{};<-2.2mm,2.2mm>*{}**@{-},
<0mm,-0.8mm>*{\circ};<0mm,0.8mm>*{}**@{},
 <2.4mm,2.4mm>*{\circ};<2.4mm,2.4mm>*{}**@{},
 <2.77mm,2.0mm>*{};<4.4mm,-0.8mm>*{}**@{-},
 <2.4mm,3mm>*{};<2.4mm,5.2mm>*{}**@{-},
     <0mm,-1.3mm>*{};<0mm,-5.3mm>*{^1}**@{},
     <2.5mm,2.3mm>*{};<5.1mm,-2.6mm>*{^2}**@{},
    <2.4mm,2.5mm>*{};<2.4mm,5.7mm>*{^2}**@{},
    <-0.38mm,-0.2mm>*{};<-2.8mm,2.5mm>*{^1}**@{},
    \end{xy}
\Ea=0.
$$
Therefore,
$$
(\ULB)^!(m,n)=
\bigoplus_{a=0}^\infty
{sgn_m\ot sgn_n [m+n-2-2a]} =
\mbox{span}
\underbrace{
\left\langle
\begin{xy}
<0mm,2mm>*{\circ};
<-0.26mm,2.26mm>*{};<-2.76mm,4.76mm>*{}**@{-},
<0.26mm,2.26mm>*{};<3mm,5mm>*{}**@{-},
 <-3mm,5mm>*{\circ};
<-3.26mm,5.26mm>*{};<-5.76mm,7.76mm>*{}**@{-},
<-2.74mm,5.26mm>*{};<0mm,8mm>*{}**@{-},
<-6.4mm,8.4mm>*{\cdot};
<-7.4mm,9.4mm>*{\cdot};
<-8.7mm,10.7mm>*{\circ};
<-9mm,11mm>*{};<-11.7mm,13.7mm>*{}**@{-},
<-8.4mm,11mm>*{};<-5.7mm,13.7mm>*{}**@{-},
<0mm,-1.6mm>*{};<0mm,1.6mm>*{}**@{-},
<0mm,-2mm>*{\circ};
<-0.26mm,-2.26mm>*{};<-2.76mm,-4.76mm>*{}**@{-},
<0.26mm,-2.26mm>*{};<3mm,-5mm>*{}**@{-},
 <-3mm,-5mm>*{\circ};
<-3.26mm,-5.26mm>*{};<-5.76mm,-7.76mm>*{}**@{-},
<-2.74mm,-5.26mm>*{};<0mm,-8mm>*{}**@{-},
<-6.4mm,-8.4mm>*{\cdot};
<-7.4mm,-9.4mm>*{\cdot};
<-8.7mm,-10.7mm>*{\circ};
<-9mm,-11mm>*{};<-11.7mm,-13.7mm>*{}**@{-},
<-8.4mm,-11mm>*{};<-5.7mm,-13.7mm>*{}**@{-},
<0mm,14mm>*{\cdot};
<-1mm,15mm>*{\cdot};
<-2mm,16mm>*{\cdot};
<3.3mm,5.3mm>*{\circ};
<0.3mm,8.3mm>*{\circ};
   \ar@{->}@(u,r) (0,9)*{};(0.3,8.0)*{},
   \ar@{->}@(ur,dr) (3.6,5.6)*{};(3.0,5.0)*{};
\end{xy}
\ \ \
\right\rangle}_{m\ \mathrm o  \mathrm u   \mathrm t\ \mathrm l
\mathrm e \mathrm g \mathrm s, \
n\ \mathrm i  \mathrm n  \ \mathrm l
\mathrm e \mathrm g \mathrm s,\
a\  \mathrm l  \mathrm o   \mathrm o\mathrm p
\mathrm s
} ,
$$
Note that the graph on the r.h.s.\ above is zero unless $\Z_{\geq 0}$-valued parameters
$m$, $n$ and $a$ satisfy inequalities,
$$
m+n+2a\geq 3, m+a\geq 1, n+a\geq 1.
$$
Hence
the dg free wheeled prop $\ULB_\infty:=B^c((\ULB)^\Koz)$
is generated  by an $\bS$-bimodule,
\Beq\label{Generators in ULieB_infty}
\mathsf w (\ULB)^\Koz(m,n)= \bigoplus_{a\geq 0}^\infty
sgn_m\ot \id_n[m-2-2a]=\mbox{span}\left\langle
\begin{xy}
 <0mm,0mm>*{\mbox{$\xy *=<3mm,3mm>
\txt{{{$_a$}}}*\frm{-}\endxy$}};
<-1.5mm,1.5mm>*{};<-5mm,5mm>*{}**@{-},
<-1mm,1.5mm>*{};<-2.7mm,5mm>*{}**@{-},
<1.5mm,1.5mm>*{};<5mm,5mm>*{}**@{-},
<1mm,1.5mm>*{};<2.7mm,5mm>*{}**@{-},
<0mm,4.5mm>*{...};
<-1.5mm,-1.5mm>*{};<-5mm,-5mm>*{}**@{-},
<-1mm,-1.5mm>*{};<-2.7mm,-5mm>*{}**@{-},
<1.5mm,-1.5mm>*{};<5mm,-5mm>*{}**@{-},
<1mm,-1.5mm>*{};<2.7mm,-5mm>*{}**@{-},
<0mm,-4.5mm>*{...};
<0mm,0mm>*{};<-5.8mm,6.3mm>*{^1}**@{},
   <0mm,0mm>*{};<-3mm,6.3mm>*{^2}**@{},
   <0mm,0mm>*{};<2.6mm,6.3mm>*{^{m\hspace{-0.5mm}-\hspace{-0.5mm}1}}**@{},
   <0mm,0mm>*{};<7.0mm,6.3mm>*{^m}**@{},
<0mm,0mm>*{};<-5.8mm,-6.3mm>*{_1}**@{},
   <0mm,0mm>*{};<-3mm,-6.3mm>*{_2}**@{},
   <0mm,0mm>*{};<3mm,-6.3mm>*{_{n\hspace{-0.5mm}-\hspace{-0.5mm}1}}**@{},
   <0mm,0mm>*{};<7.0mm,-6.3mm>*{_n}**@{},
 \end{xy}
\right\rangle_{m+n+2a\geq 3\atop m+a\geq 1, n+a\geq 1}.
\Eeq
Definition of the cobar construction given in \S 2.5 gives, after straightforward computations,
 the following
formula for the differential in $\ULB_\infty$,
\Beq\label{Differential in ULieB_infty}
\delta
\begin{xy}
 <0mm,0mm>*{\mbox{$\xy *=<3mm,3mm>
\txt{{{$_a$}}}*\frm{-}\endxy$}};
<-1.5mm,1.5mm>*{};<-5mm,5mm>*{}**@{-},
<-1mm,1.5mm>*{};<-2.7mm,5mm>*{}**@{-},
<1.5mm,1.5mm>*{};<5mm,5mm>*{}**@{-},
<1mm,1.5mm>*{};<2.7mm,5mm>*{}**@{-},
<0mm,4.5mm>*{...};
<-1.5mm,-1.5mm>*{};<-5mm,-5mm>*{}**@{-},
<-1mm,-1.5mm>*{};<-2.7mm,-5mm>*{}**@{-},
<1.5mm,-1.5mm>*{};<5mm,-5mm>*{}**@{-},
<1mm,-1.5mm>*{};<2.7mm,-5mm>*{}**@{-},
<0mm,-4.5mm>*{...};
<0mm,0mm>*{};<-5.8mm,6.3mm>*{^1}**@{},
   <0mm,0mm>*{};<-3mm,6.3mm>*{^2}**@{},
   <0mm,0mm>*{};<2.6mm,6.3mm>*{^{m\hspace{-0.5mm}-\hspace{-0.5mm}1}}**@{},
   <0mm,0mm>*{};<7.0mm,6.3mm>*{^m}**@{},
<0mm,0mm>*{};<-5.8mm,-6.3mm>*{_1}**@{},
   <0mm,0mm>*{};<-3mm,-6.3mm>*{_2}**@{},
   <0mm,0mm>*{};<3mm,-6.3mm>*{_{n\hspace{-0.5mm}-\hspace{-0.5mm}1}}**@{},
   <0mm,0mm>*{};<7.0mm,-6.3mm>*{_n}**@{},
 \end{xy}
 =(-1)^{m-1}
\begin{xy}
 <0mm,0mm>*{\mbox{$\xy *=<6mm,3mm>
\txt{{{$_{a-1}$}}}*\frm{-}\endxy$}};
<-1.5mm,1.5mm>*{};<-5mm,5mm>*{}**@{-},
<-1mm,1.5mm>*{};<-2.7mm,5mm>*{}**@{-},
<1.5mm,1.5mm>*{};<4mm,4mm>*{}**@{-},
<1mm,1.5mm>*{};<2.7mm,5mm>*{}**@{-},
<0mm,4.5mm>*{...};
<-1.5mm,-1.5mm>*{};<-5mm,-5mm>*{}**@{-},
<-1mm,-1.5mm>*{};<-2.7mm,-5mm>*{}**@{-},
<1.5mm,-1.5mm>*{};<4mm,-4mm>*{}**@{-},
<1mm,-1.5mm>*{};<2.7mm,-5mm>*{}**@{-},
<0mm,-4.5mm>*{...};
<0mm,0mm>*{};<-5.8mm,6.3mm>*{^1}**@{},
   <0mm,0mm>*{};<-3mm,6.3mm>*{^2}**@{},
   <0mm,0mm>*{};<3.0mm,6.3mm>*{^m}**@{},
<0mm,0mm>*{};<-5.8mm,-6.3mm>*{_1}**@{},
   <0mm,0mm>*{};<-3mm,-6.3mm>*{_2}**@{},
   <0mm,0mm>*{};<3.0mm,-6.3mm>*{_n}**@{},
   \ar@{->}@(ur,dr) (4.0,4.0)*{};(4.0,-4.0)*{},
 \end{xy}
 \ + \sum_{a=b+c\atop b,c\geq 0}\sum_{m=I'\sqcup I''\atop
 [n]=J'\sqcup J''}
 (-1)^{\sigma(I_1\sqcup I_2) + |I_1|(|I_2|+1)}
 \Ba{c}
 \begin{xy}
 <0mm,0mm>*{\mbox{$\xy *=<3mm,3mm>
\txt{{{$_b$}}}*\frm{-}\endxy$}};
<-1.5mm,1.5mm>*{};<-5mm,5mm>*{}**@{-},
<-1mm,1.5mm>*{};<-2.7mm,5mm>*{}**@{-},
<1.5mm,1.5mm>*{};<8mm,8mm>*{}**@{-},
<1mm,1.5mm>*{};<2.7mm,5mm>*{}**@{-},
<0mm,4.5mm>*{...};
<-1.5mm,-1.5mm>*{};<-5mm,-5mm>*{}**@{-},
<-1mm,-1.5mm>*{};<-2.7mm,-5mm>*{}**@{-},
<1.5mm,-1.5mm>*{};<5mm,-5mm>*{}**@{-},
<1mm,-1.5mm>*{};<2.7mm,-5mm>*{}**@{-},
<0mm,-4.5mm>*{...};
 <9.5mm,9.5mm>*{\mbox{$\xy *=<3mm,3mm>
\txt{{{$_c$}}}*\frm{-}\endxy$}};
<8mm,11mm>*{};<5mm,14mm>*{}**@{-},
<8.5mm,11mm>*{};<7mm,14mm>*{}**@{-},
<10.5mm,11mm>*{};<12mm,14mm>*{}**@{-},
<11mm,11mm>*{};<14mm,14mm>*{}**@{-},
<9.5mm,13.5mm>*{...};
<8.5mm,8mm>*{};<7mm,4.5mm>*{}**@{-},
<10.5mm,8mm>*{};<12mm,4.5mm>*{}**@{-},
<11mm,8mm>*{};<14mm,4.5mm>*{}**@{-},
<9.5mm,5mm>*{...};
<10.5mm,1.8mm>*{\underbrace{\ \ \ \ \ \ \  }_{J''}};
<0mm,-8mm>*{\underbrace{\ \ \ \ \ \ \ \ \  }_{J'}};
<-1mm,8.2mm>*{\overbrace{\ \ \ \ \ \ \  }^{I'}};
<9.5mm,17mm>*{\overbrace{\ \ \ \ \ \ \  \ }^{I''}};
 \end{xy}
 \Ea
 \Eeq
where $\sigma(I_1\sqcup I_2)$ is the signs of the shuffle
$[1,\ldots,m]=I_1\sqcup I_2$.

\mip
{\bf 4.3. Representations of $\ULB_\infty$ and quantum BV manifolds}. Let $(V, d)$ be a
finite-dimensional dg vector space, and $\cM_{V^*}^\hbar$ the formal $\hbar$-twisted odd symplectic manifold corresponding
to the graded commutative ring $\widehat{\odot^\bullet} (V^*\oplus V[-1])[[\hbar]]$,
$\hbar$ being the formal variable of degree $2$ (see \S 3.9).

\sip

An arbitrary morphism $\rho: \ULB_\infty \rar \cE nd_V$ is uniquely defined
by its values on the generators,
$$
\rho^{(a)}_{m,n}:=\rho\left(
\begin{xy}
 <0mm,0mm>*{\mbox{$\xy *=<3mm,3mm>
\txt{{{$_a$}}}*\frm{-}\endxy$}};
<-1.5mm,1.5mm>*{};<-5mm,5mm>*{}**@{-},
<-1mm,1.5mm>*{};<-2.7mm,5mm>*{}**@{-},
<1.5mm,1.5mm>*{};<5mm,5mm>*{}**@{-},
<1mm,1.5mm>*{};<2.7mm,5mm>*{}**@{-},
<0mm,4.5mm>*{...};
<-1.5mm,-1.5mm>*{};<-5mm,-5mm>*{}**@{-},
<-1mm,-1.5mm>*{};<-2.7mm,-5mm>*{}**@{-},
<1.5mm,-1.5mm>*{};<5mm,-5mm>*{}**@{-},
<1mm,-1.5mm>*{};<2.7mm,-5mm>*{}**@{-},
<0mm,-4.5mm>*{...};
<0mm,0mm>*{};<-5.8mm,6.3mm>*{^1}**@{},
   <0mm,0mm>*{};<-3mm,6.3mm>*{^2}**@{},
   <0mm,0mm>*{};<2.6mm,6.3mm>*{^{m\hspace{-0.5mm}-\hspace{-0.5mm}1}}**@{},
   <0mm,0mm>*{};<7.0mm,6.3mm>*{^m}**@{},
<0mm,0mm>*{};<-5.8mm,-6.3mm>*{_1}**@{},
   <0mm,0mm>*{};<-3mm,-6.3mm>*{_2}**@{},
   <0mm,0mm>*{};<3mm,-6.3mm>*{_{n\hspace{-0.5mm}-\hspace{-0.5mm}1}}**@{},
   <0mm,0mm>*{};<7.0mm,-6.3mm>*{_n}**@{},
 \end{xy}
\right)\in \Hom(\odot^n V, \wedge^m V[2-m-2a])=\odot^n V^*\ot \wedge^m V[2-m-2a])
\subset \cE nd_V(m,n).
$$
We assemble the collection of linear maps, $\{\rho^{(a)}_{m,n}\}$, into one ``generating"
degree  0 function,
$$
\bGa:=\sum_{a,m,n\geq 0} \rho_{m,n}^{(a)}\hbar^a \in \f_{\cM_V}.
$$

Let $\{e_a\}$ be an arbitrary basis in $V$, and $\{x^a, \psi_a\}$ the associated basis
in $V^*\oplus V[-1]$, $|\psi_a|=-|x^a|-1=|e_a|-1$. Then
\Beqrn
\rho^{(a)}_{m,n}(e_{b_1}, \ldots, e_{b_n})&=&
\sum_{a_1,\dots, a_m}\mu_{b_1\ldots b_n}^{\al_1\ldots \al_m}\psi_{a_1}\ldots
\psi_{a_m},\\
d(e_b)&=& \sum_b d_{b}^a \psi_a,
\Eeqrn
for some $\mu_{b_1\ldots b_n}^{\al_1\ldots \al_m}\in \K$, $d_b^a\in \K$, and we set
$$
\Ga := d+ \bGa =
 \sum_{a,b}d_b^a x^b \psi_b  + \sum_{m+n\geq 3}\frac{1}{m!n!}\sum_{{a_1,\ldots, a_m}
\atop {b_1,\ldots, b_n}} \mu^{a_1\ldots a_m}_{b_1\ldots b_n} x^{b_1}\cdots x^{b_n}
\psi_{a_1}\ldots\psi_{a_m}\in \f_{\cM_V}.
$$
It is a straightforward calculation to check using formula
(\ref{Differential in ULieB_infty})  that
the compatibility of the morphism $\rho$ with the differentials,
$$
\rho\circ \delta = d\circ \rho
$$
 is equivalent to the equation,
$$
\Delta_0 \Gamma + \frac{1}{2}[\Ga\bullet \Ga]_S=0,
$$
where  $\Delta_0=\sum\p^2/ \p x^a\p\psi_a$ and  $[\ \bullet \ ]_S$ stand for the odd Poisson
brackets on $\cM_{V^*}^\hbar$. The function $\Gamma$ satisfies the boundary conditions
(\ref{Boundary conditions for Gamma})
of the definition
\S 3.9.2. Hence we proven the following

\mip

\noindent{\bf 4.3.1. Proposition.} {\em There is a one-to-one correspondence between representations
of the dg wheeled prop $(\ULB_\infty, \delta)$ in a dg vector space $V$ and
quantum BV structures on the formal odd symplectic manifold manifold $\cM^\hbar_{V^*}$.}

\mip

Thus the category of quantum BV manifolds is controlled by a surprisingly simple
quadratic wheeled prop, $\ULB$, of unimodular Lie 1-bialgebras.

\mip

\noindent{\bf 4.3.2. Remark.} We do not know at present whether or not
the wheeled properad $\ULB$ is Koszul, i.e.\ whether or not the natural epimorphism,
$$
\pi^\circlearrowright: (\ULB_\infty, \delta) \lon (\ULB, 0),
$$
is a quasi-isomorphism. Our study of the category of quantum BV manifolds in \S 3
was partly motivated by this open problem. If it is Koszul, then {\em unimodular}\,
Poisson structures can be deformation quantized over $\Q$ with the help of the wheeled prop quantization machine
developed in \cite{Me-lec}.

\bip
\bip
\section{Wheeled dg prop of unimodular Poisson structures}

\mip

{\bf 5.1 Modular volume form}. Let $M$ be a $\Z$-graded manifold. A Poisson structure
on $M$ is a Maurer-Cartan element, $\pi\in \wedge^\bullet \cT_M$, in the Schouten
Lie algebra on $M$, that is, a total degree 2 polyvector field, satisfying the
equation $[\pi\bullet \pi]_S=0$.  If $M$ is concentrated in degree $0$, then $\pi$ must be a
bivector  field, but in general $\pi$ might have non-zero summands lying in $\wedge^n \cT_M$
with $n\neq 2$. Let $\cM$ be be the total space of the bundle, $\Omega^1_M$, of 1-forms
on $M$. Then a polyvector field $\pi$ defines a  function on $\cM$
which we denote by the same letter; the Schouten equations translate into
$\{\pi\bullet \pi \}=0$, where $\{\ \bullet\ \}$ are the odd Poisson brackets
associated with the canonical  odd symplectic structure on $\cM$ (see \S 3.3).
The Poisson structure $\pi$ gives rise to the associated
 degree 1 hamiltonian vector field, $H_\pi$, on $\cM$ which is homological,
i.e.\ $[H_\pi,H_\pi]=H_{\{\pi\bullet\pi\}}=0$.
Any volume form $\nu\in \Ber(M)$, induces, via the canonical isomorphism
$\Ber(\cM)=(\Ber(M))^{\ot 2}$, a volume form on $\cM$ which we denote by $\hat{\nu}$.

\mip

{\bf 5.1.1 Definition} \cite{We}. Let $(M, \pi)$ be a $\Z$-graded Poisson manifold.
A volume form $\nu\in \Ber(M)$ is called {\em modular}\, if the equation,
$$
\caL _{H_\pi}\hat{\nu}=0,
$$
is satisfied.
In this case $\pi$
 is called a {\em unimodular Poisson structure on $(M, \nu)$}.

 \sip

Any vector space $V$ (viewed as a linear formal manifold) admits a translation invariant Berezin volume form, $\nu_0$, which is defined uniquely up to
multiplication by a non-zero constant. A formal Poisson structure $\pi$ on $(V, \nu_0)$ is called
a {\em unimodular Poisson structure on $V$}. If $\{x^a\}$ are linear coordinates on $V$,
then a unimodular Poisson structure on $V$ is given by an ordinary Poisson structure on $V$,
$$
\pi:=\sum_{n\geq 1}\sum_{a_1, \ldots, a_n}
\pi^{a_1,\ldots, a_n}(x)\psi_{a_1}\ldots \psi_{a_n}\in\f_{\Omega^1_V},
$$
with coefficients
$\pi^{a_1,\ldots, a_n}(x)$ satisfying an extra condition
$$
\sum_b \frac{\p \pi^{ba_2,\ldots, a_n(x)}}{\p x^b}=0,\  \ \forall n\geq 1.
$$

{\bf 5.2. Wheeled dg prop of unimodular Poisson structures}. Let $I^\circlearrowright$
be the ideal in the dg wheeled prop
$\ULB_\infty$ (see \S 4.1.2) generated by loops,
$$
I^\circlearrowright:=\left\langle
\begin{xy}
 <0mm,0mm>*{\circ};<0mm,0mm>*{}**@{},
 <-0.6mm,0.44mm>*{};<-8mm,5mm>*{}**@{-},
 <-0.4mm,0.7mm>*{};<-4.5mm,5mm>*{}**@{-},
 <0mm,0mm>*{};<-1mm,5mm>*{\ldots}**@{},
 <0.4mm,0.7mm>*{};<4.5mm,5mm>*{}**@{-},
 <0.6mm,0.44mm>*{};<6mm,4mm>*{}**@{-},
   <0mm,0mm>*{};<-8.5mm,5.5mm>*{^1}**@{},
   <0mm,0mm>*{};<-5mm,5.5mm>*{^2}**@{},
   <0mm,0mm>*{};<4.5mm,5.5mm>*{^m}**@{},
 <-0.6mm,-0.44mm>*{};<-8mm,-5mm>*{}**@{-},
 <-0.4mm,-0.7mm>*{};<-4.5mm,-5mm>*{}**@{-},
 <0mm,0mm>*{};<-1mm,-5mm>*{\ldots}**@{},
 <0.4mm,-0.7mm>*{};<4.5mm,-5mm>*{}**@{-},
 <0.6mm,-0.44mm>*{};<6mm,-4mm>*{}**@{-},
   <0mm,0mm>*{};<-8.5mm,-6.9mm>*{^1}**@{},
   <0mm,0mm>*{};<-5mm,-6.9mm>*{^2}**@{},
   <0mm,0mm>*{};<4.5mm,-6.9mm>*{^{n}}**@{},
   \ar@{->}@(ur,dr) (6.0,4.0)*{};(6.0,-4.0)*{},
 \end{xy}
 \right\rangle.
$$

{\bf 5.2.1. Lemma}. $\delta \begin{xy}
 <0mm,0mm>*{\circ};<0mm,0mm>*{}**@{},
 <-0.6mm,0.44mm>*{};<-8mm,5mm>*{}**@{-},
 <-0.4mm,0.7mm>*{};<-4.5mm,5mm>*{}**@{-},
 <0mm,0mm>*{};<-1mm,5mm>*{\ldots}**@{},
 <0.4mm,0.7mm>*{};<4.5mm,5mm>*{}**@{-},
 <0.6mm,0.44mm>*{};<6mm,4mm>*{}**@{-},
   <0mm,0mm>*{};<-8.5mm,5.5mm>*{^1}**@{},
   <0mm,0mm>*{};<-5mm,5.5mm>*{^2}**@{},
   <0mm,0mm>*{};<4.5mm,5.5mm>*{^m}**@{},
 <-0.6mm,-0.44mm>*{};<-8mm,-5mm>*{}**@{-},
 <-0.4mm,-0.7mm>*{};<-4.5mm,-5mm>*{}**@{-},
 <0mm,0mm>*{};<-1mm,-5mm>*{\ldots}**@{},
 <0.4mm,-0.7mm>*{};<4.5mm,-5mm>*{}**@{-},
 <0.6mm,-0.44mm>*{};<6mm,-4mm>*{}**@{-},
   <0mm,0mm>*{};<-8.5mm,-6.9mm>*{^1}**@{},
   <0mm,0mm>*{};<-5mm,-6.9mm>*{^2}**@{},
   <0mm,0mm>*{};<4.5mm,-6.9mm>*{^{n}}**@{},
   \ar@{->}@(ur,dr) (6.0,4.0)*{};(6.0,-4.0)*{},
 \end{xy} \in I^\circlearrowright$.

Proof is a straightforward calculation based on formula (\ref{LB_infty}).

\sip

Thus $I^\circlearrowright$ is a {\em dg}\, ideal in $\ULB_\infty$,
and the quotient prop,
$$
\Po := \ULB_\infty/I^\circlearrowright,
$$
is a {\em dg}\, wheeled prop whose representations in a dg vector space $V$
are ine one-to-one correspondence with formal unimodular Poisson structures,
$\pi\in \wedge^\bullet\cT_V$, which vanish at $O\in V$.

\mip

{\bf 5.2.2. Remark.} Every free wheeled prop has a natural filtration by the number of vertices.
For applications to homological algebra and differential geometry one is often interested
in {\em completed}\, (with respect to this filtration)
{\em topological}\,  props,
 and in {\em continuous}\, morphisms between them \cite{Me-lec,Me-LieB,MV}).

\sip

In the next section we shall assume that both dg props $\ULB_\infty$ and $\Po$ are
completed with respect to the filtration by the number of vertices.

\mip

{\bf 5.3. Quasi-isomorphism theorem}. {\em
A continuous  morphism of dg wheeled topological props,
$$
F: \ULB_\infty \lon \Po,
$$
given on the generators by the formula
$$
F\left( \begin{xy}
 <0mm,0mm>*{\mbox{$\xy *=<3mm,3mm>
\txt{{{$_a$}}}*\frm{-}\endxy$}};
<-1.5mm,1.5mm>*{};<-5mm,5mm>*{}**@{-},
<-1mm,1.5mm>*{};<-2.7mm,5mm>*{}**@{-},
<1.5mm,1.5mm>*{};<5mm,5mm>*{}**@{-},
<1mm,1.5mm>*{};<2.7mm,5mm>*{}**@{-},
<0mm,4.5mm>*{...};
<-1.5mm,-1.5mm>*{};<-5mm,-5mm>*{}**@{-},
<-1mm,-1.5mm>*{};<-2.7mm,-5mm>*{}**@{-},
<1.5mm,-1.5mm>*{};<5mm,-5mm>*{}**@{-},
<1mm,-1.5mm>*{};<2.7mm,-5mm>*{}**@{-},
<0mm,-4.5mm>*{...};
<0mm,0mm>*{};<-5.8mm,6.3mm>*{^1}**@{},
   <0mm,0mm>*{};<-3mm,6.3mm>*{^2}**@{},
   <0mm,0mm>*{};<2.6mm,6.3mm>*{^{m\hspace{-0.5mm}-\hspace{-0.5mm}1}}**@{},
   <0mm,0mm>*{};<7.0mm,6.3mm>*{^m}**@{},
<0mm,0mm>*{};<-5.8mm,-6.3mm>*{_1}**@{},
   <0mm,0mm>*{};<-3mm,-6.3mm>*{_2}**@{},
   <0mm,0mm>*{};<3mm,-6.3mm>*{_{n\hspace{-0.5mm}-\hspace{-0.5mm}1}}**@{},
   <0mm,0mm>*{};<7.0mm,-6.3mm>*{_n}**@{},
 \end{xy}
 \right)=\left\{\Ba{cr}
 \begin{xy}
 <0mm,0mm>*{\circ};<0mm,0mm>*{}**@{},
 <-0.6mm,0.44mm>*{};<-8mm,5mm>*{}**@{-},
 <-0.4mm,0.7mm>*{};<-4.5mm,5mm>*{}**@{-},
 <0mm,0mm>*{};<-1mm,5mm>*{\ldots}**@{},
 <0.4mm,0.7mm>*{};<4.5mm,5mm>*{}**@{-},
 <0.6mm,0.44mm>*{};<8mm,5mm>*{}**@{-},
   <0mm,0mm>*{};<-8.5mm,5.5mm>*{^1}**@{},
   <0mm,0mm>*{};<-5mm,5.5mm>*{^2}**@{},
   <0mm,0mm>*{};<4.5mm,5.5mm>*{^{m\hspace{-0.5mm}-\hspace{-0.5mm}1}}**@{},
   <0mm,0mm>*{};<9.0mm,5.5mm>*{^m}**@{},
 <-0.6mm,-0.44mm>*{};<-8mm,-5mm>*{}**@{-},
 <-0.4mm,-0.7mm>*{};<-4.5mm,-5mm>*{}**@{-},
 <0mm,0mm>*{};<-1mm,-5mm>*{\ldots}**@{},
 <0.4mm,-0.7mm>*{};<4.5mm,-5mm>*{}**@{-},
 <0.6mm,-0.44mm>*{};<8mm,-5mm>*{}**@{-},
   <0mm,0mm>*{};<-8.5mm,-6.9mm>*{^1}**@{},
   <0mm,0mm>*{};<-5mm,-6.9mm>*{^2}**@{},
   <0mm,0mm>*{};<4.5mm,-6.9mm>*{^{n\hspace{-0.5mm}-\hspace{-0.5mm}1}}**@{},
   <0mm,0mm>*{};<9.0mm,-6.9mm>*{^n}**@{},
 \end{xy} & \mbox{for}\ a=0,\\
 0 & \mbox{otherwise},
 \Ea
 \right.
$$
is a quasi-isomorphism.}

\mip

\Proof The prop $\ULB_\infty$ is generated by the $\bS$-module
(\ref{Generators in ULieB_infty}). Let us enlarge the latter non-differential $\bS$-bimodule   to a {\em dg}\,
$\bS$-bimodule, $(E=\{E(m,n)\}, d_0)$, given by
\Beqrn
E(m,n) &:=& \bigoplus_{a\geq 0}^\infty\left(sgn_m\ot \id_n[m-2-2a]\oplus sgn_m\ot \id_n[m-1-2a]
\right)\\
&=&
\mbox{span}\left\langle
\begin{xy}
 <0mm,0mm>*{\mbox{$\xy *=<3mm,3mm>
\txt{{{$_a$}}}*\frm{-}\endxy$}};
<-1.5mm,1.5mm>*{};<-5mm,5mm>*{}**@{-},
<-1mm,1.5mm>*{};<-2.7mm,5mm>*{}**@{-},
<1.5mm,1.5mm>*{};<5mm,5mm>*{}**@{-},
<1mm,1.5mm>*{};<2.7mm,5mm>*{}**@{-},
<0mm,4.5mm>*{...};
<-1.5mm,-1.5mm>*{};<-5mm,-5mm>*{}**@{-},
<-1mm,-1.5mm>*{};<-2.7mm,-5mm>*{}**@{-},
<1.5mm,-1.5mm>*{};<5mm,-5mm>*{}**@{-},
<1mm,-1.5mm>*{};<2.7mm,-5mm>*{}**@{-},
<0mm,-4.5mm>*{...};
<0mm,0mm>*{};<-5.8mm,6.3mm>*{^1}**@{},
   <0mm,0mm>*{};<-3mm,6.3mm>*{^2}**@{},
   <0mm,0mm>*{};<2.6mm,6.3mm>*{^{m\hspace{-0.5mm}-\hspace{-0.5mm}1}}**@{},
   <0mm,0mm>*{};<7.0mm,6.3mm>*{^m}**@{},
<0mm,0mm>*{};<-5.8mm,-6.3mm>*{_1}**@{},
   <0mm,0mm>*{};<-3mm,-6.3mm>*{_2}**@{},
   <0mm,0mm>*{};<3mm,-6.3mm>*{_{n\hspace{-0.5mm}-\hspace{-0.5mm}1}}**@{},
   <0mm,0mm>*{};<7.0mm,-6.3mm>*{_n}**@{},
 \end{xy}\ ,
 \
 \begin{xy}
 <0mm,0mm>*{\mbox{$\xy *=<3mm,3mm>
\txt{{{$_a$}}}*\frm{-}\endxy$}};
<-1.5mm,1.5mm>*{};<-5mm,5mm>*{}**@{-},
<-1mm,1.5mm>*{};<-2.7mm,5mm>*{}**@{-},
<1.5mm,1.5mm>*{};<4mm,4mm>*{}**@{-},
<1mm,1.5mm>*{};<2.7mm,5mm>*{}**@{-},
<0mm,4.5mm>*{...};
<-1.5mm,-1.5mm>*{};<-5mm,-5mm>*{}**@{-},
<-1mm,-1.5mm>*{};<-2.7mm,-5mm>*{}**@{-},
<1.5mm,-1.5mm>*{};<4mm,-4mm>*{}**@{-},
<1mm,-1.5mm>*{};<2.7mm,-5mm>*{}**@{-},
<0mm,-4.5mm>*{...};
<0mm,0mm>*{};<-5.8mm,6.3mm>*{^1}**@{},
   <0mm,0mm>*{};<-3mm,6.3mm>*{^2}**@{},
   <0mm,0mm>*{};<3.0mm,6.3mm>*{^m}**@{},
<0mm,0mm>*{};<-5.8mm,-6.3mm>*{_1}**@{},
   <0mm,0mm>*{};<-3mm,-6.3mm>*{_2}**@{},
   <0mm,0mm>*{};<3.0mm,-6.3mm>*{_n}**@{},
   \ar@{->}@(ur,dr) (4.0,4.0)*{};(4.0,-4.0)*{},
 \end{xy}
\right\rangle
\Eeqrn

with the direct summand zero unless $m+n+2a\geq 3$, $m+a\geq 2$ and $n+a\geq 2$,
and with differential $d_0$ given on the generators of $E$
 by
$$
d_0
\begin{xy}
 <0mm,0mm>*{\mbox{$\xy *=<3mm,3mm>
\txt{{{$_a$}}}*\frm{-}\endxy$}};
<-1.5mm,1.5mm>*{};<-5mm,5mm>*{}**@{-},
<-1mm,1.5mm>*{};<-2.7mm,5mm>*{}**@{-},
<1.5mm,1.5mm>*{};<5mm,5mm>*{}**@{-},
<1mm,1.5mm>*{};<2.7mm,5mm>*{}**@{-},
<0mm,4.5mm>*{...};
<-1.5mm,-1.5mm>*{};<-5mm,-5mm>*{}**@{-},
<-1mm,-1.5mm>*{};<-2.7mm,-5mm>*{}**@{-},
<1.5mm,-1.5mm>*{};<5mm,-5mm>*{}**@{-},
<1mm,-1.5mm>*{};<2.7mm,-5mm>*{}**@{-},
<0mm,-4.5mm>*{...};
<0mm,0mm>*{};<-5.8mm,6.3mm>*{^1}**@{},
   <0mm,0mm>*{};<-3mm,6.3mm>*{^2}**@{},
   <0mm,0mm>*{};<2.6mm,6.3mm>*{^{m\hspace{-0.5mm}-\hspace{-0.5mm}1}}**@{},
   <0mm,0mm>*{};<7.0mm,6.3mm>*{^m}**@{},
<0mm,0mm>*{};<-5.8mm,-6.3mm>*{_1}**@{},
   <0mm,0mm>*{};<-3mm,-6.3mm>*{_2}**@{},
   <0mm,0mm>*{};<3mm,-6.3mm>*{_{n\hspace{-0.5mm}-\hspace{-0.5mm}1}}**@{},
   <0mm,0mm>*{};<7.0mm,-6.3mm>*{_n}**@{},
 \end{xy}
 =(-1)^{m-1}
\begin{xy}
 <0mm,0mm>*{\mbox{$\xy *=<6mm,3mm>
\txt{{{$_{a-1}$}}}*\frm{-}\endxy$}};
<-1.5mm,1.5mm>*{};<-5mm,5mm>*{}**@{-},
<-1mm,1.5mm>*{};<-2.7mm,5mm>*{}**@{-},
<1.5mm,1.5mm>*{};<4mm,4mm>*{}**@{-},
<1mm,1.5mm>*{};<2.7mm,5mm>*{}**@{-},
<0mm,4.5mm>*{...};
<-1.5mm,-1.5mm>*{};<-5mm,-5mm>*{}**@{-},
<-1mm,-1.5mm>*{};<-2.7mm,-5mm>*{}**@{-},
<1.5mm,-1.5mm>*{};<4mm,-4mm>*{}**@{-},
<1mm,-1.5mm>*{};<2.7mm,-5mm>*{}**@{-},
<0mm,-4.5mm>*{...};
<0mm,0mm>*{};<-5.8mm,6.3mm>*{^1}**@{},
   <0mm,0mm>*{};<-3mm,6.3mm>*{^2}**@{},
   <0mm,0mm>*{};<3.0mm,6.3mm>*{^m}**@{},
<0mm,0mm>*{};<-5.8mm,-6.3mm>*{_1}**@{},
   <0mm,0mm>*{};<-3mm,-6.3mm>*{_2}**@{},
   <0mm,0mm>*{};<3.0mm,-6.3mm>*{_n}**@{},
   \ar@{->}@(ur,dr) (4.0,4.0)*{};(4.0,-4.0)*{},
 \end{xy}
$$
$$
d_0
 \begin{xy}
 <0mm,0mm>*{\mbox{$\xy *=<3mm,3mm>
\txt{{{$_a$}}}*\frm{-}\endxy$}};
<-1.5mm,1.5mm>*{};<-5mm,5mm>*{}**@{-},
<-1mm,1.5mm>*{};<-2.7mm,5mm>*{}**@{-},
<1.5mm,1.5mm>*{};<4mm,4mm>*{}**@{-},
<1mm,1.5mm>*{};<2.7mm,5mm>*{}**@{-},
<0mm,4.5mm>*{...};
<-1.5mm,-1.5mm>*{};<-5mm,-5mm>*{}**@{-},
<-1mm,-1.5mm>*{};<-2.7mm,-5mm>*{}**@{-},
<1.5mm,-1.5mm>*{};<4mm,-4mm>*{}**@{-},
<1mm,-1.5mm>*{};<2.7mm,-5mm>*{}**@{-},
<0mm,-4.5mm>*{...};
<0mm,0mm>*{};<-5.8mm,6.3mm>*{^1}**@{},
   <0mm,0mm>*{};<-3mm,6.3mm>*{^2}**@{},
   <0mm,0mm>*{};<3.0mm,6.3mm>*{^m}**@{},
<0mm,0mm>*{};<-5.8mm,-6.3mm>*{_1}**@{},
   <0mm,0mm>*{};<-3mm,-6.3mm>*{_2}**@{},
   <0mm,0mm>*{};<3.0mm,-6.3mm>*{_n}**@{},
   \ar@{->}@(ur,dr) (4.0,4.0)*{};(4.0,-4.0)*{},
 \end{xy}
 =0.
$$
It is clear that the cohomology, $H(E)=\{H(E)(m,n)\}$, of this dg $\bS$-bimodule is equal to
$$
H(E)(m,n)=\mbox{span}\left\langle
\begin{xy}
 <0mm,0mm>*{\mbox{$\xy *=<3mm,3mm>
\txt{{{$_0$}}}*\frm{-}\endxy$}};
<-1.5mm,1.5mm>*{};<-5mm,5mm>*{}**@{-},
<-1mm,1.5mm>*{};<-2.7mm,5mm>*{}**@{-},
<1.5mm,1.5mm>*{};<5mm,5mm>*{}**@{-},
<1mm,1.5mm>*{};<2.7mm,5mm>*{}**@{-},
<0mm,4.5mm>*{...};
<-1.5mm,-1.5mm>*{};<-5mm,-5mm>*{}**@{-},
<-1mm,-1.5mm>*{};<-2.7mm,-5mm>*{}**@{-},
<1.5mm,-1.5mm>*{};<5mm,-5mm>*{}**@{-},
<1mm,-1.5mm>*{};<2.7mm,-5mm>*{}**@{-},
<0mm,-4.5mm>*{...};
<0mm,0mm>*{};<-5.8mm,6.3mm>*{^1}**@{},
   <0mm,0mm>*{};<-3mm,6.3mm>*{^2}**@{},
   <0mm,0mm>*{};<2.6mm,6.3mm>*{^{m\hspace{-0.5mm}-\hspace{-0.5mm}1}}**@{},
   <0mm,0mm>*{};<7.0mm,6.3mm>*{^m}**@{},
<0mm,0mm>*{};<-5.8mm,-6.3mm>*{_1}**@{},
   <0mm,0mm>*{};<-3mm,-6.3mm>*{_2}**@{},
   <0mm,0mm>*{};<3mm,-6.3mm>*{_{n\hspace{-0.5mm}-\hspace{-0.5mm}1}}**@{},
   <0mm,0mm>*{};<7.0mm,-6.3mm>*{_n}**@{},
 \end{xy}\right\rangle.
$$

Consider next the decreasing filtrations
$$
\ULB_\infty=F_0\ULB_\infty \supset F_1 \ULB_\infty \supset \ldots \supset F_p\ULB_\infty\supset\ldots
$$
$$
\Po=F_0\Po \supset F_1 \Po \supset \ldots \supset F_p\Po\supset\ldots
$$
 of dg props $\ULB_\infty$ and $\Po$ by the number
of vertices: the subspaces $F_p$ spanned, by definition, by decorated graphs with at least $p$ vertices.
The morphism $F$ respects the filtrations and hence induces morphism,
$\{F_r: (\cE_r \ULB_\infty, d_r) \rar (\cE_r\Po, \delta_r\}$, of the associated spectral sequence,
in particular, a morphism,
$$
F_0: (\cE_0 \ULB_\infty, d_0) \rar (\cE_r\Po, d_0\}
$$
of the initial terms. The dg $\bS$-bimodule $(\cE_0 \ULB_\infty, d_0)$
is canonically isomorphic to the following one,
$$
\cF^\circlearrowright_{no\ loops}\langle E\rangle
:=\sum_{G\in \fG^\circlearrowright_{no\ loops}} G\langle E
\rangle,
$$
with the differential induced from  $d_0$ on $E$ (hence the same notation).
As we work over a field of characteristic zero, by Kunneth and Mashke theorems
the functor $\cF^\circlearrowright_{no\ loops}$ on the category
of dg $\bS$-bimodules is exact, i.e.\
$$
H\left( \cF^\circlearrowright_{no\ loops}\langle E\rangle\right)
= \cF^\circlearrowright_{no\ loops}\langle H\left(E\right)\rangle.
$$
Therefore, the morphism $F_0$ is an isomorphism. By assumptions on $\ULB_\infty$ and
$\Po$, both filtrations are complete, exaustive and regular (degenerating at 1st term).
Hence the associated spectral sequences are
convergent by classical Complete Convergence Theorem 5.5.10 (see p.139 in \cite{Weib}).
Then, by the classical Comparison Theorem 5.2.12 (see p.\ 126  \cite{Weib}),
the morphism $F$ is a quasi-isomorphism. \hfill $\Box$

\bip
\section{$BF$  theory of quantum BV manifolds}

{\bf 6.1. Introduction.}
This section is inspired by the work of Mnev \cite{Mn} on a remarkable approach
to the homotopy transfer formulae of unimodular
$L_\infty$-algebras which is based on  the BV quantization of an extended $BF$ theory
and the associated  Feynman integrals. We apply here
 Losev-Mnev ideas
to unimodular Lie 1-bialgebras
and show that the Feynman integrals technique provides us with exactly the same
formulae for the homotopy transfer of $\ULB_\infty$-structures as the ones which one obtains
with the help of the Koszul duality technique in the wheeled props  approach to quantum BV manifolds
(see \S\S 2-4).
We believe that the established interrelation,
$$
\mbox{\sf Feynman integrals} \leftrightarrows \mbox{\sf
Morphisms of dg wheeled (co)props}
$$
is quite general.

\mip

{\bf 6.2. $BF$-theory of unimodular Lie 1-bialgebras.} Let $V$ be a finite-dimensional, and assume that its dual space $V^*$ is equipped with a structure of
 unimodular dg Lie 1-bialgebra, i.e.\ with
a degree 1 Lie brackets
 $[\ \bullet\ ]: \odot^2 V^*\rar V^*[1]$ and a degree 0 Lie co-brackets $\Delta^{\CoLie}: V^*\rar \wedge^2 V^*$
(see \S 4.2).
The dualization and degree shifting of the latter gives a map
$[\ ,\ ]: \odot^2(V[-1])\rar V[-2]$ which makes $V[-1]$ into a degree 1 Lie algebra.
Consider a  degree 2 polynomial function (called {\em action})
  on the vector space $V^*\oplus V[-1]$,
$$
\Ba{rccc}
S:& V^*\oplus V[-1]  & \lon & \K\vspace{2mm}\\
& p\oplus \om & \lon & S(p,\om):= \langle p, d\om \rangle
 \ + \ \frac{1}{2}\langle p, [\om , \om]\rangle +
\frac{1}{2}\langle [p \bullet p], \om\rangle,
\Ea
$$
where $\langle\ ,\ \rangle$ stand for the natural pairing.
A choice of a basis $\{e_a\}$ in $V$ induces linear coordinates $\{p_a: |p_a|=|e_a|\}$
on $V^*$ and linear coordinates $\{\om^a: |\om_a|= 1- |e_a|\}$ on $V[-1]$ in
which the function $S$ takes the form
$$
S(p,\om)= \sum_{a,b}\left(p_a D_b^a \om^b \pm \sum_c\frac{1}{2}
\left(p_b p_cC_a^{bc}\om_a \pm
p_a\Phi^a_{bc}\om^b\om^c\right)\right),
$$
where $ D_b^a$, $C^{bc}_a$ and $\Phi^a_{bc}$ are the structure constants of, respectively, the differential,  the odd Lie brackets
and  Lie cobrackets in the chosen basis.

\mip

Let $\cM_{V}^\hbar$ be an odd symplectic manifold corresponding to the completed
graded commutative ring $\widehat{\odot}(V\oplus V^*[1])[[\hbar]]\simeq \K[[p_a, \om^a, \hbar]]$.
\mip

{\bf 6.2.1. Lemma.} {\em
 The semidensity $e^{\frac{S(p,\om)}{\hbar}} \sqrt{D}_{p,\om}$ makes $\cM_{V}^\hbar$ into a quantum
BV manifold.}

\mip

\Proof The boundary consitions $S|_{p=0}=0$ and $S|_{\om=0}=0$ are obvious so that,
by defintion 3.9.2(i), one should only check the equation $\hbar \Delta_0 S + \frac{1}{2}\{S\bullet S\}=0$,
where $\Delta_0=\sum\frac{\p^2}{\p p_a\p \om^a}$. As $S$ is independent of $\hbar$, this is equivalent to
two equations,
$$
 \{S\bullet S\}=0\ \ \mbox{and}\ \ \ \Delta_0S=0.
$$
The first equation  follows from relations (\ref{R for LieB}).
Equations $\Delta_0\langle p, [\om \bullet \om]\rangle=0$
and $\Delta_0 \langle [p,p], \om\rangle=0$
are equivalent to unimodularity of $[\ \bullet\ ]$ and $\delta^{\CoLie}$.
Finally,  equation $\Delta_0 \langle p, d\om \rangle=0$
follows from the well-known fact that,
 for an arbitrary differential $d$, there exist
a basis in $V$ in which $d$ is given by a matrix (\ref{Differential matrix}) with zero
supertrace.
\hfill $\Box$

\mip

The quadratic form $S_{(2)}:=\langle p, d\om\rangle$ is degenerate on the vector space $V^*\oplus V[-1]$.
We shall next specify a subspace, $W \subset V^*\oplus V[-1]$, on which $S_{(2)}$
 is non-degenerate so that one can develop
a perturbative quantization of the action $S=S_{(2)}+S_{(3)}$ with $S_{(2)}$ determining the ``propagator"
of the quantum theory and
with the cubic part, $S_{(3)}:=\frac{1}{2}\langle p, [\om \bullet \om]\rangle +
\frac{1}{2}\langle [p,p], \om\rangle$,
 playing the role of ``interactions" between ``fields" $p$ and $\om$.
With this purpose
we fix an arbitrary cohomological splitting,
\Beq\label{Cohomological splitting}
V= H(V) \oplus B \oplus B[-1],
\Eeq
of the complex $V$. Let $p_a=\{p'_\ta, p''_\al, p'''_\al\}$ be an adopted to this
splitting basis of $V$ in which the differential is given
by the matrix (\ref{Differential matrix}). Put another way,  $\{p'_\ta\}_{\ta\in I'}$ is a basis
of the cohomology group $H(V,d)$, $\{p''_\al\}_{\al\in J}$  a basis of $B$,  $\{p'''_\al\}_{\al\in J}$ a basis
of $B[-1]$ and the differential $d$ is given by
$$
dp'_\ta=0, \ \ \ dp''_\al=p'''_\al,\ \ \ dp'''_\al=0.
$$
This splitting of $V$ induces associated splitting of $V^*[1]$ and hence the associated split base
of the direct sum $V\oplus V^*[1]$ which we denote as follows,
$$
\underbrace{V}_p\ \oplus\ \underbrace{V^*[1]}_\om
= \underbrace{H(V)}_{p'_\ta}\
\ \oplus \ \underbrace{B}_{p''_\al}\ \oplus\ \underbrace{B[-1]}_{p'''_\al}\ \ \oplus \ \
\underbrace{H(V)^*[1]}_{\om'^\ta}\ \oplus
\underbrace{B^*[1]}_{\om'''^\al}\ \oplus\ \underbrace{B^*[2]}_{\om''^\al}.
$$
so that
$$
d\om'^\ta=0, \ \ \ d\om''^\al=-\om'''^\al,\ \ \ d\om'''^\al=0.
$$
The linear functions on the space $V^*\oplus V[-1]$ corresponding the above basis vectors of $V\oplus V^*[1]$
we denote by the same letters $p'_\ta, p''_\al, p'''_\al, \om'^\ta, \om''^\al, \om'''^\al$.

Then the quadratic term of the action takes the form
(cf.\ (\ref{Gamma_2_chapter3}))
\Beq\label{S_2}
S_{(2)}=\langle p, d\om\rangle= -< p''', \om''' >=-\sum_{\al\in J}
p'''_\al \om'''^\al.
\Eeq
where $<\ ,\ >$ is the natural degree 2 paring between $B$ and $B^*[2]$.

\mip

Let now $\cM^\hbar_{B\oplus B[-1]}$ be the formal odd symplectic manifold corresponding to a
graded commutative algebra
$$
\widehat{\odot^\bullet}\left(B\oplus B[-1]\oplus B^*[1]\oplus B^*[2]\right)\ot \K[[\hbar]]
\simeq \K[[p'',p''', \om'', \om''', \hbar]],
$$
and
$\cM^\hbar_{H(V)}$ the odd symplectic manifold corresponding to
$$
\widehat{\odot^\bullet}\left(H(V)\oplus H(V)^*[1]\right)[[\hbar]]
\simeq \K[[p', \om', \hbar]].
$$
Cohomological splitting (\ref{Cohomological splitting}) induces
an isomorphism of odd Poisson manifolds,
$$
\cM^\hbar_{V} = \cM^\hbar_{H(V)} \times \cM^\hbar_{B\oplus B[-1]}.
$$
Following \cite{Mn} we shall show next how a perturbative Feynman type integration along a Lagrangian submanifold $\caL$
in the odd symplectic manifold $\cM^\hbar_{B\oplus B[-1]}$ transforms a  simple quantum BV structure
on $\cM^\hbar_{V}$ given by Lemma~6.2.1 into a rather non-trivial quantum BV structure on $\cM^\hbar_{H(V)}$ (in a full
accordance with Theorem 2.7.1). Let $\sqrt{D}_{B\oplus B[-1]}$  be the
semidensity on  $\cM^\hbar_{B\oplus B[-1])}$ associated with the choice of linear Darboux coordinates
made above.

\mip

{\bf 6.2.2. Lemma.}
{\em For any Lagrangian submanifold $\caL$ in $\cM^\hbar_{B\oplus B[-1]}$ and any function $f\in
\f_{\cM_V}$ one has,
$$
\bar{\Delta}_0 \int_{\caL}   f  \sqrt{D}_{B\oplus B[-1]}|_\caL =
\int_{\caL}   (\Delta_0 f)  \sqrt{D}_{B\oplus B[-1]}|_\caL,
$$
provided the integral exists. Here
 $\Delta_0=\sum_a\frac{\p^2}{\p p_a\p\om^a}$ is the odd Laplacian on $\cM^\hbar_{V}$,
$\bar{\Delta}_0 =\sum_\ta\frac{\p^2}{\p p_\ta\p\om^\ta}$ is the odd Laplacian on $\cM^\hbar_{H(V)}$ and
$\sqrt{D}_{B\oplus B[-1]}|_\caL$ stands for the restriction (in accordance with
(\ref{Semidensity_on_Lagrangian}))
of the semidensity
$\sqrt{D}_{B\oplus B[-1]}$ to a volume form on $\caL$.}

\mip

This Lemma is in fact a classical Stokes theorem in disguise. We refer to \cite{Sc} or \cite{CF2}
for its simple proof. Thus, if we can find a Lagrangian submanifold $\caL\subset
\cM_{B\oplus B[-1]}$ such that the integral   $\int_{\caL}   f  \sqrt{D}_{B\oplus B[-1]}|_\caL$
exists  for  $f=e^{\frac{S(p,\om)}{\hbar}}$  given by Lemma~6.2.1,
then we obtain a quantum BV structure on the $\hbar$-twisted
 odd symplectic manifold $\cM^\hbar_{H(V)}$
from the unimodular Lie 1-bialgebra structure on $V^*$. Formula (\ref{S_2}) suggests a natural
choice: let $\caL$ be the formal $\Z$-graded manifold associated with the vector subspace
$B[-1]\oplus B^*[1]\subset B\oplus B[-1]\oplus B^*[1]\oplus B^*[2]$. It is a submanifold
of $\cM_{B\oplus B[-1]}$ given by the equations $p''=\om''=0$. The semidensity
$\sqrt{D}_{B\oplus B[-1]}$ restricts to $\caL$ as an ordinary translation invariant
Berezin volume $dp'''d\om'''=\prod_{\al}dp'''_\al d\om'''^\al$ (see \cite{Be}).
As the  quadratic volume form $S_{(2)}=-<p''',\om'''>$ is obviously
non-degenerate on $B[-1]\oplus B^*[1]$, the integral,
$$
N:= \int_\caL e^{\frac{S_2(p,\om)}{\hbar}}\sqrt{D}_{B\oplus B[-1]}|_\caL =\int
e^{-\frac{<p''',\om'''>}{\hbar}} dp''' d\om''',
$$
is a well-defined constant\footnote{This is a ``Gaussian" integral of special type 1.2.1.2
according to Cattaneo's review \cite{Ca2} of Gaussian integrals. Strictly speaking, we should
 view here the formal parameter $\hbar$ as a purely imaginary complex number $ih$ with $h$ being
 an arbitrary positive real number; such ``Gaussian"
 integrals can be made well-defined  via a real analytic
 continuation of ordinary Gauss integrals for positive definitive quadratic forms, see \cite{Ca2}.}.
Moreover,
\Beqrn
e^{\frac{S_{\eff}(p',\om',\hbar)}{\hbar}} &:=& N^{-1}
 {\int e^{\frac{S(p',p''',\om',\om''')}{\hbar}} dp''' d\om'''}
= N^{-1}\int e^{\frac{-<p''',\om'''> + S_{(3)}(p',p''',
\om',\om''')}{\hbar}} dp''' d\om'''\\
 &=& N^{-1}\sum_{k\geq 0}\frac{\hbar^{-k}}{k!}
\int e^{-\frac{<p''',\om'''> }{\hbar}}\left(S_{(3)}(p',p''',
\om',\om''') \right)^k   dp''' d\om'''
\Eeqrn
is well-defined as an element of the algebra (\ref{h-ring}).
It can be computed
via the classical Vick theorem (see. e.g., \cite{Ca2}, ) with the propagator
$\langle\langle \om''', p'''\rangle\rangle$ (which is, by definition is equal to
the  quadratic form inverse to $S_{(2)})$ given by the matrix\footnote{
this matrix (up to the factor $\hbar^{-1}$) is precisely the coordinate representation
of the homotopy operator $h:V\rar V$ (see \S 2.7).}
$$
\langle\langle \om'''^\al, p'''_\be\rangle\rangle_0:= -\hbar\delta_\be^\al.
$$
 As
\Beqrn
\displaystyle
S_{(3)}(p',p''', \om',\om''')&=&
\frac{1}{2}\left\langle p'+p''', \left[(\om'+\om''') \bullet (\om'+\om''')\right]\right\rangle +
\frac{1}{2}\left\langle [p'+p''',p'+p'''], \om'+\om'''\right\rangle\\
&=& S_{3}(p',\om') + \langle  p',[\om'\bullet \om'''] \rangle +
\frac{1}{2}\left\langle p', \left[\om''' \bullet \om'''\right]\right\rangle
+ \langle  p''',[\om'\bullet \om'''] \rangle\\
&&  +
\frac{1}{2}\left\langle p''', \left[\om''' \bullet \om'''\right]\right\rangle
+ \langle  [p',p'''],\om'\rangle + \langle  [p',p'''],\om'''\rangle +
\frac{1}{2}\left\langle [p''',p'''], \om'\right\rangle\\
&& + \frac{1}{2}\left\langle [p''',p'''], \om''\right\rangle,
\Eeqrn
we conclude by the Wick theorem that this  integral  is equal to the formal power series,
$$
e^{\frac{S_{\eff}(p',\om',\hbar)}{\hbar}} =
\sum_{G\in \tilde{G}^\circlearrowright} G(p',\om', \hbar)
$$
where the sum runs over all possible  graphs
built from corollas of two types,
$$
[\ ,\ ] \leftrightarrow
 \begin{xy}
 <0mm,-0.55mm>*{};<0mm,-2.5mm>*{}**@{-},
 <0.5mm,0.5mm>*{};<2.2mm,2.2mm>*{}**@{-},
 <-0.48mm,0.48mm>*{};<-2.2mm,2.2mm>*{}**@{-},
 <0mm,0mm>*{\circ};<0mm,0mm>*{}**@{},
 \end{xy}
 \ \ \  , \ \ \
[\, \bullet \, ] \leftrightarrow
 \begin{xy}
 <0mm,0.66mm>*{};<0mm,3mm>*{}**@{-},
 <0.39mm,-0.39mm>*{};<2.2mm,-2.2mm>*{}**@{-},
 <-0.35mm,-0.35mm>*{};<-2.2mm,-2.2mm>*{}**@{-},
 <0mm,0mm>*{\circ};<0mm,0mm>*{}**@{},
 \end{xy}.
$$
It is well-known (see, e.g., Ch.\ 4, \S 3 in \cite{Ma2} or Proposition 2.10 in \cite{Po})
that
$$
\log  \sum_{G\in \tilde{G}^\circlearrowright} G(p',\om',\hbar)=
\sum_{G\in \tilde{G}_c^\circlearrowright}
 G(p',\om',\hbar)
 $$
where the sum  on the r.h.s.\ runs over the subset, $\tilde{G}_c^\circlearrowright
\subset \tilde{G}^\circlearrowright$, consisting of {\em connected}\, graphs.
 Thus the effective action can be written finally as
\Beq\label{S-eff}
 S_{\eff}= \sum_{G\in \tilde{G}_c^\circlearrowright}
 G(p',\om',\hbar)=\sum_{g\geq 0} \sum_{G\in \tilde{G}_{g,c}^\circlearrowright}
 \hbar^g G(p',\om')
\Eeq
where
\Bi
\item[--] the second sum runs over the subset,
  $\tilde{G}_{g,c}^\circlearrowright\subset \tilde{G}_{g,c}^\circlearrowright$,
  consisting  of all possible connected
trivalent directed graphs of genus $g$;
\item[--] $G(p',\om')$ is a linear map
$H(V)^{\ot \bullet}\rar H(V)^{\ot \bullet}$
obtained from the graph
$G$ by decorating it exactly as in Theorem 2.7.1:
vertices are decorated by the structure constants, $C_{ab}^c$ and $\Phi_{a}^{bc}$,
of the Lie and co-Lie operations in $V$, and internal edges are decorated
 with the homotopy operator
$h$; legs  are now decorated with $p''$ and $\om''$.
\Ei
By Lemmas 6.2.2 and 6.2.1, the effective action satisfies the equation,
$$
\Delta_0 e^{\frac{S_{\eff}(p',\om',\hbar)}{\hbar}}=0, \ \ \mbox{i.e.}\ \
\hbar \Delta_0 S_\eff + \frac{1}{2} \{S_\eff\bullet S_\eff\}=0,
$$
and hence makes $\cM_{H(V)}$ into a quantum BV manifold.

\mip

{\bf 6.2.3. Proposition.} {\em For any dg Lie 1-bialgebra on $V$ and any cohomological splitting of $V$
there is a canonically
 associated structure of quantum BV manifold on the cohomology, $H(V)$, given by the quantum master
function (\ref{S-eff}). Moreover, there exists a natural quasi-isomorphism of quantum BV manifolds,}
$$
\phi_\hbar: \left(\cM_{H(V)}, e^{\frac{S_\eff(p',\om',\hbar)}{\hbar}}\sqrt{D}_{p',\om'}\right) \lon
\left(\cM_V, e^{\frac{S(p,\om)}{\hbar}}\sqrt{D}_{p,\om}\right).
$$

\sip

\Proof  It remains to construct a morphism $\phi_\hbar$, which, by definition~3.9.11,
is a topological morphism
of $\K[[\hbar]]$-modules,
$$
\phi^*: \K[[p, \om,\hbar]].
\lon  \K[[p', \om',\hbar]].
$$
which in the limit $\hbar\rar 0$ induces a morphism of algebras and satisfies
the equation
\Beq\label{Morhism_BV_BF/theory}
{e^{\frac{-S_{\eff}(p',\om',\hbar)}{\hbar}}}{\bar{\Delta}_0
\left(\phi^*_\hbar(f)e^{\frac{S_{\eff}(p',\om',\hbar)}{\hbar}}\right)}
=
\phi^*\left(e^{\frac{-S(p,\om)}{\hbar}}{\Delta}_0\left(fe^{\frac{S(p,\om)}{\hbar}}\right)\right),
\Eeq
for any $f\in \K[[p, \om,\hbar]]$. In view of Lemma~6.2.2, the map (cf.\ \cite{Mn})
\Beq\label{quantum_embeding_formula}
\phi^*_\hbar(f):=N^{-1} e^\frac{-S_{\eff}(p',\om',\hbar)}{\hbar}{\underset{p''=0, \om''=0}{\int} f(p, \om, \hbar) e^{\frac{S(p,\om)}{\hbar}} dp''' d\om'''}
\Eeq
does satisfy equation (\ref{Morhism_BV_BF/theory}):
\Beqrn
{e^{\frac{-S_{\eff}(p',\om',\hbar)}{\hbar}}}{{\Delta}'_0
\left(\phi^*_\hbar(f)e^{\frac{S_{\eff}(p',\om',\hbar)}{\hbar}}\right)} &=&
N^{-1}e^{\frac{-S_{\eff}(p',\om',\hbar)}{\hbar}}
\underset{p''=0, \om''=0}{\int}\Delta_0\left( f(p, \om, \hbar) e^{\frac{S(p,\om)}{\hbar}}\right) dp''' d\om'''\\
&=&\phi^*\left(e^{\frac{-S(p,\om)}{\hbar}}{\Delta}_0\left(fe^{\frac{S(p,\om)}{\hbar}}\right)\right).
\Eeqrn
Moreover, in the limit $\hbar\rar 0$  formula (\ref{Morhism_BV_BF/theory}) gives simply the evaluation map,
$$
\lim_{\hbar\rar 0} \phi^*_\hbar(f)= f|_{\hbar=0, p''=0,p'''=0,\om''=0,\om'''=0},
$$
 and hence defines a morphism of algebras $\f_{\cM_V^0}\rar \f_{\cM_{H(V)}^0}$.
\hfill $\Box$
\mip

Formula (\ref{Morhism_BV_BF/theory}) proves Proposition 3.11.2 in the special case when the quantum
master function $S(\om,p)$ is associated with a unimodular Lie 1-bialgebra structure on a finite-dimensional vector space. However the same formula (\ref{quantum_embeding_formula}) gives obviously a well-defined perturbative power series in $\hbar$
for an arbitrary (quasi-classically) split quantum master function $S(p,\om,\hbar)$ and proves thereby Proposition 3.11.2
in general.

\bip

{\em Acknowledgement}. It is a pleasure to thank Alberto Cattaneo and Anton Khoroshkin
for valuable discussions and comments. Thanks go also to the anonymous referee for a careful reading of the paper, useful suggestions and a list of misprints.
This work was partially supported by the
G\"oran Gustafsson foundation.
\def\cprime{$'$}

  \end{document}